\documentclass{article}
\usepackage{epsfig}

\usepackage{tikz}
\usepackage{mathtools}
\usepackage{upgreek}

\usepackage{amsmath,amssymb,amscd,amsthm,bm} 
\usepackage{mathdots} 

\usepackage{caption}
\usepackage{subcaption}
 
\title{Heron triangles with two rational medians and Somos-5 sequences}  
\author{Andrew N.W. Hone~\\
School of Mathematics, 
Statistics \& Actuarial Science~\\
University of
Kent~\\ Canterbury CT2 7FS, UK.
}

\newcommand{\beq}{\begin{equation}}  
\newcommand{\eeq}{\end{equation}}  
\newcommand{\bea}{\begin{eqnarray}} 
\newcommand{\eea}{\end{eqnarray}}   
\newcommand{\bear}{\begin{array}}  
\newcommand{\eear}{\end{array}}

\newtheorem{thm}{Theorem}[section] 

\newtheorem{propn}[thm]{Proposition}

\newtheorem{lem}[thm]{Lemma}

\newtheorem{cor}[thm]{Corollary} 

\newenvironment{prf}{\trivlist \item [\hskip 
\labelsep {\bf Proof:}]\ignorespaces}{\qed \endtrivlist}

\theoremstyle{definition}
\newtheorem{definition}[thm]{Definition}
\newtheorem{remark}[thm]{Remark}

\newcommand{\F}{{\mathbb F}}

\newcommand{\Q}{{\mathbb Q}}
\newcommand{\Z}{{\mathbb Z}}
\newcommand{\C}{{\mathbb C}}
\newcommand{\R}{{\mathbb R}}
\newcommand{\Proj}{{\mathbb P}}

\newcommand{\ra}{\mathrm{a}}

\newcommand{\rx}{\mathrm{x}}
\newcommand{\ry}{\mathrm{y}}

\newcommand{\ri}{\mathrm{i}}

\newcommand\la{{\lambda}}

\newcommand\ka{{\kappa}}
\newcommand\al{{\alpha}}
\newcommand\be{{\beta}}
\newcommand\ze{{\zeta}}
\newcommand\gam{{\gamma}}

\newcommand\om{{\omega}}

\newcommand\tal{\tilde{\alpha}}
\newcommand\tbe{\tilde{\beta}}
\newcommand\tgam{\tilde{\gamma}}
\newcommand\tmu{\tilde{\mu}}
\newcommand\si{{\sigma}}

\newcommand\eps{{\epsilon}}

\newcommand{\cQ}{{\cal Q}}
\newcommand{\cP}{{\cal P}}
\newcommand{\cU}{{\cal U}}
\newcommand{\cV}{{\cal V}}
\newcommand{\cL}{{\cal L}}

\newcommand\bs{\bar{s}}
\newcommand\ba{\bar{a}}
\newcommand\bb{\bar{b}}
\newcommand\bc{\bar{c}}
\newcommand\bk{\bar{k}}
\newcommand\bl{\bar{\ell}}

\def\lc{\left\lfloor}   
\def\rc{\right\rfloor}

\begin{document} 

\maketitle
 
\begin{abstract} 
Triangles with integer length sides and integer area are known as Heron triangles. Taking rescaling freedom into account,  
one can apply the same name when all sides and the area are rational numbers. 
A perfect triangle is 
a Heron triangle with all three medians being rational, and it is a longstanding conjecture that no such triangle exists. 
However, Buchholz and Rathbun showed that there are  
infinitely many Heron triangles with two rational medians, an infinite subset of which are associated with rational points on 
an elliptic curve $E(\Q)$ with Mordell-Weil group  $\Z\times \Z/2\Z$, and they observed a connection with 
a pair of Somos-5 sequences. Here we make the latter connection more precise by providing explicit formulae 
for the integer side lengths, the two rational medians,  and the area in this infinite family of Heron triangles. 
The proof uses a combined approach to Somos-5 sequences and associated Quispel-Roberts-Thompson 
(QRT) maps in the plane, from several different viewpoints: complex analysis, real dynamics, and reduction modulo a prime.
\end{abstract}

\section{Introduction} 

\setcounter{equation}{0}

The formula 
\beq\label{heron} 
\Delta=\sqrt{s(s-a)(s-b)(s-c)}
\eeq
for the area of a triangle with sides $(a,b,c)$, where 
$$ 
s=\frac{a+b+c}{2}
$$ 
is the semiperimeter, 
is attributed to Heron of Alexandria. If $(a,b,c)$ is a triple of positive integers and the area $\Delta$ is also an integer, then this is called a Heron triangle. 
More generally, due to rescaling freedom, we say that a triangle is Heron whenever the side lengths and the area are all rational numbers. A method 
for enumerating Heron triangles was given by Schubert in \cite{schubert}, 
but a parametric formula equivalent to  
\beq\label{brahma} 
a=\frac{p^2+r^2}{p}, \qquad 
b=\frac{q^2+r^2}{q}, \qquad 
c=\pm\frac{(r^2-pq)(p+q)}{pq},  
\eeq 
for $p,q,r\in\Q$ 
and area $\Delta=rc\in\Q$ was already known to Brahmagupta in the 7th century A.D. \cite{dickson}.  
Any Pythagorean triple gives a right-angled Heron triangle, while the triangle with 
integer side lengths  $(5,5,6)$ and area $12$ arising from the choice $p=q=1, r=2$ in Brahmagupta's formula is the simplest isosceles Heron 
triangle (in the sense of having the smallest value of $a+b+c$), 
and the simplest example of a Heron triangle that is neither right-angled nor isosceles has side lengths $(15,13,14)$ and area 
$84$, being obtained by taking $p=3,q=4,r=6$ in the same formula. 
There are numerous Diophantine problems concerning Heron triangles, many of which are related to the theory 
of elliptic curves \cite{dujella, goins, heron1, heron2}. 

It is an old problem to answer the question as to whether there exists a perfect triangle: one with integer sides, medians, and area; or equivalently, is there a Heron triangle 
with three rational medians? The expectation is that there is no such triangle, but 
to prove it seems very difficult, and it is remarked in \cite{guy} that 
despite incorrect ``proofs'' in the literature, the problem remains open.  One of the first incorrect arguments is implicit in Schubert's work \cite{schubert}, 
where he claimed to present a complete parametrization of Heron triangles with one of the medians being  rational, and used this to argue that Heron 
triangles with two rational medians are impossible. However, his proposed parametrization was incomplete, and Schubert's oversight was 
pointed out 
by Dickson \cite{dickson} and  in the PhD thesis of Buchholz \cite{bthesis},  who initially 
found the case $(73,51,26)$ with area $420$ and two rational medians, of lengths 
$\frac{35}{2}$ and $\frac{97}{2}$ respectively, as well as a small number of other examples - see Table \ref{tab:table1}, in which each triangle is represented (up to scale) 
by an integer triple with $\gcd (a,b,c)=1$.

\begin{figure}
\centering
\label{triangle}
\begin{tikzpicture}
\draw (2,0) 
  -- (7,0) 
-- (5.25,1) node[above=1pt] {$b$}
--(3.5,2) 
--(2.75,1) node[left=1pt] {$c$}
-- (2,0); 
\draw (3.5,2)  
--(4.1,0.8) node[right=1pt] {$k$} 
-- (4.5,0) ;
\draw (3.3,1.7) arc [start angle=240, end angle=285, radius=5mm] node[anchor=north east] {$\be$};
\draw (3.65,1.7) arc [start angle=270, end angle=294, radius=7mm] node[below=1pt] {$\al$};
\draw (4.4,0.24) arc [start angle=95, end angle=170, radius=3mm]  node[anchor=north west] {$a$} node[anchor=south east] {$\gam $};
\end{tikzpicture}

\caption{Triangle with one labelled median}
\end{figure}
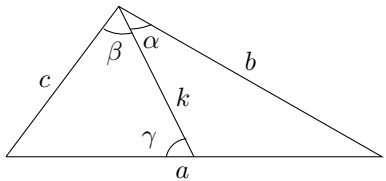

Henceforth we denote the medians that bisect sides $a,b,c$ by $k,\ell,m$, respectively, so that 
\beq\label{kmed} 
k^2=\frac{1}{4}(2b^2+2c^2-a^2), \quad  
\ell^2=\frac{1}{4}(2c^2+2a^2-b^2), \quad 
m^2=\frac{1}{4}(2a^2+2b^2-c^2), 
\eeq 
and label the angles adjacent to the median $k$ as in Fig.1. 
Then the area of the 
triangle satisfies 
$
\Delta = bk\sin\al =ck\sin\be =\frac{1}{2}ak\sin\gam, 
$
and, following \cite{schubert}, it is helpful to consider the 
half-angle cotangents 
\beq\label{mpx} 
M=\mathrm{cot}(\al/2), \quad 
P=\mathrm{cot}(\be/2),
\quad 
X=\mathrm{cot}(\gam/2),
\eeq 
which we will refer to as the Schubert parameters, using the same nomenclature and notation as in \cite{br1}.
Up to rescaling, these three parameters completely determine the triangle; clearly they are not independent, but as shown 
by Schubert they satisfy 
the  equation 
\beq\label{schub}
M-\frac{1}{M}=P-\frac{1}{P}+
2\left( X-\frac{1}{X}\right). 
\eeq  
Upon rewriting the latter as 
$
2MP(X^2-1)+MX(P^2-1)-PX(M^2-1)=0$,
we see that this 
defines an affine quartic surface in three dimensions, which we call  the Schubert surface. 
The Schubert parameters are given in terms of the area, side lengths and the median $k$ by 
the formulae 
\beq\label{schubb} 
M=\frac{4\Delta}{4bk+a^2-3b^2-c^2}, \quad 
P=\frac{4\Delta}{4ck+a^2-b^2-3c^2}, \quad 
X=\frac{4\Delta}{2ak-b^2+c^2}, 
\eeq  
which follow from the half-angle identity $\mathrm{cot}(\al/2)=\sin\al/(1-\cos\al)$ and the cosine rule, 
while the ratios of side lengths are given in terms of the  Schubert parameters by 
\beq\label{aboverc} 
\frac{a}{c}=\frac{2(X+X^{-1})}{P+P^{-1}}, \qquad
\frac{b}{c}=\frac{M+M^{-1}}{P+P^{-1}}.
\eeq

\begin{table}[h!]
  \begin{center}
    \caption{The smallest Heron triangles with two rational medians.}
    \label{tab:table1}
    \begin{tabular}{ | r|| r| r| r| r | r| r|} %
\hline
      $n$ & $a$  &  $b$ & $c$ & $k$ & $\ell$ & $\Delta$ \\
      \hline \hline
      1 & 73 & 51 & 26 & ${35}/{2}$ & ${97}/{2}$ & 420 \\
\hline
      2 & 626 &  875 & 291 & 572 &${433}/{2}$ & 55440 \\
\hline
      *    & 1241 & 4368 & 3673  & ${7975}/{2}$ & 1657 & 2042040 \\
\hline 
    **   &14384 & 14791 & 11257  & 11001 & 21177/2 & 75698280 \\
\hline
3 & 28779 & 13816 &  15155 & 3589/2 & 21937 & 23931600 \\ 
\hline
4 & 1823675& 185629& 1930456& 2048523/2& 3751059/2& 142334216640 \\
\hline
*** & 2288232 & 1976471 & 2025361 & 1641725 & 3843143/2 & 1877686881840 \\ 
\hline
**** & 22816608 & 20565641 & 19227017 & 
            16314487 & 36845705/2 & 185643608470320  \\ 
\hline 
5& 2442655864& 2396426547& 46263061& 1175099279& 2488886435/2& 2137147184560080 \\
\hline 
    \end{tabular}
  \end{center}
\end{table}

In view of the formulae (\ref{schubb}), for a Heron triangle with rational median $k$ the corresponding Schubert parameters are rational, and conversely, if 
$(M,P,X)\in\Q^3$ is a rational point on the Schubert surface (\ref{schub}), then the triangle is Heron with (at least) one rational median. Strictly speaking, we require positive rational 
solutions, since the half-angle cotangents must be positive, but the surface  (\ref{schub}) admits the obvious involutions $M\to -M^{-1}$, $P\to-P^{-1}$, $X\to-X^{-1}$, 
as well as $(M,P,X)\to(M^{-1},P^{-1},X^{-1})$, 
so 
if one of the coordinates is negative  it can be replaced by minus its reciprocal, while 
all three coordinates can simultaneously be replaced by their reciprocals, 
and we shall exploit this freedom in what follows.  The inherent subtlety in the problem of characterizing 
Heron triangles with one rational median, to which Schubert gave an incomplete solution, can be seen from the fact that the Schubert surface admits three different elliptic fibrations, obtained 
by fixing the value of any one of the parameters. For instance, setting $M=M_0\in\Proj^1(\Q)$ gives the cubic curve $2P(X^2-1)+X(P^2-1)+CPX=0$, where the constant $C=1/M_0-M_0$; so 
for generic values of $C$ the fibre is an elliptic curve, with j-invariant $(C^4 + 40C^2 + 208)^3/\big(16(C^4 + 40C^2 + 144)\big)$, and each (finite) element of the group of rational points 
on the curve 
corresponds to a Heron triangle with one rational median.
 
For a Heron triangle with two rational medians $k,\ell$, there are two associated rational points on the Schubert surface, namely the point $(M_a,P_a,X_a)$ associated with the median $k$ bisecting side $a$, as given by the formulae (\ref{schubb}), 
and the point $(M_b,P_b,X_b)$ associated with the median $\ell$ bisecting side $b$, given by the same formulae but replacing $a\to b$, $b\to c$, $c\to a$ and $k\to\ell$. As well as satisfying the equation (\ref{schub}), these two sets of 
Schubert parameters must be related by the compatibility conditions 
\beq\label{aovercab1} 
\frac{2(X_a+X_a^{-1})}{P_a+P_a^{-1}}= 
\frac{P_b+P_b^{-1}}{M_b+M_b^{-1}},  
\qquad 
\frac{M_a+M_a^{-1}}{P_a+P_a^{-1}} = 
\frac{2(X_b+X_b^{-1})}{M_b+M_b^{-1}},  
\eeq 
corresponding to the ratios of the side lengths as in (\ref{aboverc}). Thus the problem of finding a Heron triangle with two rational medians is equivalent to finding a pair of positive rational points on the Schubert surface (\ref{schub}), subject 
to the pair of constraints (\ref{aovercab1}). 
The angles $\al,\be$ as in Fig.1 must also satisfy $\al+\be<\pi$, so this imposes the additional requirements 
\beq\label{extra}
\mathrm{arccot}(M_a)+\mathrm{arccot}(P_a)<\frac{\pi}{2}, \qquad 
\mathrm{arccot}(M_b)+\mathrm{arccot}(P_b)<\frac{\pi}{2}, 
\eeq 
but once a pair of compatible positive triples has been found, these requirements can always be 
satisfied by applying $(M_a,P_a,X_a)\to( M_a^{-1},P_a^{-1},X_a^{-1})$ and/or 
$(M_b,P_b,X_b)\to( M_b^{-1},P_b^{-1},X_b^{-1})$ if necessary, since these transformations leave the 
constraints   (\ref{aovercab1}) invariant. 

There is another approach to the problem, based on the formulae 
\beq\label{abcparam} 
\begin{array}{rcl} 
a & = &\uptau \,(-2\theta^2\phi -\theta\phi^2+2\theta\phi-\phi^2+\theta+1), \\
b & = & \uptau\,(\theta^2\phi+2\theta\phi^2-\theta^2+2\theta\phi-\phi+1), \\
c & = & \uptau\, (\theta^2\phi-\theta\phi^2+\theta^2+2\theta\phi+\phi^2+\theta-\phi),
\end{array} 
\eeq 
found by  Buchholz \cite{bthesis}, which provide a rational parametrization of triangles with two rational medians 
$k,\ell$, where $\theta,\phi$ are rational numbers (constrained suitably to ensure positivity), and the parameter 
$\uptau\in\Q$ allows for the arbitrary choice of scale. Conversely, $\theta,\phi\in\Q$ 
can be written as functions of the (ratios of) side lengths, given by 
\beq\label{thetaphi}
\theta = \frac{c-a\pm2\ell}{2s}, \qquad 
\phi   = \frac{b-c\pm2k}{2s},
\eeq 
with $s=\frac{1}{2}(a+b+c)$ being the semiperimeter, as before. 

Hence an efficient method to search for 
Heron triangles with two rational medians is to run through the rational parameters $\theta,\phi$, ordered by height, 
and check whether the corresponding value of $\Delta$ is a rational number. 
More precisely, given $\theta=R/S\in\Q$ written as a fraction in lowest terms, its naive height is $H(\theta)=\max (|R|,|S|)$, 
and pairs $(\theta,\phi)\in\Q^2$  can be enumerated in order of increasing height $\tilde{H}=\max\big(H(\theta),H(\phi)\big)$, so fixing the scale 
$\uptau=1$ in (\ref{abcparam}), the side lengths $(a,b,c)$ of triangles  with two rational medians can be calculated from 
this parametrization for each pair of parameters with $\tilde{H}=1,2,3,\ldots$ and then it can be checked from 
Heron's formula (\ref{heron}) whether $s(s-a)(s-b)(s-c)$ is a perfect square, corresponding to the area being rational. 
(This method leads to duplicate triangles related to one another by different values of the scaling $\uptau$, but still 
seems more efficient than finding Heron triangles with one rational median and then checking whether a second median is rational.) 
  
The latter  method was implemented by Buchholz 
and Rathbun, initially working independently  (an independent search was also carried out by Kemnitz), 
yielding the first six rows in Table \ref{tab:table1}. In \cite{br1} they observed  remarkable properties of 
certain triangles in the latter table, with respect to their Schubert parameters, which are shown in Table \ref{tab:table2}: 
for the rows labelled by an integer $n=1,2,3,\ldots$, the factorizations are related, and in particular the parameter $M_b$ 
in row $n$ is minus the reciprocal of the parameter $P_a$ in row $n+1$ (see Table \ref{tab:table3} for details of the factorizations). The triangles labelled by asterisks do not 
seem to fit into any obvious pattern, but their observations on the other triangles (corresponding to $n=1,\ldots,5$ in Tables 1 \& 2)  
 led them to suggest that  these examples 
should extend 
to an  infinite family of triangles labelled by a positive integer $n$, 
with a  conjectured factorization of the Schubert parameters as   
\beq\label{schubs5a} 
M_a= -\frac{S_{n+1}S_{n+2}^2T_n}{S_nT_{n+1}T_{n+2}^2}, \quad 
P_a=-\frac{S_{n+1}S_{n+2}T_{n+1}T_{n+2}}{S_nS_{n+3}T_nT_{n+3}}, \quad 
X_a=2^{(-1)^{n+1}}\frac{S_nS_{n+2}^2T_{n+3}}{S_{n+3}T_nT_{n+2}^2}, 
\eeq 
\beq\label{schubs5b} 
M_b=\frac{S_{n+1}S_{n+4}T_{n+1}T_{n+4}}{S_{n+2}S_{n+3}T_{n+2}T_{n+3}}, \quad  
P_b=-\frac{S_{n+2}^2S_{n+3}T_{n+4}}{S_{n+4}T_{n+2}^2T_{n+3}}, \quad 
X_b=2^{(-1)^{n}}\frac{S_{n+1}T_{n+2}^2T_{n+4}}{S_{n+2}^2S_{n+4}T_{n+1}}, 
\eeq 
where $(S_n)$ and $(T_n)$ are integer sequences given by 
\beq\label{s5seq} 
(S_n): \quad1,1,1,2,3,5,11,37,83,274,\ldots , 
\eeq 
and 
\beq\label{s5eds} 
(T_n): \quad0,1,-1,1,1,-7,8,-1,-57,391,\ldots , 
\eeq 
(the terms above are listed starting from the index $n=0$). 
These are Somos sequences, of the kind introduced in \cite{somos}. More specifically, they are both 
Somos-5 sequences: $(S_n)$ is generated by the fifth order   
quadratic recurrence 
\beq\label{s5orig} 
S_{n+5}S_n=S_{n+4}S_{n+1}+ S_{n+3}S_{n+2}  
\eeq 
(see 
\cite{oeis}); the sequence is usually generated starting 
from five initial 1s, but here we have indexed it so that it extends symmetrically to negative $n$, with $S_{-n}=S_n$. 
As we shall see, the  sequence $(T_n)$ is closely related to $(S_n)$: it is generated by the same fifth order recurrence 
(\ref{s5orig}), and extends backwards in an antisymmetric fashion, so that  $T_{-n}=-T_n$; it is also a divisibility 
sequence, having the property that $T_n|T_m$ whenever $n|m$. (It is almost an 
elliptic divisibility sequence  in the sense of \cite{ward}, but the terms 
with even/odd index satisfy different relations of order four.) 
Henceforth we shall  refer to the sequence of triangles corresponding to the pairs of Schubert parameters (\ref{schubs5a}) and (\ref{schubs5b}) 
as the main sequence. 


\begin{table}[h!]
  \begin{center}
    \caption{Schubert parameters (with signs) for the smallest Heron triangles with two rational medians.}
    \label{tab:table2}
\scalebox{0.9}{
    \begin{tabular}{ | r|| r| r| r| r | r| r|} %
\hline
      $n$ & $M_a$  &  $P_a$ & $X_a$ & $M_b$ & $P_b$ & $X_b$ \\
\hline 
\hline   
0 & 0 & $\infty$ & $\infty$ & -3/2 & -2/3 & 2/3 \\
\hline 
    1 &  4  & 2/3 &  8/3 & 35/6&  84/5 & 7/40 \\
\hline 
      2 & 18 & -6/35 & 63/10 
& -176/105 & 360/77 &  32/99
 \\
\hline 
      *   & 
728/51 &  17 &    48/91  &  231/260 & 2431/420 &  17/55
\\
\hline  
   **   &  1395/476 &  620/153 & 63/85 &  357/95 &  4845/1736 &  1767/1360
 \\
\hline 
3 & -75/98 & 105/176 & 800/539 & 111/3080 &  275/14504 &  -147/1850 \\ 
\hline 
4 & 605/1344& -3080/111 &-363/4736 & -165585/3256 & -255189/5312 & 36480/70301 \\ \hline 
*** & 7144/2277 &    79101/24472 & 7238/7429 & 394128/101365 &49742 /11155 & 24035/27936
 \\ 
\hline
**** &  1035096/312455 & 1542840/505571 & 770431/717145 & 
 770431/218064 &  505571/117691 & 337421/412896 
\\ 
\hline
5&  105413/40 & 3256/165585 & 780330/581 & 9427792/175047 & 44428157/15618 & 4301/6001696 \\
\hline 
    \end{tabular}
}
  \end{center}
\end{table}

\begin{table}[h!]
  \begin{center}
    \caption{Prime factors of the first few Schubert parameters in the main sequence.}
    \label{tab:table3}
    \begin{tabular}{ | r|| r| r| r| r | r| r|} %
\hline
      $n$ & $M_a$  &  $P_a$ & $X_a$ & $M_b$ & $P_b$ & $X_b$ \\
\hline 
&&&&&& \\ 
0 & 0 & $\infty$ & $\infty$ & $-\tfrac{3}{2}$& $-\tfrac{2}{3}$ & $\tfrac{2}{3}$ \\
&&&&&& \\ 
    1 &  $2^2$ & $\tfrac{2}{3}$ &  $\tfrac{2^3}{3}$ & $\tfrac{5\cdot 7}{2\cdot 3}$ &  $\tfrac{2^2\cdot 3\cdot 7}{5}$ & $\tfrac{7}{2^3\cdot 5} $\\
&&&&&& \\ 
      2 & $2\cdot 3^2$ & $-\tfrac{2\cdot 3}{5\cdot 7}$ & $\tfrac{3^2\cdot 7}{2\cdot 5}$  
& $-\tfrac{2^4\cdot 11}{3\cdot5\cdot7}$ & $\tfrac{2^3\cdot3^2\cdot 5}{7\cdot11}$ & $\tfrac{2^5}{3^2\cdot11}  $
 \\
&&&&&& \\ 
	3 & $-\tfrac{3\cdot5^2}{2\cdot7^2}$ & $\tfrac{3\cdot5\cdot 7}{2^4\cdot 11}$ & $\tfrac{2^5\cdot5^2}{7^2\cdot 11}$ & 
$\tfrac{3\cdot 37}{2^3\cdot 5\cdot7\cdot 11}$ &  $\tfrac{5^2\cdot11}{2^3\cdot7^2\cdot37}$ & $-\tfrac{3\cdot7^2}{2\cdot5^2\cdot37}$ \\
&&&&&& \\ 
4 & $\tfrac{ 5\cdot11^2}{2^6\cdot3\cdot7}$ & $-\tfrac{2^3\cdot5\cdot7\cdot11}{3\cdot37}$ & $-\tfrac{3\cdot11^2}{2^7\cdot 37}$ & 
$-\tfrac{3\cdot5\cdot7\cdot19\cdot83}{2^3\cdot11\cdot37}$ & $-\tfrac{3\cdot11^2\cdot19\cdot 37 }{2^6\cdot 83 }$ &  $\tfrac{2^7\cdot3\cdot5\cdot 19}{ 7\cdot11^2\cdot 83 } $
 \\  
&&&&&& \\ 
5& $\tfrac{7\cdot11\cdot37^2}{2^3\cdot 5} $ & $ \tfrac{2^3\cdot 11 \cdot 37}{3\cdot5\cdot7 \cdot 19 \cdot 83}$  & $\tfrac {2\cdot3\cdot 5\cdot 19 \cdot 37^2}{ 7\cdot83 }$ &  
$\tfrac{ 2^4\cdot 11 \cdot 17 \cdot 23\cdot137 }{ 3\cdot19\cdot 37\cdot 83 }$ & $\tfrac{ 17\cdot23\cdot37^2\cdot   83}{ 2\cdot3\cdot 19 \cdot 137 }$ & $\tfrac{ 11\cdot17\cdot 23}{2^5\cdot37^2\cdot  137}$ 
  \\
&&&&&& \\ 
\hline 
    \end{tabular}
  \end{center}
\end{table}

Despite being provided with a theta function formula for the Somos-5 sequence $(S_n)$ by Elkies \cite{elkies}, 
Buchholz and Rathbun were unable to use this 
to prove that the Schubert parameters for this proposed infinite family of Heron triangles with two rational medians 
are given by the factorizations (\ref{schubs5a}) and (\ref{schubs5b}). Nevertheless, they were able to make 
further progress by plotting the coordinates of the sequence of parameters $(\theta,\phi)$ given by 
(\ref{thetaphi}) (with both signs taken as $+$) 
  corresponding 
to these triangles, which were empirically found to lie on one of five birationally equivalent curves ${\cal C}_{1-5}$ of genus one, repeating 
in the pattern ${\cal C}_2,{\cal C}_1,{\cal C}_4,{\cal C}_3,{\cal C}_2,{\cal C}_1,{\cal C}_5 $ with period 7, the simplest such curve being the biquadratic cubic 
\beq\label{curvec}
{\cal C}_4: \qquad 
\theta^2\phi-\theta\phi^2+\theta\phi+2\theta-2\phi-1=0. 
\eeq  
Over $\Q$, this is birationally equivalent to the elliptic curve 
\beq\label{ellipticc} 
y^2+xy=x^3+x^2-2x,
\eeq  
which has Mordell-Weil group $\Z\times \Z/2\Z$, the same curve corresponding to the theta function formula for 
the Somos-5 sequence (\ref{s5seq}) found by Elkies \cite{elkies}. 
In a subsequent paper \cite{br2}, Buchholz and Rathbun proved the 
following result. 

\begin{thm} \label{brthm}
Every rational point $(\theta,\phi)$ on the genus one curve ${\cal C}_4$ given by (\ref{curvec}), with $0<\theta<1$, 
$0<\phi<1$, $2\theta+\phi>1$  corresponds 
to a Heron triangle with two rational medians. 
\end{thm}

In subsequent work \cite{symm}, they considered the full set of discrete symmetries of the problem in terms of the 
parameters $a,b,c,k,\ell$, including sign changes e.g.\ 
$a\to -a$, $b\to -b$, etc.\,, as well as allowed permutations, such as the reflection symmetry $a\leftrightarrow b$, $k\leftrightarrow \ell$
(equivalent to changing the orientation of the triangle), and showed that, under the action of this group on the pairs  $(\theta,\phi)$, they 
obtained points on a total of eight 
isomorphic curves ${\cal C}_{1-8}$ corresponding to triangles in the main sequence; yet the four sporadic triangles, labelled with asterisks in 
Tables 1 and 2, do not give points on these curves, and we do not know if there are formulae analogous to 
(\ref{schubs5a}) and (\ref{schubs5b}) for these sporadic cases. 
More recent work on this problem has consisted of proving that all of the Heron triangles in the main sequence, 
corresponding to rational points on one of  these eight curves, have exactly two rational medians, so none of them are perfect 
triangles \cite{br3, ismail1, ismail2}. 
However, until now, many of Buchholz and Rathbun's original observations about this sequence have lacked an explanation. 

\begin{table}[h!]
  \begin{center}
    \caption{Prime factors of the semiperimeter, reduced side lengths  and area in the main sequence.}
    \label{tab:table4}
\scalebox{0.85}{
    \begin{tabular}{ | r|| r| r| r| r | r|} %
\hline
      $n$ & $s$  &  $s-a$ & $s-b$ & $s-c$  & $\Delta$ \\
\hline 
    1 &  $3\cdot 5^2 $& 2 & $2^3\cdot 3$ & $7^2$  & $2^2\cdot3\cdot 5\cdot 7$ \\ 
      2 & $5\cdot 11^2$ & $3\cdot 7$ & $2\cdot 3^3\cdot 5$ &  $2^7\cdot 7$
& $2^4\cdot3^2\cdot 5\cdot 7\cdot 11$
 \\
	3 & $11\cdot 37^2$ & $2^3\cdot 5\cdot 7^3$& $3\cdot 5^3\cdot 7\cdot  11$ & $ 2^5\cdot 3$
& $2^4\cdot3\cdot 5^2\cdot 7^2\cdot 11\cdot 37$
\\
4 & $ 7\cdot 37 \cdot 83^2$ & $ 2^9\cdot 7 \cdot 11$ & $2^3\cdot 5 \cdot 11^3\cdot 37$ &  $3^4\cdot 5 \cdot 19^2$
& $2^6\cdot3^2\cdot 5\cdot 7\cdot 11^2\cdot19\cdot  37\cdot 83$
 \\  
5& $ 2^5\cdot 7^2\cdot 83\cdot 137^2$ & $2^3\cdot 3 \cdot 19 \cdot 37$ &  $11\cdot 37^3\cdot 83$ & $3\cdot 5^2 \cdot 11 \cdot 17^2 \cdot 19\cdot 23^2$
& $2^4\cdot3\cdot 5\cdot 7\cdot 11\cdot 17\cdot19\cdot 23\cdot  37^2\cdot 83\cdot 137$
  \\
\hline 
    \end{tabular}
}
  \end{center}
\end{table}

In considering this problem afresh, we observed an elegant factorization pattern for the semiperimeter $s$, the quantities $s-a$, $s-b$, $s-c$, which we refer to as the 
reduced lengths, and hence also for the area $\Delta$ of the triangles in the main sequence (see Table \ref{tab:table4}), and we found that they could be written in terms 
of the two Somos-5 sequences. This led us not only to a proof of the formulae (\ref{schubs5a}) and (\ref{schubs5b}) for the Schubert parameters, but also to explicit 
expressions for the lengths of the sides, the two rational medians, and the area, as well as an explanation for the period 7 cycles of curves in the $(\theta,\phi)$ plane. 
Our main result is the following 

\begin{thm}\label{main} 
For each integer $n\geq 1$, the terms in the pair of Somos-5 sequences  (\ref{s5seq}) and (\ref{s5eds}) 
provide a Heron triangle with two rational medians, having integer side lengths given 
by  
\beq\label{lengthsn} 
\bear{rcl} 
a & = & |S_{n+1}S_{n+2}^3S_{n+3}T_{n+2}+S_n^2S_{n+1}T_{n+3}T_{n+4}^2|, \\
b & = &  
| S_n^2S_{n+1}T_{n+3}T_{n+4}^2-T_{n+1}T_{n+2}^3T_{n+3}S_{n+2}|, \\ 
c & = &  
|T_{n+1}T_{n+2}^3T_{n+3}S_{n+2}-S_{n+1}S_{n+2}^3S_{n+3}T_{n+2}|, 
\eear 
\eeq 
with $\gcd(a,b,c)=1$, 
rational 
median lengths 
\beq\label{lengthkln} 
\bear{rcl} 
k & = & {\scriptstyle \frac{1}{2} }
|S_{n+4}T_{n+4}(T_{n}T_{n+1}^2T_{n+2}-S_{n}S_{n+1}^2S_{n+2})|, \\
\ell & = &   {\scriptstyle \frac{1}{2} }
|S_nT_n(T_{n+2}T_{n+3}^2T_{n+4}-S_{n+2}S_{n+3}^2S_{n+4})|, 
\eear 
\eeq 
and area
\beq\label{arean} \Delta = 
|S_nS_{n+1}S_{n+2}^2S_{n+3}S_{n+4}T_nT_{n+1}T_{n+2}^2T_{n+3}T_{n+4}|
. 
\eeq
\end{thm} 

A brief outline of the paper is as follows. The next section is devoted to Somos-5 sequences and 
Quispel-Roberts-Thompson (QRT) maps: we rapidly review 
the necessary analytical fomulae from \cite{hones5}, in terms of Weierstrass functions, which are a key ingredient in our main argument,  and prove some 
determinantal identities connecting the sequences (\ref{s5seq}) and (\ref{s5eds}), before presenting  simple 
preliminary results on initial value problems and their reduction modulo a prime that will be needed later. 
We  then connect the two Somos-5 sequences with two different orbits 
of a QRT map in the plane, both of which lie on the same biquadratic curve that is isomorphic to  (\ref{curvec}), and with a 
single orbit of a QRT map on another curve related by a 2-isogeny. Section 3 contains the main results of the paper, leading to 
the proof of Theorem \ref{main}: the central result is Theorem   \ref{schubrels}, which is proved by writing 
the two sets of Schubert parameters (with signs) in terms of elliptic functions and using analytic 
arguments to verify that they lie on  
the Schubert surface as well as  satisfying the constraints (\ref{aovercab1}). However, in order to show that all the signs 
can be consistently removed by elementary transformations to end up with positive solutions of Schubert's equation, 
we need to consider the pattern of signs in the sequence (\ref{s5eds}), which turns out to have period 14, as a consequence of 
the real dynamics of one of the QRT orbits, which moves around certain segments of a curve with period 7 
(see Lemma \ref{signlem}); the latter 
pattern controls all the signs in the problem, and incidentally explains one of Buchholz and Rathbun's empirical 
observations on curves in the $(\theta,\phi)$ plane (Theorem \ref{7cycle}). The section ends with a complete description 
of the periodic dynamics of the QRT maps and associated Somos-5 sequences   over finite fields, combining and extending 
various results in the literature \cite{jogia, kanki, coprime, div, rob, swart}, which is required 
to analyse the common divisors of the side lengths. In section 4 we briefly discuss how geometrical arguments, namely  
Brahmagupta's construction, and a formula of Schubert for the tangents of half-angles in Heron triangles, lead to some 
additional identities between the Schubert parameters and other quantities involved. 
The final section contains our conclusions.

\section{Somos-5 sequences and QRT maps}

\setcounter{equation}{0}

Somos sequences are generated by quadratic recurrences of the form 
\beq\label{gensomos} 
\tau_{n+N}\tau_n=\sum_{j=1}^{\lc \frac{N}{2}\rc} \tilde{\al}_j\,\tau_{n+N-j}\tau_{n+j}, 
\eeq
where $\tal_j$ are coefficients. They encompass elliptic divisibility sequences in number theory, and as such can be regarded as nonlinear generalizations of Fibonacci, Lucas, or other linear recurrence sequences \cite{recs, ward}. 
If there are precisely two or three monomials on the right-hand side, then they are of the right shape to be generated from a  cluster algebra \cite{fz1} or an LP algebra \cite{lp}, providing one of the original 
examples of the Laurent phenomenon  \cite{fz,gale}. In addition, these special types of Somos recurrences can be obtained  as reductions of integrable partial difference equations on a three-dimensional lattice, namely the 
discrete Hirota equation \cite{hkw} or Miwa's equation \cite{fed}, which are also known by other names:  bilinear discrete KP/BKP, or the octahedron/cube recurrences. They also appear in the context 
of supersymmetric gauge theories and dimer models  \cite{bgm, efs, gk}. 

The rest of this section is devoted to presenting geometric, analytic, algebraic and arithmetic results about Somos-5 sequences, corresponding to the particular case $N=5$ of (\ref{gensomos}), as well as associated birational maps 
of the plane studied by Quispel, Roberts and Thompson (QRT maps).   

\subsection{Geometric, analytic and algebraic properties of Somos-5 sequences}

The general Somos-5 recurrence is  
\beq\label{s5gen} 
\tau_{n+5}\tau_n=\tal\,\tau_{n+4}\tau_{n+1}+\tbe\, \tau_{n+3}\tau_{n+2}. 
\eeq 
We take $\C$ as the ambient field, considering all sequences as complex-valued, but for suitable choices of the initial values and the coefficients $\tal,\tbe$ 
the recurrence produces integer sequences 
such as (\ref{s5seq}). 
One way to see this is to observe that the recurrence (\ref{s5gen}) has the Laurent property, as it arises 
from mutations in a cluster algebra \cite{fordy_marsh}, meaning that the 
iterates lie in the Laurent polynomial ring $\Z[\tal,\tbe,\tau_0^{\pm 1}, \tau_1^{\pm 1},\tau_2^{\pm 1},\tau_3^{\pm 1},\tau_4^{\pm 1}]$. Hence if 
all initial values are $\pm 1$ and the coefficients $\tal,\tbe$ are integers then $\tau_n\in\Z$ for all $n$. However, as shown in \cite{swahon}, due to the 
connection with the arithmetic of elliptic curves, a much stronger version of the Laurent property holds for this recurrence, and there are many more ways in which 
it can produce integer sequences. The geometrical structure of Somos-5 sequences is based on the following result. 
 
\begin{lem}\label{consquant} 
The recurrence (\ref{s5gen}) has two independent conserved quantities (first integrals), 
invariant 
under shifting $n\to n+1$, namely 
\beq\label{It} 
\tilde{I} = \frac{\tau_n\tau_{n+4}}{\tau_{n+1}\tau_{n+3}} +\tal\, \left(\frac{\tau_{n+1}^2}{\tau_{n}\tau_{n+2}}+
\frac{\tau_{n+2}^2}{\tau_{n+1}\tau_{n+3}}+\frac{\tau_{n+3}^2}{\tau_{n+2}\tau_{n+4}}\right) +\tbe\, \frac{\tau_{n+1}\tau_{n+3}}{\tau_n\tau_{n+4}} 
\eeq 
and 
\beq\label{Jt} 
\tilde{J} = 
\frac{\tau_n\tau_{n+3}}{\tau_{n+1}\tau_{n+2}} +\frac{\tau_{n+1}\tau_{n+4}}{\tau_{n+2}\tau_{n+3}} + 
\tal\, \left(\frac{\tau_{n+1}\tau_{n+2}}{\tau_n\tau_{n+3}}+ \frac{\tau_{n+2}\tau_{n+3}}{\tau_{n+1}\tau_{n+4}} \right) +\tbe\,
\frac{\tau_{n+2}^2}{\tau_{n}\tau_{n+4}}.
\eeq 
These two quantities are built from a 2-invariant, given by 
\beq\label{2invt}
\tilde{K}_n=\frac{\tau_n\tau_{n+4}+\tal\,\tau_{n+2}^2}{\tau_{n+1}\tau_{n+3}}, 
\eeq 
whose value repeats with period 2, so that $\tilde{K}_{n+2}=\tilde{K}_n$, with 
$$ 
\tilde{I}=\tilde{K}_n+\tilde{K}_{n+1}, \qquad 
\tal\tilde{J}+\tbe=\tilde{K}_n\tilde{K}_{n+1}. 
$$ 
\end{lem} 
\begin{prf} 
Applying (\ref{s5gen}) twice shows that $\tilde{K}_n$ satisfies 
$$
\tilde{K}_{n+1}=\tal\, \left(\frac{\tau_{n+1}^2}{\tau_{n}\tau_{n+2}}+
\frac{\tau_{n+3}^2}{\tau_{n+2}\tau_{n+4}}\right) + 
\tbe\, \frac{\tau_{n+1}\tau_{n+3}}{\tau_n\tau_{n+4}} , 
\qquad \tilde{K}_{n+2}=\tilde{K}_n,
$$ 
so the value of $\tilde{K}_n$ repeats with period 2, from which it is an immediate consequence that the sum 
$\tilde{K}_n+\tilde{K}_{n+1}=\tilde{I}$ and product  
$\tilde{K}_n\tilde{K}_{n+1}=\tal\tilde{J}+\tbe$
define two independent invariants $\tilde{I}, \tilde{J}$, which are given by (\ref{It}) and (\ref{Jt}), respectively. 
\end{prf} 

\begin{remark} \label{s4recs} 
By clearing the denominator in (\ref{2invt}), it follows that $\tau_n$ satisfies 
the 
Somos-4 relation 
$$ 
\tau_{n+4}\tau_{n} = \tilde{K}_n\, \tau_{n+3}\tau_{n+1} -\tal\,\tau_{n+2}^2, 
$$ 
with one of the coefficients  depending on the parity of $n$.
\end{remark} 

Geometrically, iteration of the Somos-5 recurrence  (\ref{s5gen})  is equivalent to iterating 
the birational map 
\beq\label{birats5} 
(\tau_0,\tau_1,\tau_2,\tau_3,\tau_4)\mapsto 
\left(\tau_1,\tau_2,\tau_3,\tau_4, \frac{\tal \tau_4\tau_1+\tbe \tau_3\tau_2}{\tau_0}\right)
\eeq 
in $\C^5$, and the existence of these two conserved quantities means that generic orbits lie on 
three-dimensional level sets given by fixed values of $\tilde{I},\tilde{J}$.  
However, the invariant $\tilde{I}$ will not play  a very significant role in what follows. The quantity 
$\tilde{J}$ is much more important, because it leads to the connection with elliptic curves: indeed, setting 
$$ 
U_n =\frac{\tau_{n+3}\tau_n}{\tau_{n+2}\tau_{n+1}}
$$ and comparing with (\ref{Jt}) shows that, for fixed $\tilde{J}$,  the pairs $(U_n,U_{n+1})$ in the plane lie 
on the biquadratic cubic curve defined by 
\beq\label{tJcurve} 
U_nU_{n+1}(U_n+U_{n+1}) +\tal(U_n+U_{n+1}) -\tilde{J}U_nU_{n+1}+\tbe =0, 
\eeq 
which (for generic values of $\tal,\tbe,\tilde{J}$) has genus one. The latter curve is birationally equivalent to  an 
elliptic curve in 
Weierstrass form (equation (\ref{origwei}) below), and this is what lies behind the analytic formula 
for the terms of a Somos-5 sequence obtained in \cite{hones5}, and described as follows.  

\begin{thm}\label{s5analytic} 
The general solution of 
(\ref{s5gen}) can be written in the form 
\beq\label{genivp}
\tau_n = A_\pm B_\pm^{\lc n/2\rc} \mu^{\lc n/2\rc^2}\si (z_0+n\kappa), 
\eeq 
where the $\pm$ subscripts apply for even/odd $n$, respectively, and $\si(z)=\si(z;g_2,g_3)$ is 
the Weierstrass sigma function associated with the elliptic curve 
\beq\label{origwei} 
y^2=4x^3-g_2x-g_3
\eeq 
with invariants defined by 
\beq\label{g2g3} 
g_2=12\tilde{\la}^2-2\tilde{J}, \qquad g_3=4 \tilde{\la}^3-g_2\tilde{\la}-\tmu^2
\eeq 
in terms of the quantities 
\beq\label{tlamu}
\tmu=(\tbe+\tal\tilde{J})^{\frac{1}{4}}, \qquad  
\tilde{\la} = \frac{1}{3\tmu^2} \left(\frac{\tilde{J}^2}{4}+\tal \right)
. \eeq
The solution corresponds to a sequence of points $\hat{{\cal P}}_0+n{\cal P}$ on the curve (\ref{origwei}), 
where  the initial point $\hat{{\cal P}}_0=\big(\wp(z_0),\wp'(z_0)\big)$ is arbitrary, and at each step it is 
translated by ${\cal P}= \big(\wp(\ka),\wp'(\ka)\big)=\big(\tilde{\la},\tmu\big)$. The other 
parameters appearing in (\ref{genivp}) are 
\beq\label{mu} 
\mu=\frac{\tmu}{\si(2\ka)} = -\si(\ka)^{-4}, 
\eeq 
and $A_+,A_-,B_+,B_-$ which are arbitrary up to the constraint that 
\beq\label{constr} 
\frac{B_+}{B_-}=-\mu^{-1}=\si(\ka)^{4}. 
\eeq 
\end{thm} 

\begin{remark} 
From the above result, the general solution of (\ref{s5gen})  
is given by fixing the 7 parameters $A_+,A_-,B_+,\ka,z_0,g_2,g_3$  (with $B_-$ given by (\ref{constr}) in terms of $B_+,\ka,g_2,g_3$), corresponding to the fact that 
the initial value problem for (\ref{s5gen})  is specified by a total of $7 = 5+2 $  parameters 
(five adjacent initial values, $\tau_0,\ldots,\tau_4$ say, plus the two coefficients $\tal, \tbe$). Moreover,  for generic initial values and 
coefficients, the initial value problem can be solved explicitly 
by using the relations 
(\ref{tlamu}) and (\ref{g2g3}) to obtain the curve (\ref{origwei}) from the values of $\tal,\tbe$ and the conserved quantity 
$\tilde{J}$ as in (\ref{Jt}); thereafter $\ka$ and $z_0$ are found by evaluating elliptic integrals, and $A_\pm,B_\pm$ can then be 
determined in terms of the initial values and values of the sigma function involving these arguments. We carry this out in detail below for the 
sequences (\ref{s5seq}) and (\ref{s5eds}).  
\end{remark} 

For what follows it is helpful to introduce another sequence $(\ra_n)$ associated with any solution 
(\ref{genivp}) of Somos-5, 
referred to as the companion EDS (elliptic divisibility sequence) in \cite{swahon}, which is 
defined by the analytic formula  
\beq\label{cEDS} 
\ra_n=\frac{\si(n\ka)}{\si(\ka)^{n^2}}. 
\eeq
The sequence of terms $\ra_n$ can be used to describe Somos relations of higher order satisfied by  $\tau_n$, 
which are summarized in the following way.  
\begin{thm}\label{compa}
The terms of the companion EDS 
(\ref{cEDS}) satisfy the Somos-4 recurrence 
\beq\label{cs4}
\ra_{n+4}\ra_n=\tmu^2\,\ra_{n+3}\ra_{n+1}-\tal\,\ra_{n+2}^2,  
\eeq 
and they 
can be written as polynomials in $\tmu,\tal,\tbe$ with integer coefficients, beginning with 
$$  
(\ra_n): 
0,1,-\tmu,\tal,\tmu\tbe, -\tal^3-\tmu^4\tbe, \tmu\tal(\tal^3+\tbe^2+\tmu^4\tbe),
-\tal^6-\tmu^4\tbe(\tal^3-\tbe^2),
-\tmu\tbe\big(2\tal^6+\tal^3\tbe(\tbe+3\tmu^4)+\tmu^8\tbe^2\big),
\ldots .
$$ 
A general Somos-5 sequence satisfies infinitely many higher Somos 
relations of odd order with coefficients determined by its companion EDS  (\ref{cEDS}), namely 
\beq\label{somosk}
\ra_2 \,\tau_{n+2j+1}\tau_n=\ra_j\ra_{j+1}\,\tau_{n+j+2}\tau_{n+j-1}-\ra_{j-1}\ra_{j+2}\,\tau_{n+j+1}\tau_{n+j}. 
\eeq   
\end{thm} 
\begin{prf} The description of the terms of the  companion EDS  as polynomials in $\tmu,\tal,\tbe$ which are generated by a 
Somos-4 recurrence was given 
in \cite{swahon}. 
A proof of the  higher Somos relations (\ref{somosk}) was provided  in 
\cite{hones5} (see also \cite{swartvdp}).
\end{prf} 
\begin{remark}
When $\tmu,\tal,\tbe\in\Z$, the companion EDS is an integer sequence, so it is a bona fide elliptic divisibility sequence in the sense of \cite{ward}. 
For $j=0,1$ the relation (\ref{somosk}) is just a tautology, while for $j=2$ it is equivalent to (\ref{s5gen}) with 
coefficients $\tal=\ra_3=\si(3\ka)/\si(\ka)^9$, $\tbe=-\ra_4/\ra_2=-\si(\ka)^{-12}\si(4\ka)/\si(2\ka)$.   
\end{remark}

Having described the general case, in the rest of this section the formula  
(\ref{genivp}) and the other results on  Somos-5 sequences given above will be specialized  
to the particular sequences 
 (\ref{s5seq}) and (\ref{s5eds}). 

\begin{propn}\label{Snform}  
The terms of the sequence (\ref{s5seq}) are 
given by the formula 
\beq\label{Sn} 
S_n = 
\begin{cases} 
 \qquad{\si (\om)}^{-1}B_+^{\lc n/2\rc} \mu^{\lc n/2\rc^2} {\si (\om +n\kappa)} , 
    & \text{for }  n \, \,\mathrm{even} \\
  \,{\si(\om+\ka)}^{-1}     B_-^{\lc n/2\rc} \mu^{\lc n/2\rc^2} {\si (\om+n\kappa)},  
 & \text{for }  n \,\,\mathrm{odd} , 
  \end{cases}
\eeq 
with 
parameters 
given by 
\beq\label{s5params}
\tmu=6^{\frac{1}{4}}=\wp'(\ka), \quad \tilde{\la}=\frac{29}{12\sqrt{6}}=\wp(\ka), \quad g_2=\frac{121}{72}, \quad g_3=-\frac{845}{1296\sqrt{6}}, 
\eeq 
and  $\mu= -\si(\ka)^{-4}$, 
\beq\label{bpm}
B_+=-\frac{\si(\ka)^4\si(\om)}{\si(\om+2\ka)}, \qquad 
B_-=-\frac{\si(\om-\ka)}{\si(\ka)^4\si(\om+\ka)}, 
\eeq 
where the numerical value 
\beq\label{kaval}
\kappa\approx - 1.052799817
\eeq
determines the point ${\cal P}=\big(\tilde{\la},\tmu\big)$ on the Weierstrass curve (\ref{origwei}) with these values of the invariants $g_2,g_3$, 
and the initial point $\hat{{\cal P}}_0=\big(\wp(z_0),\wp'(z_0)\big)$ is 2-torsion, that is 
\beq\label{2tors} 
z_0=\om, \qquad \wp(\om) =\frac{5}{12\sqrt{6}}, \qquad \wp'(\om)=0, 
\eeq 
where $\om$ is a half-period, which can be taken as the sum of real and imaginary half-periods $\om_1$, $\om_2$: 
\beq\label{omval} 
\om=\om_1+\om_2, \qquad \om_1\approx 1.849876692, \qquad \om_2\approx 1.524280920\, \ri .
\eeq  
\end{propn} 
\begin{prf} This sequence was presented as an example in \cite{hones5}:  the coefficients in (\ref{s5orig}) are  $\tal=\tbe=1$, 
while  the values of the conserved quantities (\ref{It}) and  (\ref{Jt}) are given by  $\tilde{I}=\tilde{J}=5$, 
leading to the parameter values (\ref{s5params}). 
The formulae 
in Theorem 2.7 of \cite{hones5} then show that  $z_0=\om$ corresponds to a 2-torsion point on the curve  in this case, 
and the numerical  values of $\ka$ and $\om$ are determined from (\ref{s5params}) and (\ref{2tors}) by evaluating elliptic integrals. 
\end{prf} 

From a purely algebraic point of view, it would seem more natural to apply a homothety so that everything is defined over $\Q$, i.e. rescale all the coordinates, Weierstrass functions 
and invariants by suitable powers of $\tmu=6^{\frac{1}{4}}$ in order to  work with the  
Weierstrass cubic
\beq\label{weier}
E: \quad y^2=4x^3-\frac{121}{12}x+\frac{845}{216};
\eeq 
the  expressions corresponding to this form of the curve are presented in \cite{hones5}. 
However, the analytic calculations in the sequel are easier to carry out with the choice of scale as in (\ref{s5params}). 

\begin{remark} 
The relation 
(\ref{constr}) between the ratio of the quantities in (\ref{bpm}) can be verified by using the elliptic function identity 
\beq\label{ellipid}
\frac{\si(z_0+2\ka)\si(z_0-\ka)}{\si(\ka)^4\si(z_0)\si(z_0+\ka)}=-\frac{\wp'(\ka)}{2}\,\left(\frac{\wp'(z_0)-\wp'(\ka)}{\wp(z_0)-\wp(\ka)}\right)+\frac{\wp''(\ka)}{2},  
\eeq 
which is valid for any $z_0,\ka$. Upon setting $z_0=\om$, the left-hand side of the above identity is 
\beq\label{bratio}
\frac{\si(\om+2\ka)\si(\om-\ka)}{\si(\ka)^4\si(\om)\si(\om+\ka)}=-\mu^{-1}\frac{B_-}{B_+}, 
\eeq
while the right-hand side is 
$$ 
-\frac{\wp'(\ka)}{2}\,\left(\frac{\wp'(\om)-\wp'(\ka)}{\wp(\om)-\wp(\ka)}\right)+\frac{\wp''(\ka)}{2} 
=-\frac{6^{\frac{1}{4}}}{2} \, \left(\frac{-6^{\frac{1}{4}}}{\frac{5}{12\sqrt{6}} -\frac{29}{12\sqrt{6}}  }\right) +\frac{5}{2} = 1,   
$$
where above we have substituted the values from (\ref{s5params}) and (\ref{2tors}),  as well as $\tilde{J}=\wp''(\ka)=5$. 
\end{remark}


\begin{propn}  \label{Tnform}  
The terms of the sequence 
of   (\ref{s5eds}) can be written in the form 
\beq\label{Tn} 
T_n = 
\begin{cases} 
 \quad\,{\si (\ka)}^{-1}\hat{B}_+^{\lc (n-1)/2\rc} \mu^{\lc (n-1)/2\rc^2} {\si (n\kappa)} , 
    & \text{for }  n \, \,\mathrm{odd} \\
  -{\si(2\ka)}^{-1}     \hat{B}_-^{\lc (n-1)/2\rc} \mu^{\lc (n-1)/2\rc^2} {\si (n\kappa)},  
 & \text{for }  n \,\,\mathrm{even} , 
  \end{cases}
\eeq 
where  
\beq\label{bhatpm}
\hat{B}_+=\mu, \qquad 
\hat{B}_-=-\mu^2, 
\eeq 
with  $\mu= -\si(\ka)^{-4}$
and the same value of $\ka$ as  
in Proposition \ref{Snform}. 
\end{propn} 
\begin{prf} 
For the second sequence  (\ref{s5eds}), again we have $\tal=\tbe=1$, 
while the conserved quantity (\ref{It}) takes the value $\tilde{I}=7$, and (\ref{Jt}) has the same 
value $\tilde{J}=5$. The fact that the values of $\tal,\tbe,\tilde{J}$ coincide with those for  (\ref{s5seq}) means that the two sequences are very closely related. Since 
$T_0=0$, it is convenient to specify the initial value problem with  $T_1,\ldots,T_5$,  in order to apply the general formula (\ref{genivp}); effectively this corresponds to 
shifting the index by 1 and changing the parity. The result (\ref{Tn}) then follows.
\end{prf}

Writing the terms of the sequence $(T_n)$ in the analytic form (\ref{Tn}) will be useful for comparing it with the terms of $(S_n)$ in the sequel, but 
disguises the fact that,   
up to a rescaling of terms with even/odd index, (\ref{s5eds}) coincides with the companion EDS for (\ref{s5seq}), as defined by (\ref{cEDS}) above. 
Indeed,  a simpler way to write the terms of (\ref{s5eds}) is as 
\beq\label{Tnalt} 
T_n = 
\begin{cases}\quad \,\,\,\,\ra_n =
 {\si (n\kappa)}/ {\si (\ka)}^{n^2} , 
    & \text{for }  n \, \,\mathrm{odd} \\
\tmu^{-1}\ra_n =
\tmu^{-1}
{\si (n\kappa)}/ {\si (\ka)}^{n^2} , 
 & \text{for }  n \,\,\mathrm{even} . 
  \end{cases}
\eeq
By virtue of its being the companion EDS of $(S_n)$,  rescaled according to the parity  of $n$, the sequence $(T_n)$ satisfies many identities that intertwine it with    
(\ref{s5seq}), as illustrated by  the following result.  

\begin{propn}\label{detprop} The terms of the sequences (\ref{s5seq}) and (\ref{s5eds})  satisfy the determinantal identities 
$$ 
\left|\begin{array}{ccc} 
S_{n-1} T_0 & S_n T_1 & S_{n+j-1} T_j \\ 
S_{n+j-1} T_{-j} & S_{n+j} T_{-j+1} & S_{n+2j-1} T_0 \\ 
S_{n+j+1} T_{-j-2} & S_{n+j+2} T_{-j-1} & S_{n+2j+1} T_{-2}  
\end{array} 
\right|=0  
$$ 
and 
$$ 
\left|\begin{array}{ccc} 
T_{n-1} T_0 & T_n T_1 & T_{n+j-1} T_j \\ 
T_{n+j-1} T_{-j} & T_{n+j} T_{-j+1} & T_{n+2j-1} T_0 \\ 
T_{n+j+1} T_{-j-2} & T_{n+j+2} T_{-j-1} & T_{n+2j+1} T_{-2}  
\end{array} 
\right|=0,  
$$ 
for all $j,n\in\Z$. More generally, in 
each of the infinite matrices with entries $(S_{n+2i+j}T_{-2i+j})$, $(T_{n+2i+j}T_{-2i+j})$ 
for $(i,j)\in\Z^2$,  the determinant of any  minor of size $3\times 3$ or above vanishes for all 
$n\in\Z$. 
\end{propn} 
\begin{prf}
Expanding out each of the $3\times 3$ determinants above, noting that $T_0=0$, $T_1=1$ and $T_{-n}=-T_n$, 
and removing a common factor  
leads to Somos-type relations in $n$ (for fixed $j$), that is 
\beq\label{higher}
T_1T_2S_{n+2j+1}S_{n}=T_jT_{j+1}S_{n+j+2}S_{n+j-1}-T_{j-1}T_{j+2}S_{n+j+1}S_{n+j} 
\eeq
and 
$$ 
T_1T_2T_{n+2j+1}T_{n}=T_jT_{j+1}T_{n+j+2}T_{n+j-1}-T_{j-1}T_{j+2}T_{n+j+1}T_{n+j}. 
$$ 
The first identity follows directly from (\ref{somosk}), upon replacing $\tau_i\to S_i$ and $\ra_i\to T_i$, since by 
(\ref{Tnalt}) the sequence 
$(T_n)$ is the companion EDS of $(S_n)$ up to rescaling terms of opposite parity, and only products of pairs 
of even/odd index appear. The second identity follows in the same way, by replacing $\tau_i\to T_i$ and $\ra_i\to T_i$ 
in (\ref{somosk}), since up to parity the sequence $(T_n)$ is its own companion EDS. 
The more general statement about vanishing determinants of $3\times 3$ minors follows by making the same 
replacements in the infinite matrix with $(i,j)$ entries $(\ra_{-2i+j}\,\tau_{n+2i+j})$. 
For example, when $n=0$ the rows of $(S_{n+2i+j}T_{-2i+j})$ with $i=0,\ldots 4$ include the entries 
$$ 
\begin{array}{cccccccccc}  
\cdots & -2 & 1 & -1 & 0 & 1 & -1 & 2 & 3 & \cdots \\
\cdots & 7 & -1 & -1 & 1 & -2 &  0 & 5 & -11 & \cdots \\
\cdots & 1 & -8 & 14 & -3 & -5 & 11 & -37 & 0 & \cdots \\ 
\cdots & -782 & 171 & 5 & -88 & 259 & -83 & -274 & 1217 & \cdots \\ 
\cdots &  13645 & 5005 & -14467 & 4731 & 274 & -9736 & 43127 & -22833 & \cdots 
\end{array} 
$$ 
and one of the vanishing minors that does not correspond to a quadratic (Somos-type) relation, but rather is cubic in $T_n$, is 
$$ 
\left|\begin{array}{ccc}
-1 & 1 & 2 \\ 
14 & -5 & -37 \\ 
-14467 & 274 & 43127 
\end{array}
\right|=0.
$$ 
The determinant 
of any $3\times 3$ minor of 
$(\ra_{-2i+j}\,\tau_{n+2i+j})$ 
has the form 
$$ 
D=
\left|\begin{array}{lll} 
\ra_{-2i+j} \,\tau_{n+2i+j} & \ra_{-2i+j'} \,\tau_{n+2i+j'} & \ra_{-2i+j''} \,\tau_{n+2i+j''}  \\ 
\ra_{-2i'+j} \,\tau_{n+2i'+j} & \ra_{-2i'+j'} \,\tau_{n+2i'+j'} & \ra_{-2i'+j''} \,\tau_{n+2i'+j''}  \\ 
\ra_{-2i''+j} \,\tau_{n+2i''+j}  & \ra_{-2i''+j'} \,\tau_{n+2i''+j'}  & \ra_{-2i''+j''} \,\tau_{n+2i''+j''}   
\end{array} 
\right|.
$$
The $\tau_i$ terms appearing in each column have indices $i$ with the same parity, so from the formula (\ref{genivp}) 
it follows that there is a common factor of $A_{\pm}$ that can be removed from each column; then effectively we can 
ignore these prefactors (equivalently, by a gauge transformation the even/odd index terms can always be rescaled 
separately so that $A_{\pm}\to 1$). For convenience, we introduce the notation $\si_\zeta=\si(\zeta\ka)$ for any 
$\zeta\in \C$, and let $\bar{n}=n+z_0/\ka$, so that upon substituting from  (\ref{genivp}) and (\ref{cEDS}) we find 
$$ 
D\propto 
\left|\begin{array}{ccc} 
B_0^{F_{ij}} 
C_{ij}
\, 
\si_{-2i+j} \,\si_{\bar{n}+2i+j}
 & 
B_1^{F_{ij'}}
C_{ij'} \, 
\si_{-2i+j'} \,\si_{\bar{n}+2i+j'}  
& 
B_2^{F_{ij''}}
C_{ij''}\, 
\si_{-2i+j''} \,\si_{\bar{n}+2i+j''} \\   
B_0^{F_{i'j}} 
C_{i'j}
\, 
\si_{-2i'+j} \,\si_{\bar{n}+2i'+j}  &
B_1^{F_{i'j'}}
C_{i'j'} \, 
\si_{-2i'+j'} \,\si_{\bar{n}+2i'+j'}  
& 
B_2^{F_{i'j''}}
C_{i'j''} \, 
\si_{-2i'+j''} \,\si_{\bar{n}+2i'+j''} \\   
B_0^{F_{i''j}} 
C_{i''j} 
\, 
\si_{-2i''+j} \,\si_{\bar{n}+2i''+j}  &
B_1^{F_{i''j'}}
C_{i''j'} \, 
\si_{-2i''+j'} \,\si_{\bar{n}+2i''+j'}  
& 
B_2^{F_{i''j''}}
C_{i''j''} \, 
\si_{-2i''+j''} \,\si_{\bar{n}+2i''+j''}
%
\end{array} 
\right|,
$$ 
where $C_{ij}=(-1)^{F_{ij}^2}
\si_1^{-E_{ij}}$, $E_{ij}=(-2i+j)^2 +4\lc(n+2i+j)/2\rc^2$, $F_{ij}=\lc(n+2i+j)/2\rc$, 
and $B_0=B_{\pm}$, $B_1=B_{\pm}$, $B_2=B_{\pm}$, according to the parity of $n+2i+j$, $n+2i+j'$, $n+2i+j''$, respectively. 
Depending on the parities of the latter quantities, a case by case analysis shows that the rows and columns can be rescaled appropriately so that 
the terms depending on powers of $B_\pm$, an overall sign, and the powers of $\si_1=\si(\ka)$ can be removed. The analysis relies on various identities for the 
exponents $E_{ij}$; for instance, if $n+2i+j$ and $n+2i+j'$ have the same parity then $E_{ij}+E_{i'j'}=E_{i'j}+E_{ij'}$, but if they have opposite parity this is not the case 
and it is necessary to use (\ref{constr}) to balance the powers of $\si_1$ that appear. This gives an overall factor of 
$C_{ij}C_{i'j'}C_{i''j''}B_0^{F_{ij}} B_1^{F_{i'j'}}B_2^{F_{i''j''}}$ in front, and what remains is 
$$ 
D\propto 
\left|\begin{array}{ccc} 
\si_{-2i+j} \,\si_{\bar{n}+2i+j}
 & 
\si_{-2i+j'} \,\si_{\bar{n}+2i+j'}  
& 
\si_{-2i+j''} \,\si_{\bar{n}+2i+j''} \\   
\si_{-2i'+j} \,\si_{\bar{n}+2i'+j}  &
\si_{-2i'+j'} \,\si_{\bar{n}+2i'+j'}  
& 
\si_{-2i'+j''} \,\si_{\bar{n}+2i'+j''} \\   
\si_{-2i''+j} \,\si_{\bar{n}+2i''+j}  &
\si_{-2i''+j'} \,\si_{\bar{n}+2i''+j'}  
& 
\si_{-2i''+j''} \,\si_{\bar{n}+2i''+j''}
\end{array}
\right|=\frac{1}{\si_{-2i'+j} \,\si_{\bar{n}+2i'+j}  } 
\left|\begin{array}{cc} D_{NW} & D_{NE} \\ 
D_{SW} & D_{SE} 
\end{array} \right|,
$$ 
where in the last step we have used Dodgson condensation \cite{dodgson} to expand the $3\times 3$ determinant in terms of 
its $2\times 2$ connected minors. In particular, using the standard three-term relation for the Weierstrass sigma function 
(see e.g. $\S20.53$ in \cite{ww}) we have  
$$ 
D_{NW}=\left|\begin{array}{cc} 
\si_{-2i+j} \,\si_{\bar{n}+2i+j}
 & 
\si_{-2i+j'} \,\si_{\bar{n}+2i+j'} \\ 
\si_{-2i'+j} \,\si_{\bar{n}+2i'+j}  &
\si_{-2i'+j'} \,\si_{\bar{n}+2i'+j'} 
\end{array} \right| = \si_{\bar{n}+2i+2i'} \si_{2i-2i'}\si_{\bar{n}+j+j'}\si_{j'-j},
$$ 
and similar calculations for the $NE$, $SW$ and $SE$ minors, together with the fact that $\si$ is an odd function,  yield 
$$ 
\begin{array}{rcl} D & \propto &  D_{NW}D_{SE} -D_{NE}D_{SW} \\  
&=& 
\si_{\bar{n}+2i+2i'} \si_{2i-2i'}\si_{\bar{n}+j+j'}\si_{j'-j}\times \si_{\bar{n}+2i'+2i''} \si_{2i'-2i''}\si_{\bar{n}+j'+j''}\si_{j''-j'} \\
&&  - 
\si_{\bar{n}+2i+2i'} \si_{2i'-2i}\si_{\bar{n}+j'+j''}\si_{j''-j'} \times \si_{\bar{n}+2i'+2i''} \si_{2i''-2i'}\si_{\bar{n}+j+j'}\si_{j'-j} =0, 
\end{array} 
$$ 
as required. Since any of the larger minors can be expanded in terms of $3\times 3$ minors, these all vanish as well.
\end{prf} 

\begin{remark} The vanishing of the analogous $3\times 3$ determinants involving the original Somos-5 sequence  (\ref{s5seq}) 
is proved in \cite{ma}.
\end{remark}

\subsection{Arithmetical properties of Somos-5 sequences} 

We now consider arithmetical properties of these Somos-5 sequences that will be needed in the sequel. 
The following result concerning  (\ref{s5seq})  is well known; the original proof is attributed to Bergman \cite{gale}.

\begin{lem}\label{copr}
Any five adjacent terms in the sequences  (\ref{s5seq})  and (\ref{s5eds})  are pairwise coprime, that is, 
$\gcd(S_i,S_j)=1=\gcd(T_i,T_j)$ for $|i-j|<5$.
\end{lem}
\begin{prf}
The proof is by induction. Each of the sequences is a solution 
of the recurrence 
\beq\label{s5original}
\tau_{n+5}\tau_n=\tau_{n+4}\tau_{n+1}+\tau_{n+3}\tau_{n+2}, 
\eeq
and in both cases the statement is clearly true for the first five terms, indexed by $n=1,\ldots,5$. So if we suppose the inductive hypothesis 
that $\tau_n,\tau_{n+1},\ldots,\tau_{n+4}$ are pairwise coprime, and assume that some prime $p$ is a common factor of $\tau_{n+5}$ and 
$\tau_{n+1}$, then we see from the recurrence that $p|\tau_{n+3}\tau_{n+2}$, which contradicts the fact that 
$\gcd(\tau_{n+1},\tau_{n+2})=1=\gcd(\tau_{n+1},\tau_{n+3})$, and similar contradictions arise from assuming that 
$p$ is a common factor of  $\tau_{n+5}$ and one of $\tau_{n+2},\tau_{n+3},\tau_{n+4}$. 
\end{prf}

\begin{remark}
Pairwise coprimeness of adjacent terms is a feature of clusters of Laurent polynomials in cluster algebras \cite{fz1}, and 
more generally in various birational difference equations with the Laurent property \cite{coprime}, where essentially the same argument applies. 
For the original Somos-5 sequence (\ref{s5seq}),  Robinson 
actually proved the stronger statement that  $\gcd(S_i,S_j)=1$ for $|i-j|\leq5$ \cite{rob}, but this statement is not quite true for the sequence 
(\ref{s5eds}), because 7 is a common factor of all the terms $T_{5j}$. 
\end{remark}

Robinson used elementary methods to prove the periodicity of Somos-4 and Somos-5 sequences modulo any positive integer, for the case of coefficients $\tal=\tbe=1$ with all the rational 
initial data being units in the corresponding residue ring, and made conjectures about the periods modulo a prime or a prime power. Using the connection with elliptic 
curves, several of these conjectures 
were proved in the thesis of Swart \cite{swart}, who also considered the case of general coefficients in (\ref{s5gen}), and some of these results were further strengthened by 
van der Kamp \cite{div}. 
To begin with, we would like to adapt Robinson's arguments to a slightly more 
general class of initial data, and consider periodicity modulo a prime, which is the case of most interest for us; the extension to prime powers and more general 
moduli is quite straightforward. 

It is convenient to formulate conditions on the initial data  in terms of the $p$-adic norm $|\cdot|_p$. There is a wealth of literature on rational maps 
of the projective line $\Proj^1$ over $\Q_p$ \cite{silverman}, and $p$-adic analysis is extremely useful for understanding suitable notions 
of good/bad reduction modulo a prime for maps over $\Q$. For nonlinear systems in higher dimensions, results are rather more sparse, although 
the authors of \cite{kanki} have proposed a definition of (almost) good reduction for birational maps, and used local analysis in $\Q_p$
to probe the singularity structure of certain maps in the plane.  
Here we will treat only the most relevant case of the recurrence (\ref{s5original}), which is equivalent to the birational map 
\beq\label{birat} 
(\tau_0,\tau_1,\tau_2,\tau_3,\tau_4)\mapsto 
\left(\tau_1,\tau_2,\tau_3,\tau_4, \frac{\tau_4\tau_1+\tau_3\tau_2}{\tau_0}\right)
\eeq 
in five dimensions; 
this avoids having to specify additional conditions on the coefficients, 
but essentially the same periodicity statements hold for (\ref{s5gen}) if we require that 
both of $\tal,\tbe\in\Z_p$, with at least one of them being a  $p$-adic unit. Since the choice of where to 
start indexing the sequence is arbitrary, the initial value problem for 
 (\ref{s5original}) will be specified either by the values $(\tau_0,\tau_1,\tau_2,\tau_3,\tau_4)$ or 
$(\tau_1,\tau_2,\tau_3,\tau_4, \tau_5)$. 

\begin{definition} \label{baldef}
For a prime $p$, an initial value problem for (\ref{s5original}) over $\Q$ is said to be \textit{well-balanced} $\bmod \,p$  if 
it is specified 
by either  $(\tau_0,\tau_1,\tau_2,\tau_3,\tau_4)$ or 
$(\tau_1,\tau_2,\tau_3,\tau_4, \tau_5)\in\Q^5$ such that four adjacent initial values are $p$-adic units, 
so $|\tau_1|_p=|\tau_2|_p=|\tau_3|_p=|\tau_4|_p=1$, and 
\beq\label{balance} 
|\tau_4\tau_1+\tau_3\tau_2|_p\leq |\tau_*|_p\leq 1, \qquad \tau_*\neq0,
\eeq   
where $\tau_*=\tau_0$ or $\tau_5$ accordingly.
\end{definition}

\begin{thm}\label{periods}
For any prime $p$, if an initial value problem for (\ref{s5original})  is well-balanced $\bmod \,p$, then 
the rational  Somos-5 sequence $(\tau_n)$  is defined for all $n\in\Z$, as is the reduced sequence $(\tau_n\bmod p)$, which 
is periodic in $n$.
\end{thm}

\begin{prf}
Starting from non-zero initial values, the birational map (\ref{birat}) can be iterated both forwards, to obtain $\tau_n\in \Q$ for $n\geq 0$, 
and backwards (using the inverse map) to obtain the terms with negative indices, provided that a singularity does not appear, i.e.\ unless a 
zero term appears in the sequence. However, it is still possible to continue the sequence beyond a zero (in either direction) by making use of the Laurent 
property: we will perform the analysis for iterating forwards, and then the result for the reverse direction follows from the symmetry of   
the recurrence (\ref{s5original})  under sending $n\to -n$. Suppose that, for some $n$, there are five non-zero terms $\tau_n,\tau_{n+1},\tau_{n+2},\tau_{n+3},\tau_{n+4}$, followed 
by $\tau_{n+5}=0$, which occurs when 
\beq\label{zerocond}
\frac{\tau_{n+4}\tau_{n+1}}{\tau_{n+3}\tau_{n+2}}=-1.
\eeq 
The main point is that, due to the Laurent phenomenon \cite{fz}, all subsequent terms (and also all previous terms) can be written as Laurent polynomials in these five non-zero terms, with integer 
coefficients - in fact, it has even been proved that the coefficients are positive integers \cite{ghkk,ls}, so we have $\tau_{n+j}\in\Z_{>0}[\tau_n^{\pm1},\tau_{n+1}^{\pm1},\tau_{n+2}^{\pm1},\tau_{n+3}^{\pm1},\tau_{n+4}^{\pm1}]$ 
for all $j\in\Z$. Hence the rational sequence is defined for all $n\in\Z$, simply by evaluating this sequence of Laurent polynomials at any five adjacent non-zero values - in particular, evaluating them at the five well-balanced 
initial values. 

For the reduction $\bmod \, p$, some further analysis is helpful. If all the initial values are $p$-adic units (the case considered by Robinson), then 
their reduction  
$\bmod \,p$ defines an initial value problem for 
(\ref{birat}) as a birational map $\F_p^5\to\F_p^5$; a singularity may be reached under iteration over $\F_p$, if for some $n$ the five non-zero terms $\tau_n\bmod \, p, \ldots,\tau_{n+4}\bmod\, p$ are followed by 
$\tau_{n+5}\bmod \,p=0$, but nevertheless the sequence is still defined for all $n\in\Z$ by evaluating the reduction $\bmod \, p$ of the Laurent polynomials, belonging to the ring 
$\F_p[\tau_0^{\pm1},\tau_{1}^{\pm1},\tau_{2}^{\pm1},\tau_{3}^{\pm1},\tau_{4}^{\pm1}]$ (which can be evaluated on any set of five adjacent non-zero values in $\F_p^*$). Now if a zero term never appears 
in the sequence $(\tau_n\bmod \,p)$, there are only $(p-1)^5$ possible quintuples in $(\F_p^*)^5$, so by the pigeonhole principle the sequence is periodic and this provides a crude upper bound on the period. 
However, if a zero appears somewhere, which is certainly the case for a rational initial value problem with either $|\tau_0|_p<1$ or $|\tau_5|_p<1$  appearing before/after four adjacent units, then 
the continuation of the sequence needs a more careful treatment. Note that, in either of these cases, the well-balanced condition (\ref{balance}) implies 
$$ 
\tau_5\tau_0 = \tau_4\tau_1+\tau_3\tau_2\qquad \mathrm{with} \quad |\tau_5|_p\leq 1 \quad\mathrm{or}\quad |\tau_0|_p\leq 1, 
$$ 
respectively, since in the first case the above relation defines $\tau_5$ that appears after the five initial data, and in the second case it defines $\tau_0$ that precedes it; 
but it may happen that both $\tau_0\bmod\,p=0$ and $\tau_5\bmod\,p=0$ (as in the example of the sequence (\ref{s5eds}) for $p=7$). So we consider a more general setting  where we have four 
adjacent non-zero values $\tau_{n+1},\ldots,\tau_{n+4}$ subject to the condition (\ref{zerocond}) holding in $\F_p$, followed by $\tau_{n+5}=\eps\equiv 0\bmod \, p$, and preceded 
by $\tau_n$ which is a $p$-adic integer, but may or may not be a unit; this covers both the case of a zero appearing under iteration in $\F_p$, and (up to reversing the direction of iteration) the case 
where one of the initial values is a non-unit. Thus, by iterating and then reducing $\bmod\,p$, we find 
\beq\label{calc1}
\bear{rcl} 
\tau_{n+6} & = & \tau_{n+1}^{-1}(\tau_{n+3}\tau_{n+4}+\tau_{n+2}\,\eps) 
 \equiv  \tau_{n+1}^{-1}\tau_{n+3}\tau_{n+4}; \\
\tau_{n+7}& =&(\tau_{n+1}\tau_{n+2})^{-1}\big(\tau_{n+3}^2\tau_{n+4}+(\tau_{n+2}\tau_{n+3}+\tau_{n+1}\tau_{n+4})\,\eps\big) \\ 
& = & (\tau_{n+1}\tau_{n+2})^{-1}\big(\tau_{n+3}^2\tau_{n+4}+\tau_{n}\,\eps^2\big) 
 \equiv   (\tau_{n+1}\tau_{n+2})^{-1}\tau_{n+3}^2\tau_{n+4}; \\
\tau_{n+8}& =&(\tau_{n+1}\tau_{n+2})^{-1}\tau_{n+3}\tau_{n+4}^2+\tau_{n+1}^{-1}\tau_{n+4}\,\eps+O(\eps^2) 
\equiv   (\tau_{n+1}\tau_{n+2})^{-1}\tau_{n+3}\tau_{n+4}^2; \\ 
\tau_{n+9}& =& \tau_{n+1}^{-2}\tau_{n+2}^{-1}\tau_{n+3}\big(\tau_{n+3}^2\tau_{n+4}+(\tau_{n+4}\tau_{n+1}+\tau_{n+3}\tau_{n+2})\,\eps\big)+O(\eps^2) \\ 
& = & \tau_{n+1}^{-2}\tau_{n+2}^{-1}\tau_{n+3}^3\tau_{n+4}+O(\eps^2) 
\equiv  \tau_{n+1}^{-2}\tau_{n+2}^{-1}\tau_{n+3}^3\tau_{n+4}; \\ 
\tau_{n+10}& =& \tau_{n+1}^{-3}\tau_{n+2}^{-2}\tau_{n+3}^2\tau_{n+4}\big(\tau_{n+3}\tau_{n+4}\,\eps^{-1}+\tau_{n+2}\big)(\tau_{n+4}\tau_{n+1}+\tau_{n+3}\tau_{n+2})+O(\eps) \\
&=& \tau_{n+1}^{-3}\tau_{n+2}^{-2}\tau_{n+3}^3\tau_{n+4}^2\tau_n+O(\eps) \equiv \tau_{n+1}^{-3}\tau_{n+2}^{-2}\tau_{n+3}^3\tau_{n+4}^2\tau_n. 
\eear
\eeq 
The cancellation of $\eps$ from the denominator appearing in $\tau_{n+10}$ is precisely what yields the Laurent property, and  the reduction $\bmod \,p$ 
gives four adjacent units $\tau_{n+6},\tau_{n+7},\tau_{n+8},\tau_{n+9}$, followed by $\tau_{n+10}$ which is a $p$-adic integer (and is a unit whenever $\tau_n$ is); indeed, we have 
$$|\tau_{n+9}\tau_{n+6}+\tau_{n+8}\tau_{n+7}|_p=|\tau_{n+5}|_p|\tau_{n+10}|_p\leq|\tau_{n+10}|_p\leq \max (|\tau_n|_p,|\tau_{n+5}|_p)\leq1,$$ 
so this is well-balanced (and the first inequality is strict) unless 
it happens that 
$\tau_{n+10}=0$ holds in $\Q$. Nevertheless, by shifting indices up by 5, the above calculation shows that one can write the next four terms $\tau_{n+11},\ldots,\tau_{n+14}$ 
as polynomials in $\tau_{n+10}$ with coefficients in  $\Z[\tau_{n+6}^{-1},\tau_{n+7}^{-1},\tau_{n+8}^{-1},\tau_{n+9}^{-1}]$, so 
for instance $\tau_{n+11}=\tau_{n+6}^{-1}\tau_{n+8}\tau_{n+9}+O(\tau_{n+10})$, while at the fifth step we find 
\beq\label{calc2} 
\tau_{n+15}= \tau_{n+6}^{-3}\tau_{n+7}^{-2}\tau_{n+8}^3\tau_{n+9}^2\tau_{n+5}+O(\tau_{n+10}).
\eeq 
Hence by induction  any well-balanced initial value problem determines a rational sequence $(\tau_n)$ consisting of $p$-adic integers, so the entire sequence $(\tau_n\bmod\,p)_{n\in\Z}$ is well-defined in $\F_p$, and by the 
pigeonhole principle it is periodic.   
\end{prf}

\begin{remark}  In the sequence (\ref{s5eds}), 
 $(T_0,\ldots,T_4)=(0,1,-1,1,1)$ does not provide well-balanced initial data 
for any prime  $p$, due to the initial zero, which leaves the value of $T_5$ undetermined. 
In contrast, $(T_1,\ldots,T_4)=(1,-1,1,1,-7)$ is well-balanced for any $p$, in particular for $p=7$, and comparing with (\ref{calc1}) and (\ref{calc2}) it is clear 
that 
$T_{5j}\bmod\, 7=0$ for all $j$, as asserted in the previous remark (although the actual period of the sequence $\bmod\, 7$ is 20 \cite{rob}). As another example, 
the initial data $(\tau_1,\ldots,\tau_5)=(1,-1,1,8,49)$ is not well-balanced $\bmod\,2$ because $\tau_4$ is not a unit, nor is it well-balanced   $\bmod\,7$ because 
it fails the condition (\ref{balance}) on the norm; indeed, $\tau_0=1/7$, and in fact the corresponding rational sequence does not admit reduction modulo either of 
these primes, as it has growing powers of 2 and 7 appearing as   denominators. (The growth of both Archimedean and non-Archimedean valuations for Somos-4 sequences is 
described in \cite{uchida}, and these results are relevant here because it is known that the even/odd index terms in a Somos-5 sequence each define a Somos-4, as shown  in \cite{hones5}, 
where  asymptotic results were obtained in the Archimedean case.)
Yet one more example is $(\tau_1,\ldots,\tau_5)=(1,-1,1,8,-7)$, which is well-balanced $\bmod\, 7$ and for all other primes except  $\bmod\,2$; however, an interesting feature of 
the fact that the fifth term is not a $7$-adic unit is that, although the initial values satisfy $\tau_j\bmod\,7=T_j\bmod 7$ for $j=1,\ldots,5$, the two sequences have a different reduction 
$\bmod\,7$ because e.g.\ $\tau_0=-1\not\equiv 0\,(\bmod\,7)$, and this phenomenon does not arise in the case that all the initial values are units. 
\end{remark}

For the proof of Theorem \ref{main} in the next section, we will need the following result about the particular sequences (\ref{s5seq})  and (\ref{s5eds}).

\begin{propn}\label{period23}
For the Somos-5 sequences   (\ref{s5seq})  and (\ref{s5eds}), 
$$ 
S_n\bmod 2=0 \quad \mathrm{iff} \quad n\equiv 3 \,(\bmod \, 6), \qquad 
T_n\bmod 2=0 \quad \mathrm{iff} \quad n\equiv 0 \,(\bmod \, 6),
$$
and both these reduced sequences $\bmod\,2$ have period 6, while 
$$ 
S_n\bmod 3=0 \quad \mathrm{iff} \quad n\equiv 4 \,(\bmod \, 8), \qquad 
T_n\bmod 3=0 \quad \mathrm{iff} \quad n\equiv 0 \,(\bmod \, 8),
$$ 
where both reduced sequences $\bmod\,2$ have period 16. 
Furthermore, for every prime $p$ there are infinitely many $n$ such that 
$T_n\bmod p=0$.
\end{propn}
\begin{prf} 
The initial values $(S_0,\ldots,S_4)$ give 
$(S_0\bmod 2,\ldots,S_4\bmod 2)=(1,1,1,0,1)$, 
so these are not well-balanced $\bmod \,2$. By shifting back one or two  steps and starting with $S_{-1}$ 
or $S_{-2}$ we can get something well-balanced, and use the method of reduction of Laurent polynomials, as in the proof of  
Theorem \ref{periods}. However, instead we will apply Proposition \ref{detprop}, 
making use of the fact that the sequence $(S_n)$ satisfies the Somos-7 relation 
$$ 
S_{n+7}S_n = -S_{n+5}S_{n+2}+7S_{n+4}S_{n+3}, 
$$ 
given by the case $j=3$ of (\ref{higher}), and $(T_n)$ satisfies the same recurrence. 
So if we take $ (1,1,1,0,1)$ as initial values in $\F_2$, we can iterate 
(\ref{s5orig}) three times to extend the sequence to $1,1,1,0,1,1,1,1$, then use the 
Somos-7 recurrence, which taken $\bmod\, 2$ gives 
$
S_{n+7}S_n \equiv S_{n+5}S_{n+2}+S_{n+4}S_{n+3}
$, 
and iterate it twice to extend the sequence to $1,1,1,0,1,1,1,1,1,0$, before finally applying  
the original Somos-5 relation once more to append another 1 at the end of this, so that after a 
total of 6 steps we have returned to the same initial values $ (1,1,1,0,1)$. Hence the period is 6, with a zero 
appearing precisely when $n\equiv 3\, (\bmod \,6)$, and the same pattern is repeated in 
$(T_n\bmod 2)$ except it is shifted back three steps. Taken $\bmod\,3$,  $(S_0,\ldots,S_4)$ 
provides well-balanced initial data, but we can proceed in the same way as for $p=2$. Working in $\F_3$, starting from 
$(1,1,1,2,0)$ we apply the Somos-5 recurrence three times to append the terms $2,2,1,2$ to the sequence, then 
apply the Somos-7 relation taken $\bmod\,3$, which gives 
$
S_{n+7}S_n \equiv 2S_{n+5}S_{n+2}+S_{n+4}S_{n+3} 
$, 
joining a 1 to the end of the sequence, and next we can apply Somos-5 seven more times before hitting the problem of division by zero, 
which adds the terms $2,2,0,2,1,1,1$, before another application of Somos-7 yields an extra 1, so that finally using Somos-5 again three more times 
we find that we have appended 
the 16 terms 
$$ 
2,2,1,2,1,2,2,0,2,1,1,1,1,1,2,0.
$$
The last five terms above are the initial values we started with, so the period is 16, and the zero terms appear  
when $n\equiv 4\, (\bmod \,8)$. For the sequence $(T_n\bmod\,3)$ the situation is almost identical, but the repeating 
pattern is 
$$ 
1,2,1,1,2,2,2,0,1,1,1,2,2,1,2,0, 
$$ 
with zero terms appearing when $n\equiv 0\, (\bmod \,8)$. Now if we consider any prime $p$, we have already noted that 
$(1,-1,1,1,-7)$ provides well-balanced initial data for $(T_n)$, and since $T_0=0$ and the sequence  $\bmod\,p$ is periodic 
by Theorem  \ref{periods}, it follows that $p$ is a divisor of infinitely  many terms. 
\end{prf}

\begin{remark} 
From the point of view of Theorem \ref{periods}, there is nothing special about the primes 2 and 3. However, 
a fuller understanding of the periods in Somos-5 sequences is reached from the connection with the 
underlying elliptic curve, and here it turns out that 2 and 3 (along with 17) are primes of bad reduction, so in this 
sense they are special. In particular, some additional  explanation for the values of the periods will be offered in the next section, 
in terms of the finite field dynamics of associated QRT maps, which we now introduce. 
\end{remark}

\begin{figure}
 \centering
\label{s5fig}
\epsfig{file=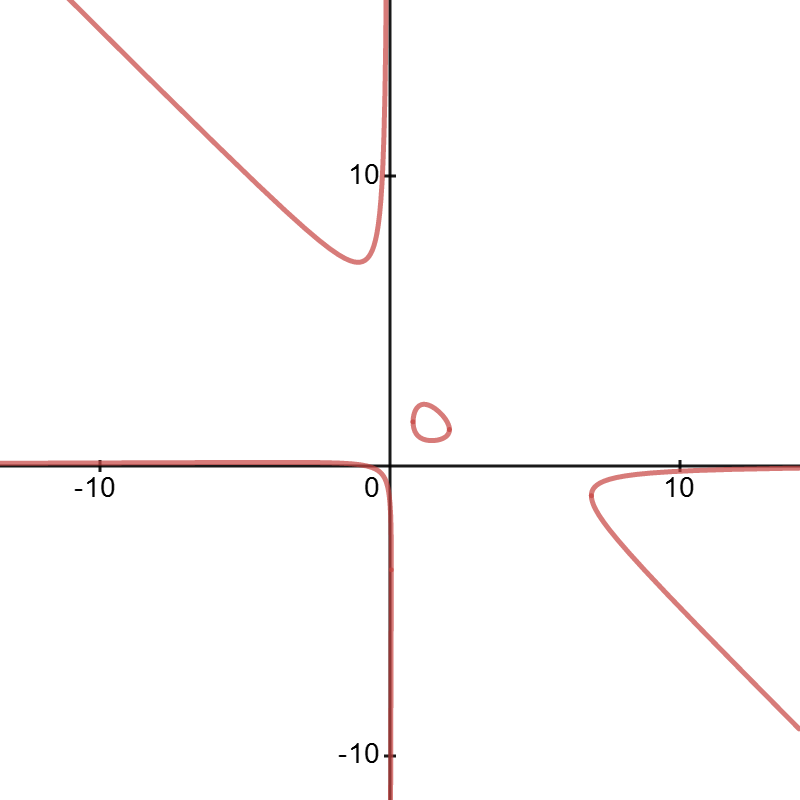, height=3.7in, width=3.7in}
\caption{The Somos-5 curve (\ref{s5curve})} 
\end{figure} 
 
\subsection{QRT maps from Somos-5 sequences} 

QRT maps, named after Quispel, Roberts and Thompson, are an 18-parameter family of birational maps of the plane that were introduced 
in   \cite{qrt} in order to unify various functional equations and maps appearing in statistical mechanics, dynamical systems and 
discrete soliton theory \cite{qrt2}. They are integrable maps in the sense of \cite{bruschi, maeda, veselov}, 
having an invariant symplectic form and a conserved quantity, and 
can be  defined intrinsically starting from families of plane biquadratic curves,  with associated elliptic fibrations of rational surfaces \cite{duistermaat, tsuda}.
The Somos-5 sequences (\ref{s5seq}) and (\ref{s5eds}) both generate particular orbits of the same QRT map, which is obtained by considering the ratios  
\beq \label{uvdef}
u_n=\frac{S_{n-2}S_{n+1}}{S_{n-1}S_n}, \qquad 
v_n  =\frac{T_{n-2}T_{n+1}}{T_{n-1}T_n}. 
\eeq 
These two sets of rational numbers both satisfy the same rational recurrence relation of order two, that is 
\beq\label{sqrtmap}
u_{n+1}u_{n-1}=1+\frac{1}{u_n}, \qquad v_{n+1}v_{n-1}=1+\frac{1}{v_n}, 
\eeq 
and the conserved quantity (\ref{Jt}) for Somos-5 can be rewritten in terms of these ratios, leading to a conserved quantity for the rational recurrence in the form 
\beq\label{s5qrtj}
\tilde{J}=u_n+u_{n-1}+\frac{1}{u_n}+\frac{1}{u_{n-1}}+\frac{1}{u_nu_{n-1}}, 
\eeq 
with the initial conditions $u_0=u_1=1$ giving $\tilde{J}=5$ in this case. (The other conserved quantity (\ref{It}) cannot be reduced to a function of these ratios.)
In other words, the sequence of points $(U,V)=(u_{n},u_{n+1})$ lies on the cubic plane curve 
\beq\label{s5curve}
{\cal C}: \qquad 
U^2V+UV^2-5UV+U+V+1=0, 
\eeq 
and the same is true for the sequence of points $(U,V)=(v_{n},v_{n+1})$; see Fig.2 
for a plot of the real curve in $\R^2$. The curve (\ref{s5curve}) is just (\ref{tJcurve}) 
with the particular values $\tal=\tbe=1$, $\tilde{J}=5$ for the coefficients: 
it is  symmetric and biquadratic, so it admits the simple involutions 
\beq\label{invols} 
\iota:  \, (U,V)\mapsto (V,U), \qquad
\iota_h: \, (U,V)\mapsto (U^\dagger,V) ,
\eeq 
where the horizontal switch $\iota_h$ is obtained by intersecting the curve with a horizontal line and  replacing each point $(U,V)$ with the other point of 
intersection  $(U^\dagger,V)$, which from Vieta's formula for the product of the roots of a quadratic is  given by 
$UU^\dagger=1+1/V$. Thus we see that the rational recurrence (\ref{sqrtmap}) corresponds to a symmetric QRT map, being given by the composition $\varphi=\iota\circ\iota_h$ 
of these two involutions, which sends $(u_{n-1},u_n)\mapsto (u_n,u_{n+1})$; and although so far it has been defined only on one particular curve, it 
lifts  to a birational map of the plane by taking the pencil of curves obtained by replacing $5\to \tilde{J}$ in  (\ref{s5curve}). 
Moreover, by construction each orbit of $\varphi$ lies on one of these curves, which generically has genus one,  and corresponds to a sequence 
of points $\hat{{\cal P}}_0+n{\cal P}$ (where $+$ denotes addition in the group law of the curve).

\begin{table}[h!]
  \begin{center}
    \caption{The first few rational numbers in the sequences (\ref{uvdef}) and (\ref{fdef}).}
    \label{uvtable}
    \begin{tabular}{ | r|| r| r| r| r | r| r| r| r| r| r|} %
\hline
      $n$ & 0  &  1 & 2 & 3 & 4 & 5 & 6 & 7 &8 & 9 \\
\hline 
\hline   
&&&&&&&&&& \\ 
$u_n$ &  $1$ & $1$ &$2$    & $\tfrac{3}{2}$ & $\tfrac{5}{6}$ & $\tfrac{22}{15}$ & $\tfrac{111}{55}$ & $\tfrac{415}{407}$& $\tfrac{3014}{3071}$ & $\tfrac{45029}{22742}$ \\
&&&&&&&&&& \\ 
    $v_n$ & $\infty$ &  $\infty$ & $0$ & $-1$ & $7$ &    $-\tfrac{8}{7}$ &  $\tfrac{1}{56}$ & $-\tfrac{399}{8}$ &  $\tfrac{3128}{57}$ & $-\tfrac{455}{22287} $\\
&&&&&&&&&& \\ 
      $f_n$ & $\infty$ & $1$ & $-1$  & $2$ & $3$ 
& $-\tfrac{5}{7}$ & $\tfrac{11}{8}$ & $-37$ & $-\tfrac{83}{57}  $ & $\tfrac{274}{391}  $
 \\
&&&&&&&&&& \\ 
\hline 
    \end{tabular}
  \end{center}
\end{table}

In what follows, an important role will be played by the ratio
\beq\label{fdef}
f_n=\frac{S_n}{T_n}, 
\eeq 
which turns out to lead to a QRT map on a different biquadratic curve, related to (\ref{s5curve}) by a 2-isogeny. 

\begin{propn}\label{fprop} 
The ratio of the two Somos-5 sequences, given by (\ref{fdef}), can be written 
as 
\beq\label{fform}
f_n =
  \begin{cases} 
-\Phi(n\ka)/\Phi(2\ka) , & \text{for }  n \, \,\mathrm{even} \\ 
\Phi(n\ka)/\Phi(\ka) , & \text{for }  n \, \,\mathrm{odd}  , 
 \end{cases}
\eeq 
where 
\beq\label{bafn} 
\Phi(z) = \frac{\si(z+\om) e^{-\eta z}}{\si(\om)\si(z)}, 
\eeq 
with $\eta=\zeta(\om)$ being the Weierstrass zeta function evaluated at the half-period $\om$ of the curve 
(\ref{origwei}) with invariants as in (\ref{s5params}). The sequence of ratios satisfies the recurrence 
\beq\label{fqrt}
f_{n+1}f_{n-1} = 
  \begin{cases}
    (1+f_n^2)/(2-f_n^2), & \text{for }  n \, \,\mathrm{even} \\
     (1+2f_n^2)/(1-f_n^2),  & \text{for }  n \,\,\mathrm{odd} , 
  \end{cases}
\eeq 
and for all $j\in\Z$ the points $(W,Z)=(f_{2j\pm1},f_{2j})$ lie on the biquadratic plane curve 
\beq\label{fcurve} 
(1-W^2)Z^2+3WZ+2W^2+1=0, 
\eeq 
corresponding to an orbit of a QRT map associated with this curve. 
\end{propn} 

\begin{figure} 
\label{fcurvepic}
\centering
\epsfig{file=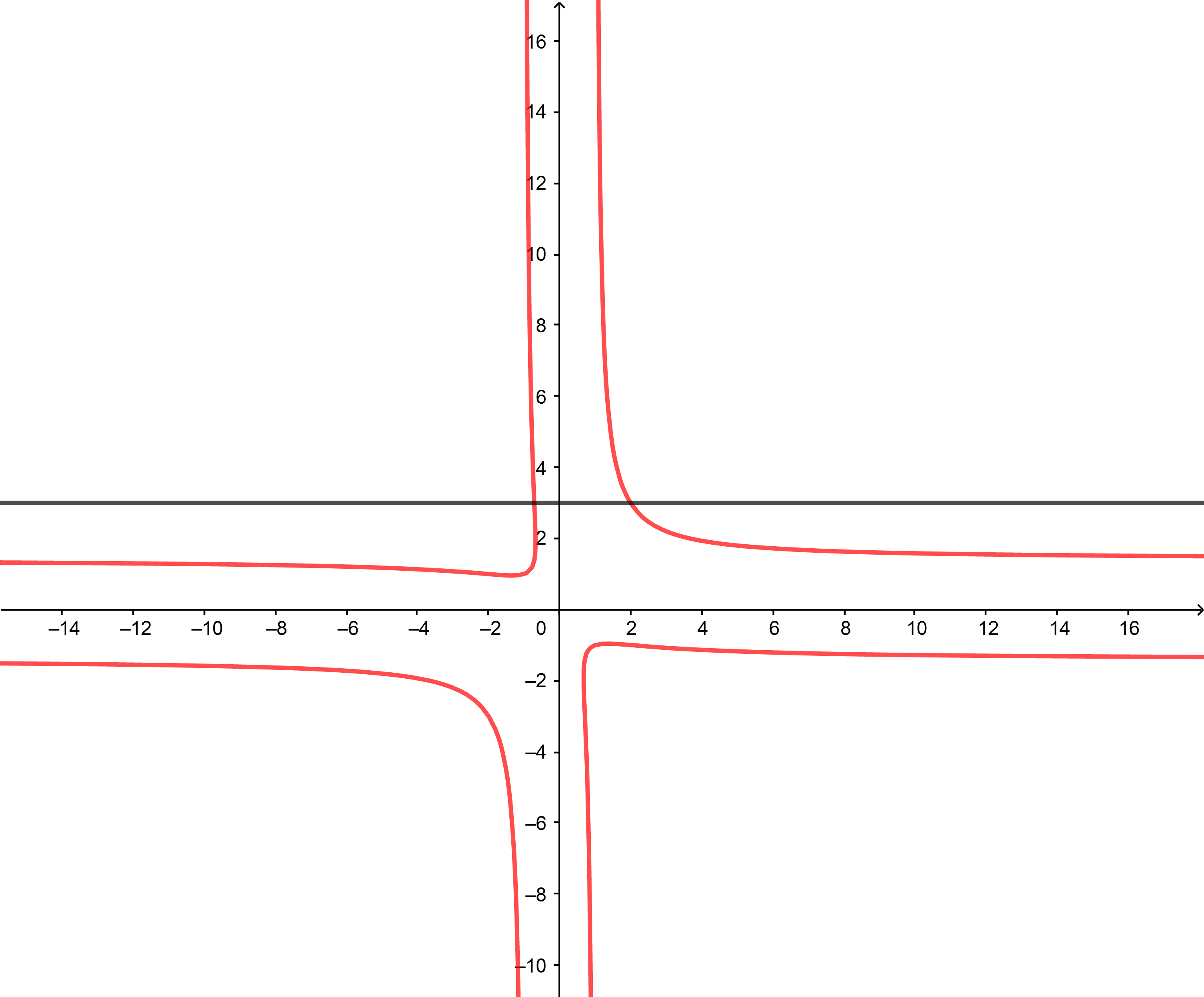, height=3.7in, width=3.7in}
\caption{The 
curve (\ref{fcurve}) in the $(W,Z)$ plane intersecting the horizontal line $Z=3$ in two points.}\end{figure}

\begin{prf}
We begin by recalling some of the properties of the function $\Phi(z)$, which  is well known as a solution 
of the simplest case of Lam\'e's equation, in the form of a Schr\"odinger equation with an elliptic potential, i.e.\ it satisfies the differential 
equation $\Phi''(z)-2\wp(z)\Phi(z)=\wp(\om)\Phi(z)$. From the quasiperiodicity of the sigma function it follows that $\Phi(z)$ is an odd function, and it is periodic 
with respect to the period $2\om$ but acquires a minus sign when shifted by the real/imaginary periods $2\om_1,2\om_2$; we record these properties, and 
its behaviour under a shift by $\om$, as follows:
\beq\label{phiprop}
\Phi(-z)=-\Phi(z), \quad \Phi(z+\om)=-\frac{e^{\eta\om}}{\si(\om)^2\Phi(z)}, \quad  \Phi(z+2\om)= \Phi(z), \quad \Phi(z+2\om_j)= -\Phi(z), \,\, \mathrm{for}\,\, j=1,2.
\eeq
Using the fact that $\Phi$ is odd, together with the standard identity 
\beq\label{sigmaddn}
\frac{\si(w+z)\si(w-z)}{\si(w)^2\si(z)^2}=\wp(z)-\wp(w),
\eeq 
valid for any $z,w\in\C$ (away from poles), 
it is apparent that 
\beq\label{phisq}
\Phi(z)^2=\wp(z)-\wp(\om),
\eeq
and since $g_2,g_3$ and $\wp(\om)$ as in (\ref{s5params}) and (\ref{2tors}) are all real,  $\Phi(z)\sim1/z$ as $z\to0$, and $\Phi$ has no real zeros, it follows 
that $\Phi(z)$ is real-valued for $z\in\R$, with $\Phi(z)>0$ for $0<z<2\om_1$ and $\Phi(z)<0$ for $-2\om_1<z<0$. Then since, from (\ref{kaval}), $\ka$ 
is negative and $-2\om_1<2\ka<0$, this allows us to compute 
\beq\label{phivals} 
\Phi(\ka)=-\sqrt{\wp(\ka)-\wp(\om)}=-\frac{\sqrt{2}}{6^{\frac{1}{4}}}, 
\quad 
\Phi(2\ka)=-\sqrt{\wp(2\ka)-\wp(\om)}=-\frac{1}{6^{\frac{1}{4}}},
\eeq 
where we have used the values of the $\wp$ function in (\ref{s5params})  and (\ref{2tors}), as well as the fact that 
$$
\wp(2\ka)=\wp(\ka)-\frac{\si(3\ka)}{\si(2\ka)^2\si(\ka)}=\wp(\ka)-\frac{\ra_3}{\ra_2^2}
=\tilde{\la}-\frac{\tal}{\tmu^2}=\frac{17}{12\sqrt{6}},
$$ 
which follows from (\ref{sigmaddn}) together with the appropriate expressions for the terms of the 
companion EDS $(\ra_n)$ defined by (\ref{cEDS}). Now for even $n=2k$ we can calculate the ratio (\ref{fdef}) using 
(\ref{Sn}) and (\ref{Tn}), to find 
$$ 
f_n=-\frac{B_+^k\,\mu^{k^2}\si(n\ka+\om)\si(2\ka)}{\hat{B}_-^{k-1}\,\mu^{(k-1)^2}\si(n\ka)\si(\om)}
=-\Phi(n\ka)\, \left(\frac{B_+}{\hat{B}_-}\, \mu^2e^{2\eta\ka}\right)^k \frac{\hat{B}_-\si(2\ka)}{\mu},
$$ 
and then note that from (\ref{bpm}) we may write 
$$
B_+=-\frac{\si(\ka)^4e^{-2\eta\ka}}{\si(2\ka)\Phi(2\ka)}=-\frac{e^{-2\eta\ka}}{\ra_2\Phi(2\ka)}=\frac{e^{-2\eta\ka}}{\tmu\Phi(2\ka)}
=-e^{-2\eta\ka},
$$
taking the value of $\Phi(2\ka)$ as in (\ref{phivals}), with $\tmu=6^{\frac{1}{4}}$, 
so by (\ref{bhatpm}) the expression in round brackets with exponent $k$ above is 
$B_+\mu^2e^{2\eta\ka}/\hat{B}_-
=1$, 
while the ratio appearing to the right of the round brackets is 
$\hat{B}_-\si(2\ka)/\mu=-\tmu=1/\Phi(2\ka)$,
thus indeed this yields $f_n=-\Phi(n\ka)/\Phi(2\ka)$ for even $n$, as required. Similarly, for odd $n=2k+1$ we find 
$$
f_n=\frac{B_-^k\si(n\ka+\om)\si(\ka)}{\hat{B}_+^k\si(n\ka)\si(\ka+\om)}=\frac{\Phi(n\ka)}{\Phi(\ka)}\left(\frac{B_-e^{2\eta\ka}}{\hat{B}_+}\right)^k, 
$$
but then by (\ref{constr}) and (\ref{bhatpm}) we see that $B_-e^{2\eta\ka}/\hat{B}_+=-\mu B_+e^{2\eta\ka}/\mu=1$, so this reduces to the 
required formula in (\ref{fform}) when $n$ is odd.  Then to obtain the recurrence (\ref{fqrt}), let us write 
$f_n=C_\pm \Phi(n\ka)$ for appropriate constants $C_\pm$ depending on the parity of $n$, as in (\ref{fform}), so that the left-hand side of 
the recurrence is given by $f_{n+1}f_{n-1}=C_\mp^2\Phi(n\ka+\ka)\Phi(n\ka-\ka)$ for even/odd $n$ respectively. Hence, by 
substituting the appropriate ratios of sigma functions and applying (\ref{sigmaddn}) to the numerator and denominator, 
this gives 
$$
f_{n+1}f_{n-1}=C_\mp^2 e^{-2\eta n\ka}\frac{\si(n\ka+\om)^2\big(\wp(\ka)-\wp(n\ka+\om)\big)}
{\si(n\ka)^2\big(\wp(\ka)-\wp(n\ka)\big)}.
$$ 
The function $\wp(z+\om)$ is an even elliptic function of order two with double poles at $z\equiv \om \bmod \Lambda$, 
where $\Lambda=2\om_1\Z\oplus 2\om_2\Z$ is the period lattice of the curve  (\ref{origwei})  with invariants 
$g_2,g_3$ as in (\ref{s5params}), hence it can be written in the form  $\wp(z+\om)=\big(A'\wp(z)+B'\big)/\big(\wp(z)-\wp(\om)\big)$ for 
suitable constants $A',B'$, and the leading (constant) term in the Taylor expansion at $z=0$ gives $A'=\wp(\om)$; the value of $B'$ can be fixed from the $O(z^2)$ term in the 
expansion, or by using the addition formula for the $\wp$ function, but will not be needed here. Thus, for another constant $B''$, we may write  
$$
f_{n+1}f_{n-1}=C_\mp^2 \Phi(n\ka)^2\frac{\big(\Phi(\ka)^2+B''\Phi(n\ka)^{-2}\big)}
{\big(\Phi(\ka)^2-\Phi(n\ka)^2\big)} 
=C_\mp^2 \frac{\big(\Phi(\ka)^2C_\pm^{-2}f_n^2+B''\big)}
{\big(\Phi(\ka)^2-C_\pm^{-2}f_n^2\big)}.
$$
Using the fact that $C_+=-\Phi(2\ka)^{-1}$, $C_-=\Phi(\ka)^{-1}$ and the values $f_1=-f_2=1$, $f_3=2$, we immediately find 
$B''=\frac{1}{2}\Phi(\ka)^4$ from the $n=2$ case of this relation, and this fixes the recurrence in the form (\ref{fqrt}) for even/odd $n$ respectively. 
If we start from the pencil of biquadratic curves $(1-W^2)Z^2+\hat{J}WZ+2W^2+1=0$, with arbitrary parameter $\hat{J}$, 
then the horizontal switch $\bar{\iota}_h:\, (W,Z)\mapsto (W^\dagger,Z)$ corresponds to the formula (\ref{fqrt}) for even $n$, which sends 
$(f_{2j-1},f_{2j})\mapsto (f_{2j+1},f_{2j})$, and the 
vertical switch $\bar{\iota}_v:\, (W,Z)\mapsto (W,Z^\dagger)$ corresponds to the case of odd $n$, which sends 
$(f_{2j+1},f_{2j})\mapsto (f_{2j+1},f_{2j+2})$, while the composition of these two involutions is a QRT map of general type, 
$\varphi_{QRT}=\bar{\iota}_v\circ\bar{\iota}_h$. The initial values $f_1=-f_2=1$ fix the value $\hat{J}=3$, giving an orbit that lies on the curve (\ref{fcurve}), which is illustrated in Fig.3, 
showing the horizontal line for the switch  $\bar{\iota}_h:\, (2,3)\mapsto \big(-\tfrac{5}{7},3\big)$ corresponding to $n=4$ in (\ref{fqrt}). 
\end{prf}


\begin{propn}\label{2isogeny} 
There is a 2-isogeny relating the  the curve (\ref{fcurve}) 
 to (\ref{s5curve}). 
\end{propn} 
\begin{prf} 
The 2-isogeny relating the curve (\ref{fcurve}) 
to (\ref{s5curve}), or equivalently to the corresponding Weierstrass curve (\ref{origwei}), 
can be seen in various ways. 
First of all, note that the function $\Phi(z)$ is not an elliptic function with respect to the original period lattice $\Lambda$, which is generated by the periods 
$2\om$ and $2\om_2$,  but it is elliptic with respect to the lattice generated by  $2\om$ and $4\om_2$; so if we set $\tau=\om_2/\om$ and start with the normalized 
lattice with generators $1,\tau$, then the new lattice has generators $1,2\tau$, and the overall effect is the period doubling $\tau\mapsto 2\tau$, corresponding to what is known 
classically as the Landen transformation (see chapter XXII in \cite{ww}, or \cite{maillard} and references). At the level of the curves, this can be seen by computing $\rho$,  the cross-ratio 
of the three roots of the cubic in 
(\ref{weier}) together with $\infty$, which allows the  
j-invariant to be calculated as 
\beq\label{badj} 
j=\frac{256(\rho^2-\rho+1)^3}{\rho^2(\rho-1)^2}=\frac{11^6}{612}.  
\eeq 
Then for the new curve related by the Landen transformation, we have that 
$$ 
\rho^*=\frac{4\sqrt{\rho}}{(1+\sqrt{\rho})^2}
$$ 
gives (an appropriate choice of) the value of the cross-ratio of the  roots of the quartic in the equation  
\beq\label{fquartic} 
\bar{y}^2=8W^4+5W^2-4, 
\eeq 
which is birationally equivalent to (\ref{fcurve}) 
via the transformation 
$$ 
Z=\frac{-3W+\bar{y}}{2(1-W^2)}, 
$$ 
from which one finds the j-invariant 
$$ 
j^*=\frac{256((\rho^*)^2-\rho^*+1)^3}{(\rho^*)^2(\rho^*-1)^2}=\frac{16(\rho^2+14\rho+1)^3}{\rho(\rho-1)^4}=\frac{46268279}{46818}.  
$$ 
Another indirect check is provided by verifying that the corresponding Hauptmoduls $\rx=1728/j$, $\ry=1728/j^*$ are 
the coordinates of a point on the modular curve 
$$ 
1953125\rx^3\ry^3-187500\rx^2\ry^2(\rx+\ry)+375\rx\ry·(16\rx^2-4027\rx\ry+16\ry^2)-64(\rx+\ry)(\rx^2+1487\rx\ry+\ry^2)+110592\rx\ry=0. 
$$ 
These checks all confirm that there is a 2-isogeny over $\C$, but to establish this over $\Q$ we provide a direct transformation of coordinates, given 
by the formulae 
\beq\label{2isog} 
x^*=W^2,  \qquad y^*=(1-W^2)WZ+\frac{3}{2}W^2,
\eeq 
which transforms (\ref{fcurve}) to the cubic 
$$ 
4(y^*)^2=8(x^*)^3 +5(x^*)^2-4x^*,
$$ 
and after shifting $x^*\to x^* -5/24$ and rescaling $x^*$ and $y^*$ by suitable powers of 2 this is seen to be equivalent to (\ref{weier}).  
\end{prf}
\begin{remark}
The primes appearing in the factorization of the denominator of (\ref{badj}), $612=2^2\cdot3^2\cdot 17$, 
are the primes of bad reduction, which we will return to at the end of the next section. 
\end{remark} 

\section{Proof of main results} 

\setcounter{equation}{0}

In order to prove our main theorem, we must first show that there is an infinite sequence of 
pairs of positive points $(M_a,P_a,X_a)$, $(M_b,P_b,X_b)$ lying on the Schubert surface (\ref{schub}), so that 
\beq\label{schubap}
M_a-{M_a}^{-1}=P_a-{P_a}^{-1}+
2\left( X_a-{X_a}^{-1}\right), 
\eeq  \beq\label{schubbp}
M_b-{M_b}^{-1}=P_b-{P_b}^{-1}+
2\left( X_b-{X_b}^{-1}\right), 
\eeq  
and these sets of Schubert parameters are compatible in the sense that 
\beq\label{aovercab} 
\frac{2(X_a+X_a^{-1})}{P_a+P_a^{-1}}= 
\frac{P_b+P_b^{-1}}{M_b+M_b^{-1}}=\frac{a}{c},  
\eeq 
\beq\label{bovercab}  
\frac{M_a+M_a^{-1}}{P_a+P_a^{-1}} = 
\frac{2(X_b+X_b^{-1})}{M_b+M_b^{-1}}=\frac{b}{c},  
\eeq 
which serves as the definition of the ratios $a/c$, $b/c$, and hence (up to scale) defines the associated Heron triangle with two rational medians.  
As the main initial step, we begin by giving a proof of the empirical observation of Buchholz and Rathbun in \cite{br1} 
that there is an infinite sequence of Schubert 
parameters with signs, given in terms of the two Somos-5 sequences by (\ref{schubs5a}) and (\ref{schubs5b}), or equivalently in terms of 
the rational sequences $(u_n)$, $(v_n)$ and $(f_n)$ by 
\beq\label{mpxfuva}
M_{a,n}=-\frac{f_{n+1}f_{n+2}^2}{f_n} , 
\quad 
P_{a,n}=-\frac{1}{u_{n+2}v_{n+2}}, 
\quad 
X_{a,n}=2^{(-1)^{n+1}}\frac{f_{n}f_{n+2}^2}{f_{n+3}},  
\eeq 
\beq\label{mpxfuvb}
M_{b,n}={u_{n+3}v_{n+3}},  
\quad 
P_{b,n}=-\frac{f_{n+2}^2f_{n+3}}{f_{n+4}} , 
\quad 
X_{b,n}=2^{(-1)^{n}}\frac{f_{n+1}}{f_{n+2}^2f_{n+4}}.  
\eeq 
Having proved that these formulae give points on the Schubert surface satisfying the necessary constraints, we will then show that the pattern of signs varies 
coherently with $n$  in such a way that replacing any negative parameter by minus its reciprocal will preserve the constraints and hence provide, for each $n$, 
two compatible sets  of positive Schubert parameters. Then finally we will be able to verify the formulae for $a,b,c,k,\ell$ and $\Delta$ in Theorem \ref{main}.    

Before tackling the Schubert parameters with signs, we introduce the sequence of quantities 
\beq\label{bars}
\begin{array}{rl} 
\bs_n 
=-S_{n+3}S_{n+4}^2T_n^2T_{n+1}, 
& 
\ba_n 
=-S_{n+2}T_{n+1}T_{n+2}^3T_{n+3}, \\ 
\bb_n 
=S_{n+1}S_{n+2}^3S_{n+3}T_{n+2}, 
& 
\bc_n 
=S_n^2S_{n+1}T_{n+3}T_{n+4}^2,
\end{array} 
\eeq 
which we will refer to as the signed lengths. After taking absolute values, for each $n$ there is an equality of sets of positive 
numbers: 
$$ 
\{ \, |\bs_n|, |\ba_n|, |\bb_n|, |\bc_n|\, \} = \{ \, s, s-a, s-b, s-c\, \}.
$$ 
When $n=1$ the choice of signs in (\ref{bars}) ensures that %
$\bs_1,\ba_1,\bb_1$ and $\bc_1$ 
are all positive, and coincide with $s,s-a,s-b$ and $s-c$, respectively (see Table \ref{lengtable}), but it turns out that  in general  
the signed lengths  correspond to a permutation of the latter four quantities with signs, in  a pattern that repeats with period 14. 
Nevertheless, the pattern of permutations and signs respects the linear relation 
$$s=(s-a)+(s-b)+(s-c)$$ that holds between their associated positive counterparts. 

\begin{table}[h!]
  \begin{center}
    \caption{The first few signed lengths.}
    \label{lengtable}
    \begin{tabular}{ | r|| r| r| r| r | r| r| r| r| } %
\hline
      $n$ & 0  &  1 & 2 & 3 & 4 & 5 & 6 & 7  \\
\hline 
\hline   
&&&&&&&& \\ 
$\bs_n$ &  $0$ & $75$ &$-605$    & $-15059$ & 
$1784251$ & $-2442672736$ & $25972376704$ & $2633103021849$
 \\
&&&&&&&& \\ 
    $\ba_n$ & $1$ &  $2$ & $21$ & $13720$ & $-39424$ &    $-16872$ &  
$-6010068429$ & $-424782960455490$ 
\\
&&&&&&&& \\ 
      $\bb_n$ & $-2$ & $24$ & $270$  & $-28875$ & $1969880$ 
& $-46246189$ & $-330416266542$ & $812450885698024$ 
 \\
&&&&&&&& \\ 
      $\bc_n$ & $1$ & $49$ & $-896$  & $96$ & $-146205$ 
& $-2396409675$ & $362398711675$ & $-385034822220685$ 
 \\
&&&&&&&& \\ 
\hline 
    \end{tabular}
  \end{center}
\end{table}

\begin{lem}\label{linlem}
The signed lengths (\ref{bars}) satisfy 
\beq\label{lin} 
\bs_n=\ba_n+\bb_n+\bc_n
\eeq  
for all $n\in\Z$. 
\end{lem}

\begin{prf}
The equation (\ref{lin}) is a degree 6 linear relation between products of terms of the Somos sequences $(S_n)$ and $(T_n)$, 
but has a different form compared with the identities for minors proved in Proposition \ref{detprop}. 
Upon rearranging 
and dividing both sides by $\ba_n$,  
it can be rewritten as 
\beq\label{vfid} 
1-f_{n+1}f_{n+2}^2f_{n+3} =(v_{n+2}v_{n+3})^2\left( \frac{f_n^2f_{n+1}+f_{n+3}f_{n+4}^2}{f_{n+2}} \right).
\eeq 
The latter identity relates the rational sequences $(v_n)$ and $(f_n)$ corresponding to particular orbits  of the two different QRT maps 
(\ref{sqrtmap}) and (\ref{fqrt}); alternatively, it could be written in terms of  $(u_n)$ and $(f_n)$, since the 
ratio $u_n/v_n$ can be expressed in terms of the $f_n$. It is equivalent to an identity 
between elliptic functions,  since $f_n$ is given by the formula (\ref{fform}), 
and $v_n$ is given by a ratio of Somos-5 terms as in (\ref{uvdef}), which are themselves given in analytic form by 
(\ref{Tn}).  
Hence, in order to prove it, we  
set $z=(n+2)\ka$, so that the left-hand side above can be written as 
\beq\label{lhs}
1-\frac{\Phi(z-\ka) \Phi(z+\ka)\big(\wp(z)-\wp(\om)\big)}{\Phi(\ka)^2\Phi(2\ka)^2}, 
\eeq
while the right-hand side is given by 
\beq\label{rhs}
\frac{\tmu^4\big( \wp(z)-\wp(2\ka)\big)^2\, \Phi(z+\om)\,
\Big( \Phi(z-\ka)\big(\wp(z-2\ka)-\wp(\om)\big) +  \Phi(z+\ka)\big(\wp(z+2\ka)-\wp(\om)\big)\Big)}
{e^{\eta\om}\si(\om)^{-2}\Phi(\ka)\Phi(2\ka)} ,
\eeq 
where we have used (\ref{sigmaddn}) to simplify the ratio of sigma functions in $v_{n+2}v_{n+3}$, 
together with (\ref{phisq}) and the expression for the reciprocal of $\Phi$ in (\ref{phiprop}).  
Both (\ref{lhs}) and (\ref{rhs}) are even elliptic functions of $z$, so to verify their equality it is sufficient 
to check that they have poles in  the same places with the same singular part of the Laurent expansion around each pole, and agree at one finite value, 
since their difference is then an elliptic function without poles and therefore constant, and if they take the same finite value somewhere then this constant must be $0$.
The left-hand side has double poles for $z\equiv 0 \bmod\Lambda$, and for $z\to 0$ we have 
$\Phi(z-\ka)\Phi(z+\ka) =-\Phi(\ka)^2 +O(z^2)$ and $\wp(z)=1/z^2+O(z^2) $, so its Laurent 
expansion around the origin is $\Phi(2\ka)^{-2}/z^2+O(1)$. On the right-hand side, note that 
$$\Phi(z\pm \ka)\big(\wp(z\pm2\ka)-\wp(\om)\big) = \pm\Phi(\ka)\big(\wp(2\ka)-\wp(\om)\big) 
+C'z
+O(z^2), 
$$ 
where $C'=\big(\Phi'(\ka)\big(\wp(2\ka)-\wp(\om)\big) +\Phi(\ka)\wp'(2\ka) \big)$, 
so the sum of these two terms 
gives an odd function with Taylor expansion $2C'z +O(z^3)$ at the origin, while $\Phi(\om)=0$ and 
we have $\Phi(z+\om)=\Phi'(\om)z+O(z^3)$. Hence to leading order, the expansion of the right-hand side around $z=0$ 
is $C/z^2+O(1)$, where 
$$
C=\frac{2e^{-\eta\om}\si(\om)^2\tmu^4\Big(\Phi'(\ka)\big(\wp(2\ka)-\wp(\om)\big) +\Phi(\ka)\wp'(2\ka) \Big)\Phi'(\om)}{\Phi(\ka)\Phi(2\ka)}.
$$ 
Now from (\ref{phisq}) it follows that $2\Phi\Phi'=\wp'$, so $\Phi'(\ka)=\frac{1}{2}\wp'(\ka)/\Phi(\ka)$, 
while putting $\Phi(z)=1/z+O(z)$ as $z\to0$ into the identity $\Phi(z)\Phi(z+\om)=-e^{\eta\om}\si(\om)^{-2}$ gives 
$\Phi'(\om)=-e^{\eta\om}\si(\om)^{-2}$. Thus we have 
$$ 
C\Phi(2\ka)^2=-2\Phi(2\ka)\tmu^4\Big(\frac{1}{2}\Phi(\ka)^{-2}\wp'(\ka)\big(\wp(2\ka)-\wp(\om)\big)+\wp'(2\ka)\Big)
=\tmu^3\left(\tmu\Phi(2\ka)^2\Phi(\ka)^{-2}-2\tmu^{-3}\right), 
$$
using (\ref{phisq}) to replace the $\big(\wp(2\ka)-\wp(\om)\big)$ term inside the large brackets, together with 
$\wp'(\ka)=\tmu$ and $\wp'(2\ka)=\si(4\ka)/\si(2\ka)^4=\ra_4/\ra_2^4=-\tbe/\tmu^3=-\tmu^{-3}$ from the 
formulae for the companion EDS, with $\tbe=1$ in this case. Then substituting in $\tmu=6^{\frac{1}{4}}$ and the 
values of $\Phi$ in (\ref{phivals})  gives 
$$ 
C\Phi(2\ka)^2=6\, \left(\frac{1/\sqrt{6}}{2/\sqrt{6}}\right)-2 =1, 
$$ 
so the singular parts of the Laurent expansions as $z\to 0$ are the same on each side of (\ref{vfid}). 
Then (\ref{lhs}) has poles at precisely two other places, namely simple poles for $z\equiv \pm \ka\bmod\Lambda$ with 
residues 
$$\mp\frac{\Phi(2\ka)  \big(\wp(\ka)-\wp(\om)\big) }{\Phi(\ka)^2 \Phi(2\ka)^2 }=\mp\Phi(2\ka)^{-1}, 
$$
respectively, while in (\ref{rhs}) we see simple poles  at the same places, with the residues being  
$$
\frac{e^{-\eta\om}\tmu^4\big( \wp(\ka)-\wp(2\ka)\big)^2\, \Phi(\pm\ka+\om)\, \big(\wp(\ka)-\wp(\om)\big) }
{\Phi(\ka)\Phi(2\ka)}=\mp \frac{\tmu^4\big( \wp(\ka)-\wp(2\ka)\big)^2\Phi(\ka)^2}{\Phi(\ka)^2\Phi(2\ka)}
=\mp\frac{6\left(\frac{29}{12\sqrt{6}} -\frac{17}{12\sqrt{6}} \right)^2}{\Phi(2\ka)},
$$
i.e.\ equal to $=\mp\Phi(2\ka)^{-1}$, the same as for (\ref{lhs}). The expression (\ref{rhs}) also contains the terms 
$\wp(z\mp2\ka)$ with double poles at $z\equiv \pm2\ka\mod\Lambda$, but these poles are cancelled by the prefactor 
$\big( \wp(z)-\wp(2\ka)\big)^2$ which has double zeros at these points; and there is also the term $\Phi(z+\om)$ 
which gives  a simple pole for $z\equiv \om\mod\Lambda$, but the factor that appears after it, evaluated at $z=\om$, 
yields 
$$ 
 \Phi(\om-\ka)\big(\wp(\om-2\ka)-\wp(\om)\big) +  \Phi(\om+\ka)\big(\wp(\om+2\ka)-\wp(\om)\big)
=\big(\Phi(\ka+\om)-\Phi(\ka-\om)\big)\big(\wp(2\ka+\om)-\wp(\om)\big)=0, 
$$
from the periodicity of $\Phi$ and $\wp$ under shifts by $2\om$ and the fact that these are odd/even functions, respectively, 
so this simple pole is cancelled by a zero. The values $n=-4,-3,-2,-1,0$ are all singular cases of (\ref{vfid}), due to the presence of the term $f_0=\infty$: these give the values of $z$  corresponding to the poles, together 
with the removable singularities at $z=\pm2\ka$; but it is easy to check directly that (\ref{lin}) is satisfied for these values of $n$. The first nonsingular value is  $n=1$, which  corresponds to setting  $z=3\ka$, and 
from the values in Table \ref{uvtable} it is 
straightforward to check that the left-hand side and the right-hand side of   (\ref{vfid}) are  both  equal to $13$ in this case. 
Hence the functions (\ref{lhs}) and (\ref{rhs}) coincide, and the result follows.
\end{prf} 

It appears that we have to prove four identities for the two sets of Schubert parameters with signs: two copies of the equation for the Schubert surface, and two constraints between the two sets of parameters. Moreover, 
(\ref{aovercab}) and (\ref{bovercab}) each contain two equalities, so upon replacing the ratios $a/c$ and $b/c$ by appropriate combinations of signed lengths, this gives a further two identities that must be verified. 
However, there is a symmetry to the problem which cuts the amount of work down by half. 

\begin{lem}\label{symlem} 
Under the involution $n\to -n-4$, the Schubert parameters with signs transform as 
\beq\label{sym} 
 M_{a,-n-4}=P_{b,n}, \qquad P_{a,-n-4}=-1/M_{b,n}, \qquad X_{a,-n-4}=1/X_{b,n},
\eeq 
and the signed lengths transform as 
\beq\label{signsym}
\bar{s}_{-n-4}=\bar{c}_n, \quad \bar{a}_{-n-4}=-\bar{a}_n, \quad \bar{b}_{-n-4}=-\bar{b}_n, \quad \bar{c}_{-n-4}=\bar{s}_n.
\eeq 
\end{lem}
\begin{prf}
As already noted previously, the sequences (\ref{s5seq}) and (\ref{s5eds}) extend to all $n\in\Z$ in a way that is respectively symmetric/antisymmetric about $n=0$, so that 
$$
S_{-n}=S_n, \qquad T_{-n}=-T_n. 
$$
For the corresponding rational sequences defined by (\ref{uvdef}) and (\ref{fdef}), this implies immediately that 
$$
u_{1-n}=u_n, \qquad v_{1-n}=v_n, \qquad f_{-n}=-f_n,
$$
and then for  $n\to -n-4$ it follows from  (\ref{mpxfuva}) and (\ref{mpxfuvb}) that the Schubert parameters transform 
according to (\ref{sym}).
For the signed lengths, we have 
$$
\bar{s}_{-n-4}=-S_{-n-1}S_{-n}^2T_{-n-4}^2T_{-n-3}=- S_{n+1}S_{n}^2(-T_{n+4})^2(-T_{n+3})=\bar{c}_n,
$$
and similarly for the other three. 
\end{prf}

Now starting from the Schubert equation (\ref{schubap}) for $(M_a,P_a,X_a)$ and replacing $M_a\to P_b$, $P_a\to-M_b^{-1}$, $X_a\to X_b^{-1}$ gives 
$
P_b-P_b^{-1}=-M_b^{-1}+M_b+2(X_b^{-1}-X_b)$, 
which is just a rearrangement of  (\ref{schubbp}), and if at the same time we replace $P_b\to M_a$, $M_b\to -P_a^{-1}$, $X_b\to X_a^{-1}$ then it is 
clear that these two copies of Schubert's equation are switched.  
Similarly, applying this involution to the first equality in (\ref{aovercab}) gives the first equality in (\ref{bovercab}), up to an overall minus sign, and vice versa. 
Furthermore, we can apply 
this symmetry to the second equality   in (\ref{aovercab}) by interpreting the right-hand side suitably in terms of the signed lengths, so that (omitting the index $n$) we 
may write 
\beq\label{signida} 
\frac{P_b+P_b^{-1}}{M_b+M_b^{-1}}=\frac{\bar{s}-\bar{a}}{\bar{s}-\bar{c}},  
\eeq 
and then applying the involution to the Schubert parameters on the left, as well as $\bar{s}\to \bar{c}$, $\bar{a}\to-\bar{a}$, $\bar{c}\to\bar{s}$ on the right, 
 this becomes 
$$
-\frac{M_a+M_a^{-1}}{P_a+P_a^{-1}}=\frac{\bar{c}+\bar{a}}{\bar{c}-\bar{s}},  
$$
but then applying the first equality in (\ref{bovercab})  together with Lemma \ref{linlem}, this implies 
\beq\label{signidb} 
\frac{2(X_b+X_b^{-1})}{M_b+M_b^{-1}}=\frac{\bar{s}-\bar{b}}{\bar{s}-\bar{c}}, 
\eeq 
which is just the second equality in (\ref{bovercab}), with the right-hand side written in terms of the signed lengths. Hence we see that for the  the Schubert parameters with signs 
and the signed lengths, the involution 
$n\to -n-4$ interchanges the two copies of Schubert's equation, and the two pairs of equalities given by the constraints   (\ref{aovercab}) and (\ref{bovercab}), so it is equivalent to 
switching $a\leftrightarrow b$ in each triangle. Thus it is sufficient to prove  (\ref{schubap}) and the two equalities in (\ref{aovercab}) for all integer values of $n$, and the other 
relations follow by symmetry.

\begin{thm}\label{schubrels} The two sets of Schubert parameters with signs, that is $(M_{a,n},P_{a,n},X_{a,n})$, $(M_{b,n},P_{b,n},X_{b,n})$ given by 
(\ref{mpxfuva}) and (\ref{mpxfuvb}), satisfy the relations (\ref{schubap}), (\ref{schubbp}),  (\ref{aovercab}) and (\ref{bovercab}) for all $n\in\Z$, where 
the quantities on the right-hand sides of (\ref{aovercab}) and (\ref{bovercab}) should be interpreted in terms of the signed lengths as 
$$ 
\frac{\bar{s}_n-\bar{a}_n}{\bar{s}_n-\bar{c}_n} \qquad and \qquad \frac{\bar{s}_n-\bar{b}_n}{\bar{s}_n-\bar{c}_n}, 
$$ 
respectively.
\end{thm} 
\begin{prf}
From (\ref{mpxfuva}) we may set $z=(n+2)\ka$ and write 
\beq\label{maform}
\bear{l}
M_{a,n} = -C_+C_-\tgam \Phi(z-\ka)\Phi(z)^2\Phi(z-2\ka+\om), \\
P_{a,n} =  -\mu^{-2} \si(z-\ka+\om)\si(z+\om) \si(z-\ka)\si(z) /\big( \si(z-2\ka+\om)\si(z+\ka+\om) \si(z-2\ka)\si(z+\ka)\big) 
, \\
X_{a,n}  = 2^{\mp 1} C_{\pm}^3C_{\mp}^{-1}\tgam\Phi(z-2\ka)\Phi(z)^2\Phi(z+\ka+\om), 
\eear 
\eeq 
for even/odd $n$, respectively, where we have used (\ref{phiprop}) with the same notation as in the proof of
Proposition (\ref{fprop}), 
and for convenience we have written the residue of $\Phi(z)^{-1}$ at $z=-\om$ as  
\beq\label{tgamdef}
\tgam=-e^{-\eta\om}\si(\om)^2=\underset{ z\to -\om}{\rm lim}\frac{(z+\om)}{\Phi(z)},
\eeq  
namely the multiplier that appears when the reciprocal of $\Phi$ is replaced by the same function shifted by $\om$.
The above formula for $P_{a,n}$ is a consequence of (\ref{Sn}) and (\ref{Tn}), 
which imply that 
$$ 
u_n = \left(\frac{B_\mp}{B_{\pm}}\right) \mu^{1\mp 1}\, \frac{\si(z-2\ka+\om)\si(z+\ka+\om)}{\si(z-\ka+\om)\si(z+\om)}, \quad 
v_n = \left(\frac{\hat{B}_\pm}{\hat{B}_{\mp}}\right) \mu^{1\pm 1}\, \frac{\si(z-2\ka)\si(z+\ka)}{\si(z-\ka)\si(z)}, 
$$ 
for even/odd $n$, and the prefactors in brackets cancel in the product $u_nv_n$, due to (\ref{constr}) and (\ref{bhatpm}). 
From these expressions, we see that $M_{a,n}$ has poles at $z\equiv 0,\ka,2\ka+\om$ modulo the period lattice $\Lambda$, while 
$M_{a,n}^{-1}$ has poles at $z\equiv \om,\ka+\om,2\ka$,
$P_{a,n}$ has poles at $z\equiv -\ka,-\ka+\om, 2\ka, 2\ka+\om$, while 
$P_{a,n}^{-1}$ has poles at $z\equiv 0,\om,\ka,\ka+\om$,
and 
$X_{a,n}$ has poles at $z\equiv 0,2\ka,2,-\ka+\om$, while 
$X_{a,n}^{-1}$ has poles at $z\equiv \om,2\ka+\om,-\ka$. Then to verify that $(M_{a,n},P_{a,n},X_{a,n})$ satisfies 
(\ref{schubap}) for all $n$, it is sufficient to check that the two sides of the equation define the same elliptic functions of $z$, by checking that they  agree in the 
singular parts of their Laurent expansions around the poles at all these places, and at one finite value. This is equivalent to saying that the 
formulae (\ref{maform}) define an analytic embedding of the elliptic curve in the Schubert surface, and this does not depend on the parity of $n$ because the 
expressions for   $M_{a,n}$ and $P_{a,n}$ are manifestly the same for even/odd $n$, while for the coefficient in front of the $z$-dependent part of the 
formula for $X_{a,n}$ we find the identity 
$\frac{1}{2}C_+^3C_-^{-1}=2C_-^3C_+^{-1}$,   
that is equivalent to $(C_+/C_-)^4=\big(\Phi(\ka)/\Phi(2\ka)\big)^4=4$, which  follows from (\ref{phivals}). To show that these 
parameter triples 
satisfy  (\ref{schubap}), it is convenient to rewrite the equation, collecting the terms as 
\beq\label{schubr}
(M_{a,n}-2X_{a,n})+(2X_{a,n}^{-1}-M_{a,n}^{-1})=P_{a,n}-P_{a,n}^{-1}. 
\eeq 
The first bracketed expression on the left-hand side above is given as a function of $z$ by 
$$ 
F(z)=\Phi(\ka)^{-1}\Phi(2\ka)^{-1}\Phi(z)^2G(z), 
$$ 
where 
$$
G(z)=\tgam \big( \Phi(z-\ka)\Phi(z-2\ka+\om)+2\Phi(z-2\ka)\Phi(z+\ka+\om)\big)=
\frac{\Phi(z-\ka)}{\Phi(z-2\ka)} +2\frac{\Phi(z-2\ka)}{\Phi(z+\ka)}
,
$$ 
and around $z=0$ we have $\Phi(z)^2=1/z^2+O(1)$ and 
$$ 
G(z) =\left(\frac{\Phi(-\ka)}{\Phi(-2\ka)} +2\frac{\Phi(-2\ka)}{\Phi(\ka)}\right)
+\left(\frac{\Phi'(-\ka)}{\Phi(-2\ka)}-\frac{\Phi(-\ka)\Phi'(-2\ka)}{\Phi(-2\ka)^2} +2\frac{\Phi'(-2\ka)}{\Phi(\ka)}-2\frac{\Phi(-2\ka)\Phi'(\ka)}{\Phi(\ka)^2}\right)z+O(z^2),
$$ 
so from  (\ref{phivals}) we see that $G(0)=0$, and as $z\to 0$, $G(z)=C^* z+O(z^2)$, 
with 
$$ 
C^*=\Phi'(\ka)\left(-\Phi(2\ka)^{-1}+2\Phi(2\ka)\Phi(\ka)^{-2}\right) 
+\Phi'(2\ka)\left(2\Phi(\ka)^{-1}+\Phi(\ka)\Phi(2\ka)^{-2}\right) =-\sqrt{2}\tmu^{-1}, 
$$
using $\Phi'(2\ka)=\frac{1}{2}\wp'(2\ka)/\Phi(2\ka)$ and the values of $\Phi(\ka)$, $\Phi(2\ka)$ and $\wp'(2\ka)$ as before. 
Hence $F(z)$ has a simple pole at $z=0$ with residue $-\Phi(\ka)^{-1}\Phi(2\ka)^{-1}\sqrt{2}\tmu^{-1}=-\tmu=-6^{\frac{1}{4}}$, while 
on the right-hand side of (\ref{schubr}), $-P_{a,n}^{-1}$ also has a simple pole there with residue 
$$
\mu^2\frac{\si(-2\ka+\om)\si(\ka+\om)\si(-2\ka)\si(\ka)}{\si(-\ka+\om)\si(\om)\si(-\ka)},
$$ 
but using  the  oddness of the sigma  function and its quasiperiodicity, e.g. $\si(\ka+\om)=-e^{2\eta\ka}\si(\ka-\om)$, we may rewrite 
this residue as $\mu^2\si(2\ka)^2\Phi(2\ka)=-\tmu$ by (\ref{mu}), so these two residues agree. 
Now a calculation of the effect of shifting the argument in the Schubert parameters by $\om$ shows that, because the expressions for $M_a$ and $X_a$ are both quartic 
in $\Phi$ with prefactors $\mp C_+C_-\tgam$, respectively, they satisfy 
$M_{a,n}(z+\om)=(C_+C_-\tgam)^2\tgam^{-4}M_{a,n}(z)^{-1}=\Phi(\ka)^{-2}\Phi(2\ka)^{-2}\tgam^{-2}M_{a,n}(z)^{-1}$, and similarly 
$X_{a,n}(z+\om)=\Phi(\ka)^{-2}\Phi(2\ka)^{-2}\tgam^{-2}X_{a,n}(z)^{-1}$. However, by (\ref{phisq}) we see that
$$
\tgam^2= \underset{ z\to -\om}{\rm lim}\frac{(z+\om)^2}{\wp(z)-\wp(\om)}=\frac{2}{\wp''(\om)}=\frac{2}{6\wp(\om)^2-\frac{1}{2}g_2}=-\frac{1}{2}\tmu^4=-3, 
$$
and the quasiperiodicity of the sigma function implies that $P_a$ is unchanged under shifting by this half-period. So overall, using  (\ref{phivals}) once more, we find  
\beq\label{half} 
M_{a,n}(z)^{-1}=-M_{a,n}(z+\om), \qquad P_{a,n}(z)=P_{a,n}(z+\om), \qquad X_{a,n}(z)^{-1}=-X_{a,n}(z+\om).  
\eeq 
In particular, if we consider simple poles at $z=-\om$, this implies that the residue of the second set of bracketed terms on the left-hand side of (\ref{schubr}) is also equal to $-\tmu$, and is 
the same as the residue of the term   $-P_{a,n}^{-1}$ on the right-hand side. 
Next we consider $z=\ka$, and verify that both  $M_{a,n}$ on the  left-hand side  and  $-P_{a,n}^{-1}$ on the right-hand side of (\ref{schubr}) have the same residue $\tmu$ there. 
Also, at $z=-\ka-\om$, we see that the residue at the simple pole in $-2X_{a,n}$ is given by 
$$ 
-2C_+C_-\tgam\Phi(-3\ka-\om)\Phi(-\ka-\om)^2=-2\tgam^{-2}\Phi(\ka)^{-3}\Phi(2\ka)^{-1}\Phi(3\ka)^{-1}=\frac{4}{f_3}\tmu^{-4}\Phi(\ka)^{-4}\Phi(2\ka)^{-1}=-\frac{\tmu}{2},
$$
and this is balanced by  $\mathrm{res}\, P_{a,n}|_{z=-\ka-\om}$, for which we find the same value 
$$ 
\mu^{-2}
\frac{\si(\ka)\si(2\ka)\si(\ka+\om)\si(2\ka+\om)}{\si(3\ka)\si(\om)\si(3\ka+\om)}
=\mu^{-2}\frac{\si(\ka)^2\si(2\ka)^2\Phi(\ka)\Phi(2\ka)}{\si(3\ka)^2\Phi(3\ka)}=\frac{\ra_2^2\Phi(2\ka)}{\ra_3^2f_3}=-\frac{\tmu}{2}.
$$
Using (\ref{half}), we see that the total residue at the simple pole at $z=2\ka$ on the left-hand side of (\ref{schubr}) comes from the combination 
$-M_{a,n}^{-1}-2X_{a,n}$, being given by 
$$ 
-C_+C_-\tgam \Phi(\ka+\om)\Phi(2\ka+\om)^2-2C_+C_-\tgam \Phi(2\ka)^2\Phi(3\ka+\om)=\tgam^{-2}\Phi(\ka)^{-2}\Phi(2\ka)^{-3}+2\Phi(\ka)^{-2}\Phi(2\ka)f_3^{-1}=\frac{\tmu}{2},
$$ 
and we find the same value for $\mathrm{res}\, P_{a,n}|_{z=2\ka}$ on the right-hand side. The fact that the residues balance at the other poles at places congruent to 
$z=\ka+\om, -\ka,2\ka+\om$ modulo $\Lambda$ then follows immediately from the symmetry (\ref{half}), and it is easy to see that the set of finite Schubert parameters 
$(4,\tfrac{2}{3},\tfrac{8}{3})$ for $n=1$, corresponding to the value $z=3\ka$, is a point on the Schubert surface, so this completes the verification that (\ref{schubr}) holds as an 
identity between elliptic functions of $z$, and in particular shows that  (\ref{schubap}) is satisfied for all $n$, and thus by the preceding lemma the second sequence of Schubert 
parameters satisfies (\ref{schubbp}) as well.

Now to prove the first of the equalities in (\ref{aovercab}), we rewrite it as 
\beq\label{1st}
2(X_{a,n}+X_{a,n}^{-1})(M_{b,n}+M_{b,n}^{-1})
=
(P_{a,n}+P_{a,n}^{-1})(P_{b,n}+P_{b,n}^{-1}),
\eeq 
where by Lemma \ref{symlem} we see that the analytic expressions for $M_{b,n},P_{b,n},X_{b,n}$ are obtained by replacing $z\to-z$ in the formulae for $-P_{a,n}^{-1},M_{a,n},X_{a,n}^{-1}$, respectively, so 
from (\ref{maform}) we find
\beq\label{mbform} 
\bear{l}
M_{b,n}=\mu^2\si(z-\ka+\om)\si(z+2\ka+\om)\si(z+2\ka)\si(z-\ka)/\big(\si(z+\ka+\om)\si(z+\om)\si(z)\si(z+\ka)\big)
,\\ 
P_{b,n}=-C_+C_-\tgam\Phi(z+\ka)\Phi(z)^2\Phi(z+2\ka+\om), \\ 
X_{b,n} =-C_+C_-\tgam \Phi(z+\ka)\Phi(z+\om)^2\Phi(z-2\ka+\om)
.
\eear
\eeq 
On each side of the relation (\ref{1st}), there are poles at all the same values of $z$ that were considered in the case of (\ref{schubap}), as well as at points congruent to $z=-2\ka$ and $z=-2\ka-\om$ modulo $\Lambda$.
At $z=0$, there are poles of order 3 on each side, so that we should have $2X_aM_b\sim P_bP_a^{-1}$, so we need to show that 
$$
2\big(X_{a,n}+O(z^2)\big)(M_{b,n}+M_{b,n}^{-1})=\big(P_{b,n}+O(z^2)\big)(P_{a,n}^{-1}+P_{a,n}), 
$$
where the terms $M_{b,n}^{-1}, P_{a,n}$ are corrections of $O(z)$, giving a term with a simple pole when they are multiplied by the triple pole 
in the first bracket on the left/right-hand side, respectively. The expansion around the triple pole is rather arduous, but the problem of showing that the two sides balance can be 
simplified by noting that, from (\ref{schubap}) we have (omitting index $n$) $P_a^{-1}= 2X_a-M_a+O(z)$, which we can use to replace the first term in the second bracket 
on the right 
above, and similarly in the second bracket on the left we can use $M_b=P_b-2X_b^{-1}+O(z)$ from (\ref{schubbp}), so that at leading order, the balancing of the two sides 
is equivalent to $2X_a( P_b-2X_b^{-1})\sim P_b(2X_a-M_a)$, and we can cancel the term $2X_aP_b$ from each side. This may look like the problem has become more difficult, because
we are left with a leading order pole of order 4 on each side, but in fact the singular terms that remain require that, as $z\to 0$, 
\beq\label{singz}
 2X_{a,n}(-2X_{b,n}^{-1}+2M_{b,n}^{-1})\sim P_{b,n}(-M_{a,n}+2P_{a,n}), 
\eeq 
where all omitted terms are $O(1)$, and the corrections  inside each bracket above, namely $2M_{b,n}^{-1}$ and $2P_{a,n}$, are both $O(z)$.
It turns out that, due to (\ref{mbform}), the leading term on the left-hand side of (\ref{singz}) is an even function, namely 
$$
-4X_{a,n}X_{b,n}^{-1}=-4X_{a,n}(z)X_{a,n}(-z)=-4(C_+C_-)^2\Phi(z)^4\hat{F}(z), 
$$
where, using the standard result that any even elliptic function  is given by a rational function of $\wp(z)$ \cite{ww}, we have
$$
\hat{F}(z)=\frac{\Phi(z+2\ka)\Phi(z-2\ka)}{\Phi(z+\ka)\Phi(z-\ka)}
=\left(\frac{\Phi(2\ka)}{\Phi(\ka)}\right)^2\,\frac{\big(\wp(z)-\wp(\ka)\big)\big(\wp(z)-\wp(2\ka+\om)\big)}{\big(\wp(z)-\wp(2\ka)\big)\big(\wp(z)-\wp(\ka+\om)\big)}
=
\frac{1}{2}\big(1+\hat{C}z^2+O(z^4)\big).
$$ 
Then, using $\wp(z)=1/z^2+O(z^2)$ as $z\to 0$ and the addition formula for the Weierstrass $\wp$ function, the coefficient $\hat{C}$ is found to be 
$$ 
\bear{rcl}
\hat{C} &
= & 
\wp(2\ka)-\wp(\ka)+\wp(\ka+\om)-\wp(2\ka+\om) \\
& = & 2\big(\wp(2\ka)-\wp(\ka)\big) +\frac{1}{4}\wp'(\ka)^2/\big(\wp(\ka)-\wp(\om)\big)^2-\frac{1}{4}\wp'(2\ka)^2/\big(\wp(2\ka)-\wp(\om)\big)^2
\\
& = & 2\tmu^{-2}\left(\frac{17}{12}-\frac{29}{12}\right)+
\frac{1}{4}\tmu^2/\left(\tmu^{-4}\left(\frac{29}{12}-\frac{5}{12}\right)^2\right)-\frac{1}{4}\tmu^{-6}/\left(\tmu^{-4}\left(\frac{17}{12}-\frac{5}{12}\right)^2\right)=0,
\eear 
$$
so that the leading order part on the left is 
$$
-4(C_+C_-)^2\Phi(z)^4\hat{F}(z)=-4(C_+C_-)^2\big(\wp(z)-\wp(\om)\big)^2\hat{F}(z)
=-2(C_+C_-)^2\big(1/z^4-2\wp(\om)/z^2+O(1)\big)\big(1+O(z)^4\big), 
$$
while on the right the leading term is another even function, that is 
$$
-M_{a,n}P_{b,n}=-M_{a,n}(z)M_{a,n}(-z)=-(C_+C_-)^2\frac{\Phi(z)^4}{\hat{F}(z)}=-(C_+C_-)^2\frac{\big(\wp(z)-\wp(\om)\big)^2}{\frac{1}{2}\big(1+O(z^4)\big)}, 
$$ 
which gives the same even order singular terms as appear on the left. Thus, for the poles at $z=0$, it remains to check that the residues balance on each side, i.e.\ 
for the remaining correction terms in (\ref{singz}) we must have
$4 X_{a,n}M_{b,n}^{-1}\sim 2P_{b,n}P_{a,n}$, which is a consequence of 
$$\underset{z\to 0}{\mathrm{lim}}M_{b,n}^{-1}(z)/z=\underset{z\to 0}{\mathrm{lim}} P_{a,n}(z)/z \qquad 
\mathrm{and} \qquad 
\underset{z\to 0}{\mathrm{lim}}4X_{a,n}(z)z^2=\underset{z\to 0}{\mathrm{lim}} 2P_{b,n}(z)z^2
;$$ 
the first limit follows from (\ref{mbform}), 
and the second is 
 the identity 
$4C_+C_-\Phi(-2\ka)/\Phi(\ka)=-2C_+C_-\Phi(\ka)/\Phi(2\ka)$. 
Having verified $z=0$, there is an analogous balance of order 3 poles at $z=\om$, which follows immediately by applying the 
symmetry (\ref{half}) to (\ref{1st}). For the balance at $z=\ka$, note that $X_a$ and $P_b$ are both regular there, and we find 
$X_{a,n}(\ka)=2=P_{b,n}(\ka)^{-1}$, so to balance the simple poles on each side of (\ref{1st}) requires 
$2M_b^{-1}\sim P_a^{-1}$, and similar calculations to those done previously show that there is the same residue 
$-\tmu$ on each side of this relation. At $z=-\ka$ there are simple poles in $X_a^{-1},M_b,P_a,P_b$, so their reciprocals have 
simple zeros, and we must verify  $2 X_{a,n}^{-1}M_{b,n}= P_{a,n}P_{b,n}+O(1)$ 
by checking that the double poles balance and 
the residues are the same on each side. Both $\Phi(z+\ka)$ and $1/\si(z+\ka)$ have leading order expansions of 
the form $(z+\ka)^{-1}+O\big((z+\ka)\big)$ in the neighbourhood of $z=-\ka$, so we can write 
$$
2X_{a,n}^{-1}
=\frac{\tilde{F}(z)}{z+\ka}+O\big((z+\ka)\big)=
\tilde{F}(-\ka)\left(\frac{1}{z+\ka}+\frac{d}{dz}\log\tilde{F}(z)|_{z=-\ka}\right)+O\big((z+\ka)\big),
$$  
where $\tilde{F}(z)=-2C_+C_-\tgam^{-2}\Phi(z-2\ka)^{-1}\Phi(z)^{-2}$ is regular as $z\to-\ka$, and write leading order expansions of 
$M_{b,n},P_{a,n},P_{b,n}$ of exactly the same form but with $\tilde{F}$  replaced by appropriate regular functions in each case. 
So at leading order (the coefficient of the double pole) 
we have to verify that the product of the two regular functions on the left equals the product of the 
two regular functions on the right, evaluated at $z=-\ka$, and once this is done the residues are verified by checking that 
the sum of the logarithmic derivatives of the two regular functions on each side takes the same value on each side. 
Now at leading order a short calculation shows that 
$$
\underset{z\to-\ka}{\mathrm{lim}}\frac{2X_{a,n}^{-1}}{P_{b,n}}=\frac{1}{2}=
\underset{z\to-\ka}{\mathrm{lim}}\frac{P_{a,n}}{M_{b,n}},
$$
while by calculating the logarithmic derivative terms on each side, we require that 
$$
\underset{z\to-\ka}{\mathrm{lim}}\frac{d}{dz}\log\left(\frac{2X_{a,n}^{-1}}{P_{b,n}}\right)= 
\underset{z\to-\ka}{\mathrm{lim}}\frac{d}{dz}\log\left(\frac{P_{a,n}}{M_{b,n}} \right), 
$$
leading to a relation involving the Weierstrass zeta function, namely  
$$ 
\frac{\Phi'(3\ka)}{\Phi(3\ka)}+5\frac{\Phi'(\ka)}{\Phi(\ka)}
=
\ze(3\ka-\om)+\ze(3\ka)
-2\ze(\ka-\om)- \ze(\ka+\om)-3\ze(\ka).
$$ 
Then using the fact that $\frac{d}{dz}\log\Phi(z)=\ze(z+\om)-\ze(z)-\ze(\om)$ and 
the quasiperiodicity relation $\ze(z+2\om)=\ze(z)+2\ze(\om)$, this rearranges to yield the identity 
$
\ze(3\ka)+\ze(\ka)+4\ze(\om)-4\ze(\ka+\om)=0$,  
which is verified by rewriting it as 
$$ \bear{rcl}
\big(\ze(3\ka)-\ze(2\ka)-\ze(\ka)\big)
&+ &
\big(\ze(2\ka)-2\ze(\ka)\big)  \, \,- \,\,  4
\big(\ze(\ka+\om)-\ze(\ka)-\ze(\om)\big) \\
& =& \frac{1}{2}\left(\frac{\wp'(2\ka)-\wp'(\ka)}{\wp(2\ka)-\wp(\ka)}\right) 
+
\frac{1}{2}\left(\frac{\wp''(\ka)}{\wp'(\ka)}\right) 
-2
\left(\frac{\wp'(\ka)-\wp'(\om)}{\wp(\ka)-\wp(\om)}\right) \\
&=&
\frac{1}{2}\left(\frac{-\tmu^{-3}-\tmu}{\tmu^{-2}(17/12-29/12)}\right)
+\frac{1}{2}\left(\frac{6\tmu^{-4}(29/12)^2-121/144}{\tmu^{-1}}\right)
-2\left(\frac{\tmu}{\tmu^{-2}(29/12-5/12)}\right)=0, 
\eear 
$$
as required, where we used a standard identity for the zeta function, as well as 
$\wp''(\ka)=6\wp(\ka)^2-\frac{1}{2}g_2$. 
At $z=2\ka$ there are simple poles on each side of (\ref{1st}), coming from the terms $X_a,P_a$, and we have 
to verify $2X_{a,n}(M_{b,n}+M_{b,n}^{-1})\sim P_{a,n}(P_{b,n}+P_{b,n}^{-1})$; but 
$2X_a$ and $P_a$ both have residue $\tmu/2$ at this point, and for the regular part we find the same factor of 
$-13/6$ on each side, as  $z=2\ka$ corresponds to setting  
$n=0$, giving $M_{b,0}=u_3v_3=-3/2$, which is the reciprocal of $P_{b,0}=-f_2^2f_3/f_4=-2/3$, so overall 
the residues are the same. Similarly, at  $z=-2\ka$ there is a balance of simple poles with 
$2M_{b,n}^{-1}(X_{a,n}+X_{a,n}^{-1})\sim P_{b,n}^{-1}(P_{a,n}+P_{a,n}^{-1})$, 
where both $2M_b^{-1}$ and $P_b^{-1}$ have residue $\tmu$, and for the corresponding value $n=-4$, 
we have $X_{a,-4}=3/2=P_{a,-4}^{-1}$, giving the same overall multiplier $13/6$ on each side. 
The balances of poles at the other points congruent to $\pm\ka+\om,\pm 2\ka+\om$ follow 
from the symmetry (\ref{half}), and it is easy to check that (\ref{1st}) is satisfied for $n=1$, so this verifies 
that it holds as an identity of elliptic functions for all $z$, hence in particular is true for all $n\in\Z$; 
the first equality in (\ref{bovercab}) is then given for free, due to Lemma \ref{symlem}. 

Finally, to verify the second equality in  (\ref{aovercab}), we can use Lemma \ref{linlem} to rewrite it in the form 
\beq\label{lid} 
\left(\frac{\bb_n}{\ba_n}+1\right)(P_{b,n}+P_{b,n}^{-1})= \left(\frac{\bs_n}{\ba_n}-1\right)(M_{b,n}+M_{b,n}^{-1}), 
\eeq  
and then we can substitute in 
\beq\label{ratforms}
\frac{\bb_n}{\ba_n}=-f_{n+1}f_{n+2}^2f_{n+3}, \qquad 
\frac{\bs_n}{\ba_n}=\frac{f_{n+3}f_{n+4}^2v_{n+2}^2v_{n+3}^2}{f_{n+2}},
\eeq 
and express it as yet another identity between elliptic functions of $z=(n+2)\ka$. On each side, we find there are poles 
at places congruent to the points $z=0,\om,\pm\ka,\pm\ka+\om, -2\ka,-2\ka+\om$ modulo $\Lambda$. Moreover, under a shift by the 
half-period $\om$ 
it follows from (\ref{half}) and (\ref{mbform}) that $M_{b,n}\to M_{b,n}$ and $P_{b,n}\to -P_{b,n}^{-1}$, while a short calculation using 
(\ref{ratforms}) shows that $\bb_n/\ba_n\to\ba_n/\bb_n$ and $\bs_n/\ba_n\to -\bc_n/\bb_n$, hence overall 
the relation (\ref{lid}) is invariant under this transformation, and therefore it is sufficient to verify the balance of poles at the real values 
$z=0,\pm\ka,-2\ka$, and also check that the identity holds at one other point where it is finite-valued, conveniently chosen 
as $z=2\ka$. So the proof can be completed in the same way as for the other identities, by checking expansions in $z$, but we now
prefer to use a slightly different method which, though formally equivalent, is more arithmetical in nature and easier to apply.
Note that the values of $z$ to be checked correspond to taking $n=-4,-3,-2,-1,0$, and each place where there is a pole 
indicates the presence of $T_0=0$ appearing as a denominator; so at leading order we set $T_0=\eps$, replace all the other 
quantities with their finite values determined by ratios of non-zero terms from the sequences $(S_n)$, $(T_n)$, and consider Laurent expansions in $\eps$,
recovering the appropriate singular behaviour when $\eps\to 0$. In particular, we need to replace 
$f_0=S_0/T_0\to 1/\eps$,  $v_{-1}=T_{-3}T_0/(T_{-2}T_{-1})\to \eps$, 
$v_{0}=T_{-2}T_1/(T_{-1}T_{0})\to -1/\eps$, $v_{1}=T_{-1}T_2/(T_{0}T_{1})\to 1/\eps$, 
$v_{2}=T_{0}T_3/(T_{1}T_{2})\to -\eps$, while all other occurrences of $f_n$, $u_n$ and  $v_n$ that arise, determining  
the terms that appear in the identity (\ref{lid}), correspond to finite and 
non-zero  values, which can be substituted directly. For $n=0$, using the values of $M_{b,0}$ and $P_{b,0}$ as before, we see that both sides of 
the identity are finite and take the value $13/6$. 
When $n=-1$ we find simple poles on each side of (\ref{lid}), with the balance 
$$
-f_0f_1^2f_2(P_{b,-1}+P_{b,-1}^{-1})\sim \left(\frac{f_2f_3^2(v_1v_2)^2}{f_1}-1\right)M_{b,-1}^{-1}
$$ 
in the limit $\eps\to0$,
where we replace $f_0\to\eps^{-1}$, $M_{b,-1}^{-1}=(u_2v_2)^{-1}\to-\frac{1}{2}\eps^{-1}$ and $v_1v_2\to -1$, 
so substituting in $P_{b,-1}=1/2$ and the other finite values of $f_1=1$, $f_2=-1$, $f_3=2$ yields the same residue 
$-5/2$ as the coefficient of $\eps^{-1}$ on each side. The balance for $n=-4$ is similar, with simple poles on each side, 
since $\bs_{-4}/\ba_{-4}=f_{-1}v_{-2}^2(f_0v_{-1})^2/f_{-2}\to - 1$, $M_{b,-4}^{-1}=1/(u_{-1}v_{-1})\to \frac{1}{2}\eps^{-1}$, $P_{b,-4}^{-1}=-f_0/(f_{-1}f_{-2}^2)\to \eps^{-1}$, 
so from $\bb_{-4}/\ba_{-4}=-2$ we find the same residue $-1$ on each side of the balance 
$(\bb_{-4}/\ba_{-4}+1)P_{b,-4}^{-1}\sim (\bs_{-4}/\ba_{-4}-1)M_{b,-4}^{-1}$. 
For the cases of $n=-3,-2$, where each side of (\ref{lid}) has double poles, and poles of order 4, respectively, this simple analysis of leading order 
terms is only sufficient to balance the two leading terms in each case. To balance the poles of lower order (i.e.\ the residues on each side when $n=-3$, and the 
remaining singular terms at order $\eps^{-3},\eps^{-2},\eps^{-1}$ when $n=-2$) it is necessary to calculate higher corrections. If we treat $\eps$ as a local parameter around the 
point $(\infty,1)$ on the curve (\ref{fcurve}), and fix $f_0=\eps^{-1}$, then from the equation of the curve we find the points 
$(f_{\pm1},f_0)$ with expansions 
$$
f_{\pm1}=\pm 1 +\frac{3}{2}\eps\pm\frac{21}{8}\eps^2+3\eps^3+\cdots, 
$$
and subsequently  we obtain 
$$ 
f_{\pm2}=\mp 1 +\frac{1}{2}\eps\mp\frac{13}{8}\eps^2+\frac{1}{2}\eps^3+\cdots,
$$
 either by using the curve or from the map (\ref{fqrt}).
Similarly we can obtain $v_{-1}=\eps-\frac{7}{2}\eps^2+\cdots$, 
$v_0=\eps^{-1}+\frac{5}{2}+\cdots$,
$v_1=-\eps^{-1}+\frac{5}{2}+\cdots$,
$u_0=1+\eps+\cdots$, $u_1=1-\eps+\cdots$,
and further higher corrections to these other terms that appear in (\ref{lid}) when $n=-3,-2$, using the 
equation for the curve (\ref{s5curve}) and/or the map   (\ref{sqrtmap}), and in this way the remaining balances 
of poles on each side of (\ref{lid}) are checked.  
Once again, there is no need to verify  the second equality in  (\ref{bovercab}), because of Lemma \ref{symlem}.
\end{prf}

\begin{figure}
 \centering
\label{quals5fig}
\epsfig{file=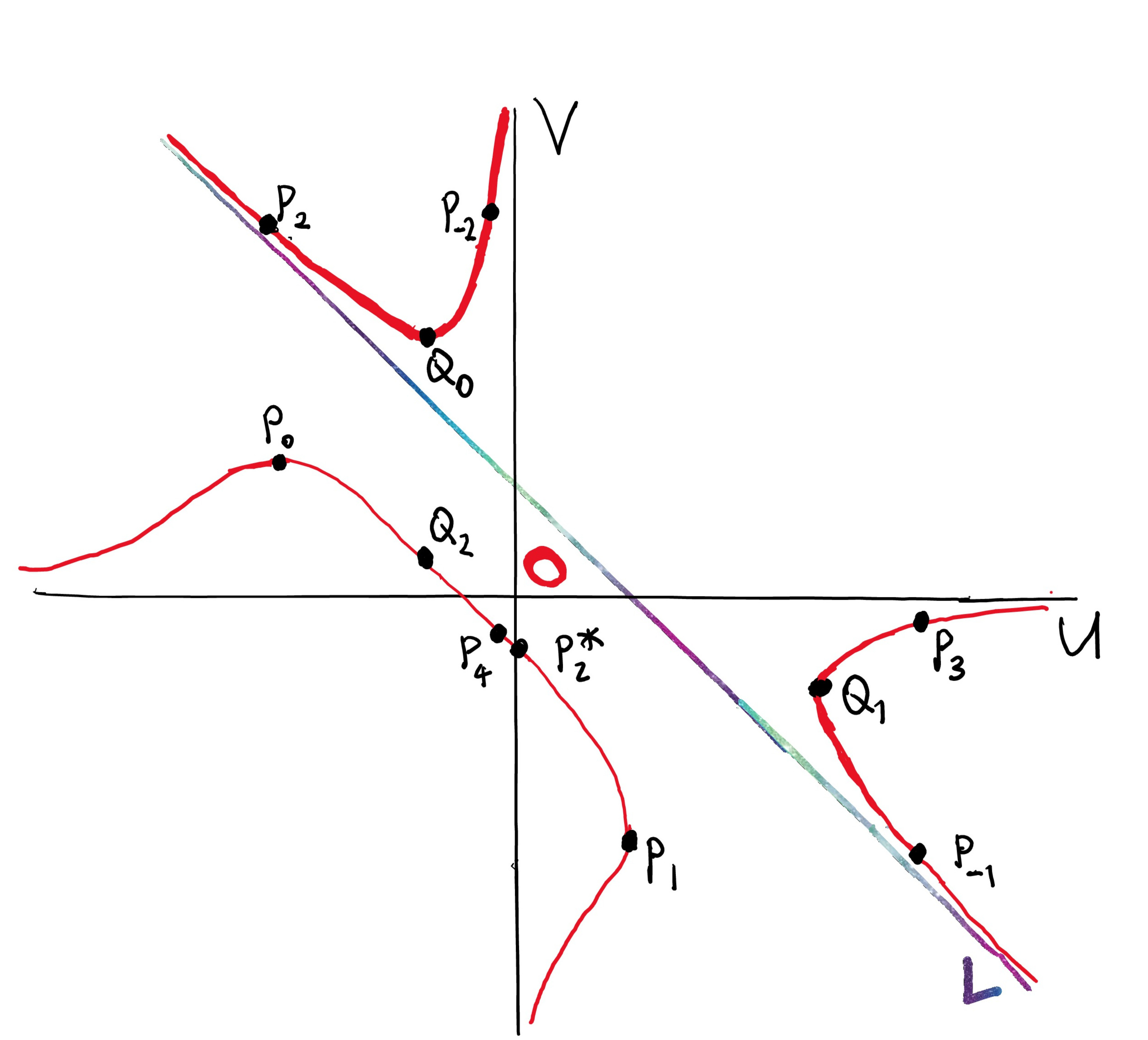, height=3.7in, width=3.7in}
\caption{Qualitative sketch of the Somos-5 curve (\ref{s5curve}).} \end{figure} 
 
In order to understand how the signs of the quantities (\ref{mpxfuva}) and (\ref{mpxfuvb}) change with $n$, we need 
the following result. 

\begin{lem}\label{signlem}
For $n\geq 1$ the signs of the terms in the Somos-5 sequence $(T_n)$, as in (\ref{s5eds}), repeat in the pattern 
\beq\label{signpat} 
+-++-+--+--+-+ 
\eeq 
with period 14. 
\end{lem}
\begin{prf}
The periodic sign pattern for the integer sequence  (\ref{s5eds}), as above,  consists of a block of 7 followed by the same block of 7 but with all the 
signs reversed, and can be deduced by induction 
from the signs 
of the associated sequence of rational numbers $v_n$ defined by the right-hand relation in (\ref{uvdef}), which for 
$n\geq 3$ repeat the pattern 
\beq\label{vsignpat} 
-+-+-+-
\eeq 
with period 7 (cf.\ Table \ref{uvtable}). Indeed, given $T_1=-T_2=T_3=1$, the pattern  (\ref{signpat}) follows from   (\ref{vsignpat}) 
by writing $T_{n+1}=v_nT_nT_{n-1}/T_{n-2}$. So it remains to prove the period 7 sign pattern of the 
rational sequence $(v_n)$, which is achieved by considering the corresponding orbit in $\R^2$ of the QRT map defined by 
(\ref{sqrtmap}), that is 
\beq\label{uvqrt} 
\varphi: \, (U,V)\mapsto \big(V, U^{-1}(1+V^{-1})\big). 
\eeq 
As described above, the orbit lies on the (real) curve (\ref{s5curve}) in the $(U,V)$ plane, shown in Fig.3. 
The real curve 
has four connected components, of which only the part in the positive quadrant is compact, 
but 
this particular orbit lies outside the positive quadrant, being
restricted to the three unbounded components; in contrast, the orbit  
 corresponding  to the sequence $(u_n)$ is given by ratios of the terms of the Somos-5 sequence (\ref{s5seq}), which are all positive, 
and hence lies on the compact oval. 
The relevant properties of the orbit corresponding to $(v_n)$ are not so easy to see from Fig.3, 
which is drawn to scale,
so in Fig.4 
we have produced a more schematic drawing of the real curve (\ref{s5curve}) which highlights the essential features. 
From 
$$
\frac{dV}{dU}=-\frac{2UV+V^2-5V+1}{U^2+2UV-5U+1}, 
$$ 
by taking the resultant of the numerator of the above with the equation of the curve we see that 
the points with horizontal slope have $U$ values given by the four real roots of the quartic $U^4+2U^3-7U^2-2U+7$. In 
particular, there is a local 
maximum at the point ${\cal P}_0$ with coordinates 
$$
\left(-\tfrac{1}{2}\left(\sqrt{2}+1+\sqrt{15+6\sqrt{2}}\right), -\tfrac{1}{2}\left(2\sqrt{2}-5+\sqrt{21-12\sqrt{2}}\right)\right)\approx (-3.63019,0.08211)$$  
and a local minimum at the point ${\cal Q}_0$ with coordinates 
$$
\left(\tfrac{1}{2}\left(\sqrt{2}-1-\sqrt{15-6\sqrt{2}}\right), \tfrac{1}{2}\left(2\sqrt{2}+5+\sqrt{21+12\sqrt{2}}\right) \right)\approx (-1.06909,6.99523); $$ 
the other two stationary points lie on the compact oval.  
Then we consider the orbit of the map (\ref{uvqrt}) starting from the point $\cP^*_3=(v_3,v_4)=(-1,7)$, and it will be convenient to introduce the notation 
$\cP^*_n=(v_n,v_{n+1})$, $\cP_n=\varphi^n(\cP_0)$ and  $\cQ_n=\varphi^n(\cQ_0)$ for $n\in\Z$, as well as letting 
$\overline{\cQ_j\cP_k}$ denote the segment of the real curve  (\ref{s5curve}) connecting the points  $\cQ_j$ and $\cP_k$, and taking 
$\overline{{\cal A}\cP_k}$ to mean the part of the curve with asymptote ${\cal A}$ starting from the point $\cP_k$, where 
we may have ${\cal A}=\cU$ (the $U$ axis), $\cV$ (the $V$ axis) or $\cL$ (the line $U+V=5$). The asymptotes as well as the points 
$\cP_0$, $\cQ_0$ and some of their images/preimages under the QRT map are shown in Fig.4. 
With this notation, for the first 7 iterates we have 
$\cP_3^*\in\overline{\cQ_0\cP_{-2}}$,  
$\cP_4^*\in\overline{\cQ_1\cP_{-1}}$,  
$\cP_5^*\in\overline{\cQ_2\cP_{0}}$,  
$\cP_6^*\in\overline{\cV\cP_{1}}$,  
$\cP_7^*\in\overline{\cL\cP_{2}}$,  
$\cP_8^*\in\overline{\cU\cP_{3}}$,  
$\cP_9^*\in\overline{\cP_2^*\cP_{4}}$, where $\cP_2^*=(0,-1)$. 
Under the action  of the QRT map $\varphi$, which is given by the composition of the two involutions $\iota_h$ and $\iota$ as in (\ref{invols}), 
we have that $\varphi (\overline{\cP_2^*\cP_{4}})\subset \iota (\overline{\cQ_1\cP_{3}})=\overline{\cQ_0\cP_{-2}} $, 
while $\varphi (\overline{\cQ_0\cP_{-2}})=\overline{\cQ_1\cP_{-1}}$ and $\varphi (\overline{\cQ_1\cP_{-1}})=\overline{\cQ_2\cP_{0}}$. 
Also, $\varphi (\overline{\cQ_2\cP_{0}})\subset \overline{\cV\cP_{1}}$,  $\varphi (\overline{\cV\cP_{1}})=\overline{\cL\cP_{2}}$, 
$\varphi (\overline{\cL\cP_{2}})=\overline{\cU\cP_{3}}$, and  $\varphi (\overline{\cU\cP_{3}}) =\overline{\cP_2^*\cP_{4}}$. Thus it follows by induction that 
$\cP_{n+3}^*\in\overline{\cQ_0\cP_{-2}}$,  
$\cP_{n+4}^*\in\overline{\cQ_1\cP_{-1}}$,  
$\cP_{n+5}^*\in\overline{\cQ_2\cP_{0}}$,  
$\cP_{n+6}^*\in\overline{\cV\cP_{1}}$,  
$\cP_{n+7}^*\in\overline{\cL\cP_{2}}$,  
$\cP_{n+8}^*\in\overline{\cU\cP_{3}}$,  
$\cP_{n+9}^*\in\overline{\cP_2^*\cP_{4}}$ holds for all $n\geq 0$. The pattern of signs of the coordinates in these 7 consecutive regions is $(-,+)$, 
$(+,-)$, $(-,+)$, $(+,-)$, $(-,+)$, $(+,-)$, $(-,-)$, which gives the sign pattern (\ref{vsignpat}) for the rational sequence $(v_n)$, as required. 
\end{prf}

Since $S_n>0$ for all $n$, the sign pattern of the sequence $(f_n)$ is clearly the same as (\ref{signpat}), consisting of 
two blocks of 7 that differ by an overall sign, so $\mathrm{sgn}(f_{n+7})=-\mathrm{sgn}(f_{n})$ (where $\mathrm{sgn}$ denotes the sign function). Four of the 
Schubert parameters with signs are given by monomials in the $f_n$ with a homogeneous degree that is even (2 or $-2$), 
hence their sign pattern repeats with period 7, and the other two have a sign that is determined by the sequence $(v_n)$, also 
varying with period 7. As for the signed lengths, their signs are determined by $(T_n)$, so they each have a pattern that varies with period 14, also made up of two blocks of 7 related by a sign flip.   This can be summarized by the following statement.

\begin{cor}\label{signcor} 
The Schubert parameters with signs $(M_{a,n},P_{a,n},X_{a,n})$, $(M_{b,n},P_{b,n},X_{b,n})$ display 5 distinct 
combinations of signs, which for $n\geq 1$ repeat in a sequence with period 7, beginning with 
\beq\label{schubsigns}
\bear{ll}
n=1: \qquad & (+,+,+),(+,+,+) ,\\ 
n=2:            & (+,-,+),(-,+,+) ,\\ 
n=3:            & (-,+,+),(+,+,-) ,\\ 
n=4:            & (+,-,-),(-,-,+) ,\\ 
n=5:            & (+,+,+),(+,+,+) ,\\ 
n=6:            & (+,-,+),(-,+,+) ,\\ 
n=7:            & (-,+,+),(-,-,+) .
\eear 
\eeq 
The signed lengths repeat the following sign patterns for $n\geq 1$:
\beq\label{signleng}
\bear{ll}
\bs_n:\qquad & +--+-++-++-+-- ,\\
\ba_n:           & +++-------++++ ,\\
\bb_n:           & ++-+--+--+-++- ,\\
\bc_n:           & +-+--+--+-++-+ .
\eear
\eeq 
\end{cor}

We are now ready to prove most of the statements in Theorem \ref{main}. For the cases $n\equiv 1$ or $5\,(\bmod\, 7)$, both triples 
of Schubert parameters given by (\ref{mpxfuva}) and (\ref{mpxfuvb}) are positive and satisfy the constraints, so 
they produce a Heron triangle with two rational medians and integer sides $(a,b,c)$ whose ratios are given by  
$$ 
\frac{a}{c}=\frac{\bs_n-\ba_n}{\bs_n-\bc_n}=\frac{\bb_n+\bc_n}{\ba_n+\bb_n}, \quad
\frac{b}{c}=\frac{\bs_n-\bb_n}{\bs_n-\bc_n}=\frac{\bc_n+\ba_n}{\ba_n+\bb_n},
$$ 
and by (\ref{signleng}) the signed lengths are either all positive (when $n\equiv 1,12\, (\bmod\, 14)$) or all negative (when $n\equiv 5,8\, (\bmod\, 14)$), so we can set $a=\pm(\bb_n+\bc_n)$, $b=\pm(\bc_n+\ba_n)$, $c=\pm(\ba_n+\bb_n)$ accordingly, 
which verifies the formulae (\ref{lengthsn}) for the side lengths. Moreover, it follows that the semiperimeter is  
$s=\pm(\ba_n+\bb_n+\bc_n)=\pm\bs_n$, and the reduced lengths are $s-a=\pm\ba_n$, $s-b=\pm\bb_n$, 
$s-c=\pm\bc_n$, so the expression (\ref{arean}) for the area follows immediately from Heron's formula.   
When $n\equiv 2$ or $6\, (\bmod\, 7)$, we need to replace $P_a\to-P_a^{-1}$ and 
$M_b\to-M_b^{-1}$ in order to have two triples of positive coordinates of points on the Schubert surface, 
but this change of signs is compatible with the two constraints, in the sense that it introduces an overall minus 
sign on both sides of the first equality in each of (\ref{aovercab}) and (\ref{bovercab}), so that the side 
ratios of the corresponding Heron triangle with two rational medians are given by 
$$ 
\frac{a}{c}=-\frac{\bs_n-\ba_n}{\bs_n-\bc_n}=-\frac{\bb_n+\bc_n}{\ba_n+\bb_n}, \quad
\frac{b}{c}=-\frac{\bs_n-\bb_n}{\bs_n-\bc_n}=-\frac{\bc_n+\ba_n}{\ba_n+\bb_n}. 
$$ 
For these values of $n$, from (\ref{signleng}) we see that $\ba_n$ and $\bb_n$ are both positive, 
and  $\bs_n$ and $\bc_n$ are both negative, or vice versa, so we can take integer side lengths 
 $a=\mp(\bb_n+\bc_n)$, $b=\mp(\bc_n+\ba_n)$, $c=\pm(\ba_n+\bb_n)$, which verifies  (\ref{lengthsn}), 
while in this case the semiperimeter and reduced lengths are permutations (up to sign) of the signed lengths, as 
(depending on the value of $n\bmod 14$) we have 
$s=\mp\bc_n$, $s-a=\pm\bb_n$, $s-b=\pm\ba_n$, $s-c=\mp\bs_n$, which confirms the area formula (\ref{arean}). 
The analysis of the other three combinations of signs in (\ref{schubsigns}) proceeds similarly. For the case of 
$n\equiv3\, (\bmod\, 7)$ 
only the second constraint  (\ref{bovercab}) acquires an overall minus sign when the negative Schubert parameters 
$M_a,X_b$
are replaced by $-M_a^{-1},-X_b^{-1}$, and we find $s=\mp \bb_n$, $s-a=\pm\bc_n$, $s-b=\mp \bs_n$, $s-c=\pm\ba_n$ 
(according to whether $n\equiv 3$ or $10\, (\bmod\, 14)$). 
When $n\equiv4\, (\bmod\, 7)$, the replacements $P_a\to-P_a^{-1}$, $X_a\to-X_a^{-1}$, $M_b\to-M_b^{-1}$,  $P_b\to-P_b^{-1}$ result in an overall change of sign in the second constraint only, as the signs in (\ref{schubsigns}) 
and (\ref{signleng}) are the exact opposite of those in the previous case, so we have 
$s=\pm \bb_n$, $s-a=\mp\bc_n$, $s-b=\pm \bs_n$, $s-c=\mp\ba_n$. Finally, when $n\equiv 0\,  (\bmod \,7)$, 
replacing $M_a\to-M_a^{-1}$, $M_b\to-M_b^{-1}$,  $P_b\to-P_b^{-1}$ again only introduces a minus sign 
in the second constraint, and we also find $s=\pm \bb_n$, $s-a=\mp\bc_n$, $s-b=\pm \bs_n$, $s-c=\mp\ba_n$ 
in this case. Thus in each case we have a pair of positive triples of Schubert parameters satisfying the necessary constraints. 
If we further require that these should correspond to the half-angle cotangents of appropriate angles in the triangle (cf.\ Fig.1), 
then it may be necessary to apply one or both of the transformations
\beq\label{extrasym}
(M_a,P_a,X_a)\to( M_a^{-1},P_a^{-1},X_a^{-1}), \qquad  
(M_b,P_b,X_b)\to( M_b^{-1},P_b^{-1},X_b^{-1}), 
\eeq 
in order to satisfy the conditions (\ref{extra}).

Having shown that the Schubert parameters with signs in the main sequence for $n\geq 1$ can be consistently transformed 
to a set of positive Schubert parameters, thus providing a sequence of Heron triangles with two rational medians $k,\ell$, 
whose integer sides $(a,b,c)$ and area $\Delta$ are given by the  formulae in Theorem \ref{main}, 
it remains to verify the expressions   (\ref{lengthkln}) for the medians, and also show that $\gcd(a,b,c)=1$, which requires a bit more work. We begin 
by defining signed median lengths, given by 
\beq\label{signkl}
\bk_n=\frac{1}{2}S_{n+4}T_{n+4}(T_nT_{n+1}^2T_{n+2}-S_nS_{n+1}^2S_{n+2}), \qquad 
\bl_n = \frac{1}{2}S_{n}T_{n}(S_{n+2}S_{n+3}^2S_{n+4}-T_{n+2}T_{n+3}^2T_{n+4}),
\eeq
where the overall signs have been chosen so that these initially coincide with the positive median lengths, i.e.\ when 
$n=1$ we have $k=\bk_1=35/2$, $\ell=\bl_1=97/2$.  In order to show that these quantities agree with the median lengths 
up to a sign, we need to prove that $\bk_n$ satisfies a signed version of one of the identities in (\ref{schubb}), and that 
$\bl_n$ satisfies a signed version of one of the analogous identities for $\ell$, obtained by replacing $a\to b$, $b\to c$, $c\to a$ 
and $k\to \ell$ on the right-hand side of each formula. We start by picking the middle identity for $P=P_a$, which becomes 
\beq\label{kid} 
P_{a,n}=\frac{4\Delta_n}{4(\bs_n-\bc_n)\bk_n+(\bs_n-\ba_n)^2-(\bs_n-\bb_n)^2-3(\bs_n-\bc_n)^2},
\eeq  
obtained by setting $P_a\to P_{a,n}$ on the left-hand side, and replacing $a,b,c$ and $k$ by their signed versions 
on the right-hand side, as well as inserting the signed area 
$$ 
\Delta_n=S_nS_{n+1}S_{n+2}^2S_{n+3}S_{n+4}T_nT_{n+1}T_{n+2}^2T_{n+3}T_{n+4}
$$
in the numerator. As for $\ell$, the direct analogue of the first identity in  (\ref{schubb}) is 
$$ 
M_b = \frac{4\Delta}{4c\ell+b^2-3c^2-a^2},
$$ but instead, for reasons that will shortly become clear, we would like to use another expression for $M_b$, namely 
$$ 
M_b = \frac{4c\ell-b^2+3c^2+a^2}{4\Delta},
$$
where the latter is seen to be equivalent to the former due to the relation 
$16\Delta^2=16c^2\ell^2-(b^2-3c^2-a^2)^2$, which follows from Heron's formula and the expression 
for $\ell^2$ in (\ref{kmed}). So as the signed analogue of the latter identity for $M_b$, we take 
\beq\label{ellid}
M_{b,n} = \frac{4(\bs_n-\bc_n)\bl_n-(\bs_n-\bb_n)^2+3(\bs_n-\bc_n)^2+(\bs_n-\ba_n)^2}{4\Delta_n}.
\eeq  

\begin{lem}\label{kllem} The signed median lengths $\bk_n$, $\bl_n$ given by (\ref{signkl}) satisfy the 
relations (\ref{kid}) and (\ref{ellid}) for all $n\in\Z$.
\end{lem}
\begin{prf}
To begin with, we rewrite (\ref{kid}) as 
\beq\label{newkid} 
P_{a,n}^{-1}=\frac{4(\bs_n-\bc_n)\bk_n+(\bs_n-\ba_n)^2-(\bs_n-\bb_n)^2-3(\bs_n-\bc_n)^2}{4\Delta_n},
\eeq 
and consider the symmetry $n\to-n-4$, as in Lemma \ref{symlem}. On the left-hand side 
we have $P_{a,-n-4}^{-1}=-M_{b,n}$, while on the right we use 
the transformations (\ref{signsym}) together with the linear relation (\ref{lin}), 
as well as $\bk_{-n-4}=\bl_n$ and $\Delta_{-n-4}=\Delta_n$, 
to see that (up to an overall minus sign on both sides), the relation (\ref{newkid}) is transformed to (\ref{ellid}). 
Therefore it will be sufficient to prove the above relation involving $\bk_n$ alone, and we can proceed 
as in the proofs of Lemma \ref{linlem} and Theorem \ref{schubrels}, by substituting in the analytic 
formulae for each of the terms and regarding it as an identity between elliptic functions of $z=(n+2)\ka$, 
so that the left-hand side has simple poles at the places $z\equiv 0,\om,\ka,\ka+\om\bmod\Lambda$. 
However, we also wish to exploit the additional symmetry under shifting by the half-period $\om$, as 
we have $P_{a,n}\to P_{a,n}$ under this symmetry, and using the  results of our previous calculations 
we see that overall the 
right-hand side of (\ref{newkid}) is also left invariant by this transformation. Thus it is sufficient to check only the poles 
at $z=0,\ka$ as well as one other value where the relation is finite, and the case $n=1$ (corresponding to $z=3\ka$) 
where $P_{a,1}=2/3$ is readily 
verified. Then since we only have residues at two simple poles to check, corresponding to the values $n=-2,-1$, we can use 
the simplified method with a local parameter $\eps$, as at the end of the proof of Theorem \ref{schubrels}. Using 
(\ref{mpxfuva}) and 
(\ref{lin}) we can rewrite 
(\ref{newkid}) as 
\beq\label{nnewkid} 
-u_{n+2}v_{n+2}= \frac{4\left(1+\frac{\bb_n}{\ba_n}\right)\frac{\bk_n}{\ba_n}
+\left(\frac{\bs_n}{\ba_n}-1\right)^2-\left(1+\frac{\bc_n}{\ba_n}\right)^2-3\left(1+\frac{\bb_n}{\ba_n}\right)^2}{4{\Delta_n}/{\ba_n^2}}.
\eeq 
On the right-hand side above we can make use of the expressions (\ref{ratforms}), as well as 
$$\bc_n/\ba_n=-f_n^2f_{n+1}v_{n+2}^2v_{n+3}^2/f_{n+2}$$ and 
$$
\frac{\bk_n}{\ba_n}=-\frac{f_{n+4}v_{n+2}v_{n+3}^2(1-f_nf_{n+1}^2f_{n+2})}{2f_{n+2}}, 
\qquad 
\frac{\Delta_n}{\ba_n^2}=f_nf_{n+1}f_{n+3}f_{n+4}v_{n+2}^2v_{n+3}^2.
$$
Now when $n=-2$, at leading order on the left-hand side we have $-u_0v_0\sim \eps^{-1}$, 
while on the right-hand side the leading order in the denominator is 
$4\Delta_{-2}/\ba_{-2}^2\sim 4\eps^{-4}$. In the numerator, $\bb_{-2}/\ba_{-2}\sim\eps^{2}$ and 
$\bk_{-2}/\ba_{-2}\sim \frac{1}{2}\eps^{-3}$, hence the first term gives the leading order contribution 
$4(1+\bb_{-2}/\ba_{-2})\bk_{-2}/\ba_{-2}\sim 2\eps^{-5}$, but for the difference of squares that follows we need a 
correction at next-to-leading order, so that from $\bs_{-2}/\ba_{-2}=\eps^{-3}+\frac{1}{2}\eps^{-2}+\cdots$, 
$\bc_{-2}/\ba_{-2}=\eps^{-3}-\frac{1}{2}\eps^{-2}+\cdots$ we find 
$(\bs_{-2}/\ba_{-2}-1)^2-(\bc_{-2}/\ba_{-2}+1)^2\sim 2\eps^{-5}$, and the final term $-3(1+\bb_{-2}/\ba_{-2})^2$ 
gives a lower order contribution at $O(\eps^{-4})$; thus overall at leading order the right-hand side gives 
$(2\eps^{-5}+2\eps^{-5})/(4\eps^{-4})=\eps^{-1}$, as required.  
When $n=-2$, the left-hand side is $-u_1v_1\sim-\eps^{-1}$, while the denominator of the 
right-hand side is $\Delta_{-1}/\ba_{-1}^2\sim 8\eps^{-1}$, and the numerator contains the  terms
$\bb_{-1}/\ba_{-1}\sim\eps^{-1}$, $\bk_{-1}/\ba_{-1}\sim-\eps^{-1}$, 
$\bs_{-1}/\ba_{-1}\sim-4$, $\bc_{-1}/\ba_{-1}\sim-\eps^{-1}$, 
so that overall these combine to give 
$(-4\eps^{-2}-\eps^{-2}-3\eps^{-2})/(8\eps^{-1})=-\eps^{-1}$ at leading order, in agreement with the left-hand side. 
Note also that for the values $n=-4,-3$ and $0$, the left-hand side is finite while the other side has removable singularities: some of the terms in the numerator/denominator on the  right-hand  side of (\ref{nnewkid}) are singular, but overall these cancel to 
give a finite value in the limit $\eps\to 0$; so the relation holds as an  identity between elliptic functions of $z=(n+2)\ka$, and in particular for all $n\in\Z$.
\end{prf} 

Given that $\bk_n$ satisfies (\ref{kid}), for each $n$ we can compare this with the corresponding positive Schubert 
parameter $P_a>0$, given by 
$$ 
P_a=\frac{4\Delta}{4ck+a^2-b^2-3c^2}=\frac{4ck-a^2+b^2+3c^2}{4\Delta},
$$ 
where we use either the first or the second rational expression above involving $a,b,c,k$ and $\Delta$, 
according to whether $P_{a,n}=\pm P_a^{\pm 1}$, determined by $n\bmod 7$ as in Corollary \ref{signcor}. Now 
$\mathrm{sgn}(\Delta_n)=\mathrm{sgn}(T_nT_{n+1}T_{n+3}T_{n+4})$, which repeats the pattern 
$$ 
\Delta_n: \qquad ++--+++
$$ with period 7, 
while our previous analysis also showed that, for each $n$, $\bs_n-\bc_n=\mathrm{sgn}(\bs_n-\bc_n)\, c$, where the sign 
pattern is 
$$ 
 \bs_n-\bc_n: \qquad ++-+--+--+-++-
$$ 
with period 14. For the cases $n\equiv 1,2,5,6,0 \, (\bmod \,7)$ we see that $\Delta_n=\Delta$, and then according to whether 
$P_{a,n}>0$ (when $n\equiv 1,5,0$) or $P_{a,n}<0$ (when $n\equiv 2,6$) we can directly compare the right-hand side of (\ref{kid}) with 
either the first formula for $P_a$ above, or compare $-P_{a,n}^{-1}$ with the second formula above, respectively; then since the squared terms 
can always be identified, that is 
$a^2=(\bs_n-\ba_n)^2$ etc., we have $ck=(\bs_n-\bc_n)\bk_n$ in the first case, giving $\bk_n=\mathrm{sgn}(\bs_n-\bc_n)\, k$, 
but $ck=-(\bs_n-\bc_n)\bk_n$ in the second case, giving $\bk_n=-\mathrm{sgn}(\bs_n-\bc_n)\, k$. However, the cases 
$n\equiv 3,4 \, (\bmod \,7)$, when $\Delta_n=-\Delta$, are different because in those cases we need to  apply the first of the transformations 
in (\ref{extrasym}) to ensure that (\ref{extra}) holds. So for  $n\equiv 3 \, (\bmod \,7)$ that means comparing $P_{a,n}^{-1}$ with the second formula for 
$P_a$ above, yielding  $\bk_n=-\mathrm{sgn}(\bs_n-\bc_n)\, k$, and for  $n\equiv 4 \, (\bmod \,7)$ it requires that $-P_{a,n}$ should be compared with the first formula for 
$P_a$, hence $\bk_n=\mathrm{sgn}(\bs_n-\bc_n)\, k$. Thus overall we see that $\bk_n$ is related to the median length $k$ by an overall sign, which varies with period 14 in 
the following pattern: 
\beq\label{kpattern} 
\bk_n: \qquad +-++-++-+--+-- .
\eeq
Similarly, $\bl_n$ is equal to the median length $\ell$ up to a sign with pattern 
\beq\label{lpattern} 
\bl_n: \qquad +-++-+--+--+-+ , 
\eeq
and this verifies that the formulae (\ref{lengthkln})  hold. 

Based on computer experiments, it was observed  in \cite{br1} that for the triangles in the main sequence,  the pairs $(\theta,\phi)$ corresponding to 
the parametrization (\ref{abcparam}), given by (\ref{thetaphi}) with $+$ signs in both equations, cycle through one of five isomorphic plane curves 
in a pattern that repeats with period 7. Applying the symmetry $a\leftrightarrow b$, $k\leftrightarrow\ell$, one obtains an alternative pair of parameters 
\beq\label{altthetaphi}
\tilde{\theta} = \frac{c-b+2k}{2s}, 
\qquad 
\tilde{\phi}   = \frac{a-c+2\ell}{2s},
\eeq 
and Buchholz and Rathbun noted that these pairs cycle with period 7 through the same set of curves but in a different order, namely ${\cal C}_1,{\cal C}_2,{\cal C}_3,{\cal C}_4,{\cal C}_1,{\cal C}_2,{\cal C}_5$. 
A more detailed study of the allowed discrete symmetries in \cite{br3} showed that by applying appropriate permutations of $a,b,c,k,\ell$ and changes of sign, one could also obtain pairs of coordinates 
on three more (isomorphic) curves ${\cal C}_6,{\cal C}_7,{\cal C}_8$. However, until now there was no explanation for the period 7 behaviour with respect to the index $n$ in the main sequence, which 
we provide here. 

\begin{thm}\label{7cycle} 
For $n\geq 1$, the rational parameters 
$$ 
\theta = \frac{c-a+2\ell}{a+b+c}, \qquad 
\phi   = \frac{b-c+2k}{a+b+c}
$$ 
corresponding to the main 
sequence of Heron triangles with sides $(a,b,c)$ and medians $k,\ell$ lie on one of five birationally equivalent plane curves of genus one, which repeat in a pattern with period 7. 
The same is true for any such sequence of parameters obtained from these by the action of the discrete symmetry group that leaves the equations 
$4k^2=2b^2+2c^2-a^2$, $4\ell^2=2c^2+2a^2-b^2$ invariant.   
\end{thm}
\begin{prf}
We begin by considering the pair of quantities  
$$ 
\theta_n = \frac{\ba_n-\bc_n+2\bl_n}{2\bs_n}, \qquad 
\phi_n   = \frac{\bc_n-\bb_n+2\bk_n}{2\bs_n},
$$ 
obtained by replacing $a,b,c,k,\ell$ by their signed counterparts $\bs_n-\ba_n,\bs_n-\bb_n,\ldots$ etc.\ For $n=1$ this coincides with $\theta,\phi$. Moreover, from the analytic 
parametrization of the two Somos-5 sequences, this gives a pair of independent elliptic functions of $z=(n+2)\ka$. Then by a standard result in the theory of elliptic functions, 
the pairs $(\theta_n,\phi_n)$ lie on a plane curve ${\cal C}$ of genus one, which is birationally equivalent to the original elliptic curve
(this follows from the Riemann-Roch theorem, but there is also an explicit  classical method for constructing the equation of ${\cal C}$, based on the expressions for 
$\theta_n,\phi_n$ in terms of $\wp,\wp'$ \cite{akhiezer}). Under the change $n\to n+1$, the original pair $(\theta,\phi)$ corresponds to a different pair of elliptic 
functions, because each term $a,b,\ldots$ that appears is related to $\bs_n-\ba_n,\bs_n-\bb_n,\ldots$ by a sign which changes with $n$; so the pair of 
coordinates for the next value of $n$ will lie on a different curve, ${\cal C}^*$ say. Now all of these signs repeat with period 14, but in each block of 7 the signs are flipped 
with respect to the previous block, and only the ratio of terms appears in $\theta,\phi$, so their overall pattern of signs repeats with period 7, hence there are at most 
7 different curves repeating with this period. However, a closer examination of (\ref{signleng}), (\ref{kpattern}) and (\ref{lpattern}) reveals that there are really only 5 
different pairs of functions appearing in each block of 7, because all the signs for $n\equiv1\, (\bmod\,7)$ and  $n\equiv5\, (\bmod\,7)$ are the opposite of each other, and similarly for 
 $n\equiv2\, (\bmod\,7)$ and  $n\equiv6\, (\bmod\,7)$. The same argument applies to any pair $(\theta,\phi)$ obtained from this one by the action of the discrete symmetry group.
\end{prf}

In order to show that the side lengths given by  (\ref{lengthsn}) have greatest common divisor 1 for all $n\geq 1$, 
it is necessary to reconsider the orbits of the QRT map $\varphi$ defined by (\ref{uvqrt}), associated with the rational sequences $(u_n)$, $(v_n)$, and examine their 
reduction modulo a prime $p$, which corresponds to considering the map over the finite field $\F_p$. (For a detailed treatment of QRT maps over finite fields, see \cite{jogia}.) A general QRT 
map (over any field) is a birational map $\varphi:\,\Proj^1\times\Proj^1\to\Proj^1\times\Proj^1$ given by a composition of two involutions, and as already mentioned, each 
orbit lies on a curve belonging to a pencil of biquadratic curves in the plane; generic curves ${\cal C}$ have genus one, and the orbit gives  a sequence of points $\hat{{\cal P}}_0+n{\cal P}\in{\cal C}$, with 
each iteration of $\varphi$ corresponding to addition of the point ${\cal P}$ in the group law of the curve. For the particular case at hand, the curve is given by (\ref{s5curve}), which 
has identity element ${\cal O}=(\infty,\infty)$, and $\iota$ in (\ref{invols}) is the elliptic involution that sends any point $\hat{\cal P}\to-\hat{\cal P}$. 
Under the QRT map $\varphi=\iota\circ\iota_h$,  ${\cal P}=(\infty, 0)$ is the point being added at each step, and the orbit corresponding to the Somos-5 sequence (\ref{s5seq})  
via (\ref{uvdef}) is 
\beq\label{orbit1} 
\varphi^n(u_0,u_1)=(u_n,u_{n+1})=(1,1)+n{\cal P},
\eeq  
where the initial point $\hat{{\cal P}}_0=(u_0,u_1)=(1,1)\in{\cal C}$ is 2-torsion, 
while the orbit associated with the other sequence (\ref{s5eds}) is 
\beq\label{orbit2} 
\varphi^n(v_0,v_1)=(v_n,v_{n+1})=n{\cal P},
\eeq
with initial point ${\cal O}=(v_0,v_1)=(\infty,\infty)$. 

The above description is valid over $\C$ or any subfield where the orbit is defined, in particular over $\Q$, 
in which case ${\cal P}$ is an infinite order element 
in the Mordell-Weil group of the curve ${\cal C}(\Q)$, which is generated by ${\cal P}$ and the 2-torsion point $\hat{{\cal P}}_0=(1,1)$. Under reduction $\bmod \,p$ the same description holds provided 
that $p$ is a prime of good reduction, so that the curve ${\cal C}(\F_p)$ is non-singular, with the main difference being that now the curve has finitely 
many points satisfying the Hasse-Weil bound $\left| p+1-\#{\cal C}(\F_p)\right|\leq2\sqrt{p}$, so 
 ${\cal P}$ has some finite order which we denote by $\mathrm{ord}_p({\cal P})$. Hence it follows that both orbits (\ref{orbit1}) and (\ref{orbit2}) over $\F_p$ are periodic with 
the same period    $\mathrm{ord}_p({\cal P})$. 

For primes of bad reduction, which in this case are 2,3 and 17, the situation is slightly more complicated. For the primes $p=2,3$, we find that 
the singular curve ${\cal C}(\F_p)$ is reducible, since $-5\bmod 2=-5\bmod 3=1$, so in both cases  we can factorize (\ref{s5curve}) as 
\beq\label{sing23}
(UV+1)(U+V+1)=0.
\eeq
When $p=2$, the two orbits  $(u_n,u_{n+1})$ and  $(v_n,v_{n+1})$  are the same up to a shift of starting point, repeating the sequence 
$$ 
(\infty,\infty),(\infty,0),(0,1),(1,1),(1,0),(0,\infty)
$$ 
with period $6=2\times(2+1)$, corresponding to the fact that each (genus zero) irreducible component of (\ref{sing23}) is isomorphic to  
$\Proj^1(\F_2)$, so contains 3 points, and the orbit alternates between points on each component, i.e.\ $(\infty,0),(1,1),(0,\infty)\in\{\, UV+1=0\,\}$ 
and $(\infty,\infty),(0,1),(1,0)\in\{\, U+V+1=0\,\}$. Similarly, when $p=3$ both orbits have period $8=2\times (3+1)$, repeating the sequence 
$$ 
(\infty,\infty),(\infty,0),(0,2),(2,1),(1,1),(1,2),(2,0),(0,\infty),
$$ 
which again alternates between the two irreducible components of the curve (\ref{sing23}). For $p=17$, the curve (\ref{s5curve}) 
 is irreducible and singular: this is a case of non-split multiplicative reduction \cite{lmfdb}, so it has genus zero and contains $\#{\cal C}(\F_{17})=p+2=19$ points, and after removing the singular point 
$(5,5)$, ${\cal C}(\F_{17})\setminus\{\, (5,5)\,\}$ has the structure of an abelian group of order 18. The two orbits both have period 9,  corresponding to the two cosets of the cyclic subgroup of order 9 generated 
by ${\cal P}$; the coordinates of the points in these two orbits can be read off from Table \ref{uvtable17}, which also includes part of the associated periodic sequence $(f_n\bmod \,17)$ (of period 18).  

Although it is not essential to our main argument, we can now offer a brief explanation for how the periods modulo a prime arise for the corresponding Somos-5 sequences. Setting $t= \mathrm{ord}_p({\cal P})$ and 
working $\bmod\,p$, from the relation $u_{n+t}=u_n$ we have 
$$ 
 \frac{S_{n-2+t}S_{n+1+t}}{S_{n-1+t}S_{n+t}}= \frac{S_{n-2}S_{n+1}}{S_{n-1}S_{n}}
$$  
for all $n$, and clearly the period of the Somos-5 sequence must be a multiple of $t$. 
(We will write everything for the sequence $(S_n)$, but the same treatment applies to $(T_n)$ or any other Somos-5 sequence.) 
Somos-5 sequences have a 3-parameter group of gauge transformations which leave the ratios $u_n$, and hence the recurrence (\ref{s5orig}), unchanged: one can replace 
$S_n\to A^*_{\pm}(B^*)^n \,S_n$, i.e. rescale all even/odd index terms by an arbitrary non-zero scalar  $A^*_{\pm}$, and rescale each term by the powers of another arbitrary quantity $B^*\neq 0$; this 
symmetry group has a  natural interpretation in terms of the quiver that defines the associated cluster algebra \cite{fh}. Hence we see that $S_n$ and $S_{n+t}$ should be related by a gauge transformation 
of this kind.  
\begin{propn}
Under shifting by $t= \mathrm{per}(u_n\bmod \,p)$, the period of the associated QRT map, 
the Somos-5 sequence reduced $\bmod \,p$ satisfies
\beq\label{quasi}
S_{n+t}=A^*_\pm(B^*)^n\,S_n, \qquad 
A^*_+, A^*_-B^*, (B^*)^2\in\F_p^*
\eeq
for even/odd $n$, respectively, where $A^*_+,A^*_-,B^*$ are constants independent of $n$.
\end{propn} 
\begin{prf}
To prove this gauge transformation formula directly, define the ratio $r_n=S_{n+t}/S_n$, and observe that $u_{n+t}=u_n$ for all $n$ implies that 
$r_{n+2}/r_n=r_{n+3}/r_{n+1}=r_{n+4}/r_{n+2}=r_{n+5}/r_{n+3}=(B^*)^2$ say, a constant in 
$\F_p^*$ independent of $n$, 
and this can be solved separately for even/odd $n$ to give $r_n=A^*_{\pm}(B^*)^n$, 
where $A^*_+=r_0$, $A^*_-B^*=r_1\in\F_p^*$ (assuming none of the terms in these ratios is zero).
\end{prf} 
\begin{cor} Let $t$ be 
the order of the point ${\cal P}\in {\cal C}(\F_{p})$ 
if $p$ is a prime of good reduction (or equivalently, the period of 
$(u_n\bmod\,p)$, appropriately reinterpreted in the case of bad reduction).
Then the  period  
of the Somos-5 sequence $\bmod \,p$ is a multiple of $t$, 
given by 
\beq\label{per}
\mathrm{per}(S_n\bmod \,p)=\ell^*\, t
\leq2(p-1)t,
\eeq
where (up to a possible factor of 2) $\ell^*$ 
denotes the lowest common multiple of $\mathrm{ord}_p\big((B^*)^2\big)$ and/or
the orders in  $\F_p^*$ 
of one or two other combinations of $A^*_+,A^*_-,B^*$, 
depending on the parity of $t$.
\end{cor}
\begin{prf} 
Iterating the  quasiperiodicity relation  (\ref{quasi}) $j$ times  when $t$ is even  yields 
$$ 
S_{n+jt}=(A^*_\pm)^j(B^*)^{nj+tj(j-1)/2}\,S_n,
$$ 
for even/odd $n$, so 
requiring 
$S_{n+jt}=S_n$ for all $n$ imposes the conditions 
\beq\label{leven}
(B^*)^{2j}=1, \qquad (A_+^*)^j=(A_-^*B^*)^j=(B^*)^{-tj(j-1)/2}=\pm 1, 
\eeq 
where the choice of $j$ satisfying the first condition for a given $B^*$ fixes the sign on the right. Then $\ell^*$ is the smallest 
positive value of $j$ that satisfies all these conditions. An odd value of $\ell^*$ is 
only possible in the case that $(B^*)^2$, $A_+^*$ and $A_-^*B^*$ are all 
quadratic non-residues in $\F_p^*$, in which case only the plus sign can occur on the right-hand side above, and 
$$ 
\ell^*= \mathrm{lcm} \big(\mathrm{ord}_p(A^*_+),\mathrm{ord}_p(A^*_-B^*),\mathrm{ord}_p((B^*)^2)\big), 
$$
but otherwise  there must be an even $\ell^*|\mathrm{lcm} \big(\mathrm{ord}_p(A^*_+),\mathrm{ord}_p(A^*_-B^*),\mathrm{ord}_p((B^*)^2)\big)$, where it may be possible to divide out some factors of 2 from the first two terms, depending 
on the orders of these elements.
For odd $t$, iterating (\ref{quasi}) $j$ times produces  
$$ 
S_{n+jt}=(A^*_\pm)^{\lc\frac{j+1}{2}\rc}(A^*_\mp)^{\lc\frac{j}{2}\rc}(B^*)^{nj+tj(j-1)/2}\,S_n,
$$ 
and then by requiring $S_{n+jt}=S_n$ for all $n$, in the case of odd $j$ (exploiting the freedom to replace 
$B^*\to -B^*$ and $A^*_-\to -A^*_-$ simultaneously) this leads to 
$$
A^*_+=A^*_-, \qquad (B^*)^j=(A^*_+)^j=1\implies \ell^* =  
\mathrm{lcm} \big(\mathrm{ord}_p(A^*_+),\mathrm{ord}_p(B^*)\big),
$$
corresponding to a situation where $\tilde{K}_0=\tilde{K}_1$ and all the terms satisfy the same Somos-4 recurrence 
(cf.\ equation (\ref{2invt}) and \cite{swartvdp}).  
For an even value of $j$ with $t$ odd, we find instead 
$$(B^*)^j=1, \qquad (A^*_+A^*_-B^*)^{j/2}=(B^*)^{\big(1-t(j-1)\big)j/2}=1,$$
so now, given the values of $(B^*)^2$ and $A^*_+A^*_-B^*$, the smallest possible value 
of $j/2$ can be found, yielding 
\beq\label{lodd}
\ell^*=2\, \mathrm{lcm} \big(\mathrm{ord}_p(A^*_+A^*_-B^*),\mathrm{ord}_p((B^*)^2)\big).
\eeq 
The upper bound on the period comes from 
Fermat's little theorem. 
\end{prf}
An advantage of the relation (\ref{quasi}) is that it allows the period of the Somos-5 sequence to be computed from (\ref{per}) without calculating so many terms: in general it is 
sufficient to find the minimum $t$ such that $(u_0,u_1)=(u_t,u_{t+1})$, giving $t= \mathrm{ord}_p({\cal P})$, which requires the $t+5$ adjacent terms $S_{-2},\ldots,S_{t+2}$, and 
then $A^*_\pm,B^*$ can be obtained from the ratios $S_t/S_0$, $S_{t+1}/S_1$ and $S_{t+2}/S_2$ (provided $S_0S_1S_2\neq0$, otherwise  one of the ratios
$S_{t-1}/S_{-1}$, $S_{t-2}/S_{-2}$  can be used instead). For example, in the case of the sequence (\ref{s5seq}) taken $\bmod\,23$, we find that $(u_n \bmod \,23)$ repeats the 
pattern $1,1,2,13,20,3,20,13,2$ with period 9, so  $t= \mathrm{ord}_{23}({\cal P})=9$ 
which is odd, while $A^*_+\equiv S_9/S_0 \bmod\,23=S_{10}/S_1\bmod\,23=21$ and $S_{11}/S_2\bmod\,23=20$, 
thus we find $(B^*)^2=20/21=13=6^2$ in $\F_{23}$, and taking $B^*=6$ gives $A^*_-=21/6=15\neq A^*_+$; 
then $\mathrm{ord}_{23}(13)=11$ and $\mathrm{ord}_{23}(21\times15\times 6)=\mathrm{ord}_{23}(21^2)=11$, so $\ell^*=22$ 
by (\ref{lodd}), 
hence $\mathrm{per}(S_n\bmod \,23)=9\times 22=198$, in agreement with the value found by Robinson in \cite{rob}. As another example, for the sequence (\ref{s5eds}) taken $\bmod\,61$, 
the corresponding terms  $(v_n \bmod \,61)$ have a repeating pattern 
$$\infty,\infty,0,60, 7,25,12,34,11,32,2,1,1,2,32,11,34,12,25,7,60,0$$ 
with even period $t= 22$, and from the terms $1,1$ in the middle it is apparent that this is the same as the orbit $(u_n \bmod \,61)$ but shifted; thus from $T_{21}/T_{-1}\bmod\,61=14\equiv A_-^*(B^*)^{-1}$, 
$T_{23}/T_{1}\bmod\,61=60\equiv A_-^*B^*$, $T_{24}/T_{2}\bmod\,61=14\equiv A_+^*(B^*)^2$ we find $A_-^*=13=14^2$, $A_+^*=48=14^{-1}$, $B^*=14$ in $\F_{61}$, and 
$\mathrm{ord}_{61}(14)=6$, so $\ell^*$ divides 
$\mathrm{lcm} \big(\mathrm{ord}_{61}(14^{-1}),\mathrm{ord}_{61}(14^3),\mathrm{ord}_{61}(14^2)\big)$=6, 
the smallest value of $j$ satisfying all the requirements in (\ref{leven}), hence  
 the period of the sequence $(T_n\bmod \,61)$ is $22\times 6=132$, and this is the same as the period of the original Somos-5 sequence 
$(S_n\bmod \,61)$ as found in \cite{rob}. 

We are now ready to finish off the proof of the main theorem.
\begin{lem}\label{gcdstn}
For all $n\in\Z$, the terms of the sequences (\ref{s5seq}) and (\ref{s5eds}) with the same index are coprime, that is 
$$ 
\gcd(S_n,T_n)=1.
$$ 
\end{lem}
\begin{prf} Suppose that for some $n$ the terms $S_n$ and $T_n$ have a common divisor. Then for the reduced sequences in $\F_p$ it follows 
that $S_n\equiv 0\equiv T_n$, and so $u_n=u_{n+1}=\infty=v_n=v_{n+1}$. If $p$ is a prime of good reduction for the curve ${\cal C}$ given by 
(\ref{s5curve}), then in the group law of ${\cal C}(\F_p)$ this implies that 
$(1,1)+n{\cal P}={\cal O}=n{\cal P}$, which is a contradiction. If $p$ is a prime of bad reduction,  then we have $p=2,3$ or 17, and $p=2$ or $p=3$  are 
impossible by Proposition \ref{period23}, while Table \ref{uvtable17} shows that 17 is never a divisor of $S_n$, since the sequence $(u_n\bmod\,17)$ has period 9 and 
remains finite. 
\end{prf}

\begin{table}[h!]
  \begin{center}
    \caption{The analogue of Table \ref{uvtable} in the finite field $\F_{17}$.}
    \label{uvtable17}
    \begin{tabular}{ | r|| r| r| r| r | r| r| r| r| r| r|} %
\hline
      $n$ & 0  &  1 & 2 & 3 & 4 & 5 & 6 & 7 &8 & 9 \\
\hline 
\hline   
&&&&&&&&&& \\ 
$u_n$ &  $1$ & $1$ &$2$    & $10$ & $15$ & $6$ & $15$ & $10$ & $2$ & $1$ \\
&&&&&&&&&& \\ 
    $v_n$ & $\infty$ &  $\infty$ & $0$ & $16$ & $7$ &    $11$ &  $7$ & $16$ &  $0$ & $\infty$\\
&&&&&&&&&& \\ 
      $f_n$ & $\infty$ & $1$ & $16$  & $2$ & $3$ 
& $9$ & $12$ & $14$ & $6$ & $\infty$
 \\
&&&&&&&&&& \\ 
\hline 
    \end{tabular}
  \end{center}
\end{table}

\begin{propn}\label{abcngcd}
For all $n\in\Z$, the signed lengths satisfy $\gcd(\ba_n,\bb_n,\bc_n)=1$. 
\end{propn}
\begin{prf} First of all, note that from Lemma \ref{linlem} we have 
$\gcd(\bs_n,\ba_n,\bb_n,\bc_n)=\gcd(\ba_n,\bb_n,\bc_n)$. Up to an overall sign, each of the signed lengths 
is a degree six monomial in the two sets of five adjacent values $S_n,\ldots,S_{n+4}$ and $T_n,\ldots,T_{n+4}$. 
 Now suppose that a prime $p$ is a divisor of $\bc_n$.  If $p|S_n$, then by Lemma \ref{copr}
it cannot divide $S_{n+1},S_{n+2}$ or $S_{n+3}$, so if it divides $\bb_n$ then it must divide $T_{n+2}$, 
thus (by Lemma \ref{copr} again) it is coprime to $S_{n+3},S_{n+4},T_n$ and $T_{n+1}$, hence it cannot divide $\bs_n$.
Similarly, if $p|T_{n+4}$ and it is a divisor of $\ba_n$ then, by the same lemma, it is coprime to $T_{n+1},T_{n+2}$ and 
$T_{n+3}$, so it must divide $S_{n+2}$, but in that case again it is coprime to $\bs_n$. Applying the same lemma once more, we see that if 
$p|S_{n+1}$ then $p|\bb_n$, so if $p|\bs_n$ as well then it must divide either $T_n$ or $T_{n+1}$, but the first case  is impossible because 
$p$ must also divide $\ba_n$ by (\ref{lin}), yet it has to be coprime to $S_{n+2},T_{n+1},T_{n+2},T_{n+3}$; so this leaves the second case, 
which requires  
$p|\gcd(S_{n+1},T_{n+1})$, contradicting Lemma  \ref{gcdstn}. Finally, if $p|T_{n+3}$ then an analogous argument leads to 
$p|\gcd(S_{n+3},T_{n+3})$, another contradiction. 
\end{prf}

The formulae (\ref{lengthsn}) imply that 
$\gcd(a,b,c)=\gcd(\bb_n+\bc_n,\bc_n+\ba_n,\ba_n+\bb_n)$, and any prime divisor of the latter three linear combinations 
of the signed lengths must be a divisor of $2\ba_n=(\bc_n+\ba_n)+(\ba_n+\bb_n)-(\bb_n+\bc_n)$, and similarly of $2\bb_n$ 
and $2\bc_n$, so by the above proposition the only possible common divisor of $(a,b,c)$ is 2, but Proposition  \ref{period23} 
shows that 
precisely one of $\bb_n+\bc_n,\bc_n+\ba_n,\ba_n+\bb_n$ is even (cf.\ Table \ref{lengtable}), so 
$\gcd(a,b,c)=1$
and this completes the proof of Theorem \ref{main}. 

\section{Brahmagupta angles and geometrical identities} 
\setcounter{equation}{0}

By virtue of the fact that there is a 2-isogeny relating the curves (\ref{s5curve}) and (\ref{fcurve}), there are infinitely many identities between 
elements of the associated function fields, which correspond to relations between terms of the sequences $(f_n)$ and $(v_n)$, or equivalently $(u_n)$, 
since $u_n=v_nf_{n-2}f_{n+1}/(f_{n-1}f_n)$.   Here we present some identities that arise naturally from the geometry of Heron triangles, which in 
particular  leads 
to an appealing way to represent and visualize the triangles in the main sequence. 

We begin by considering Brahmagupta's construction of rational Heron triangles, which is based on concatenating two rational Pythagorean right 
triangles with a common height of length $2r$, where $a,b$ are the hypotenuses of the two  triangles. 
So if $d,e$ are the bases of the two triangles, then we have 
\beq\label{bracon}
4r^2+d^2=a^2, \qquad 4r^2+e^2=b^2, \qquad \pm(d\pm e)=c,
\eeq
where the choice of sign inside the bracket above depends on whether the two triangles are joined back-to-back along their height, or are overlapping, and in the latter 
case there is an overall sign outside depending on which of $d$ or $e$ is the larger. 

The case of two overlapping right triangles with $d>e$ is the situation relevant to the triangle for $n=1$ in the main sequence, and we denote the acute angles 
at the base of each triangle by $\psi_a,\psi_b$, respectively. These two angles completely determine the Heron triangle, up to scale, and we will refer to them as the Brahmagupta angles. 
Then upon comparing with Fig.1 it is clear that 
$$ 
\psi_a=\pi-\be-\gam, \qquad \psi_b=\pi-\gam'+\al', 
$$ 
with 
$P_a=\mathrm{cot}(\be/2)$, $X_a=\mathrm{cot}(\gam/2)$ as before, and 
$M_b=\mathrm{cot}(\al'/2)$, $X_b=\mathrm{cot}(\gam'/2)$. So from $2r=a\sin\psi_a=b\sin\psi_b$ and trigonometric
identities we can write 
\beq\label{pqform} 
r=a\, \frac{(P_a+X_a)(P_aX_a-1)}{(P_a^2+1)(X_a^2+1)}
=b\, \frac{(M_b-X_b)(M_bX_b+1)}{(M_b^2+1)(X_b^2+1)},
\eeq
while we have $c=d-e$ with $d=a\cos\psi_a$, $e=b\cos\psi_b$, and then we reproduce Brahmagupta's formulae (\ref{brahma}) by taking 
$$ 
p=a\cos^2\left(\frac{\beta+\gam}{2}\right), \qquad q=b\sin^2\left(\frac{\gam'-\al'}{2}\right).
$$
Once again we can rewrite the trigonometric functions in terms of the appropriate Schubert parameters, to obtain 
\beq\label{rform} 
\frac{p}{a}=\frac{(P_aX_a-1)^2}{(P_a^2+1)(X_a^2+1)}, \qquad 
\frac{q}{b}= \frac{(M_b-X_b)^2}{(M_b^2+1)(X_b^2+1)}
\eeq
The above formulae have been written with the appropriate sign choices for the $(73,51,26)$ triangle, which has Brahmagupta parameters 
$p=49/13$, $q=588/13$, $r=210/13$. However,  each of the lengths and 
Schubert parameters in the two different expressions for $r$ in (\ref{rform}) can now be replaced by their signed counterparts, immediately yielding another relation, which is 
equivalent to an identity between elliptic functions.
%
%
\begin{thm}\label{brahmathm} For all $n\in\Z$ 
the identity 
$$ 
(\bs_n-\ba_n)\sin\psi_{a,n}=(\bs_n-\bb_n)\sin\psi_{b,n} 
$$ 
holds for the signed Brahmagupta angles $\psi_{a,n}$, $\psi_{b,n}$, which are defined in terms of the Schubert 
parameters with signs by 
\beq\label{sbrahma}
\sin\psi_{a,n}=\frac{(P_{a,n}+X_{a,n})(P_{a,n}X_{a,n}-1)}{(P_{a,n}^2+1)(X_{a,n}^2+1)}, \qquad 
\sin\psi_{b,n}= \frac{(M_{b,n}-X_{b,n})(M_{b,n}X_{b,n}+1)}{(M_{b,n}^2+1)(X_{b,n}^2+1)}.
\eeq 
\end{thm}

\begin{figure}
\centering
\begin{subfigure}{.5\textwidth}
  \centering
  \includegraphics[width=.7\linewidth]{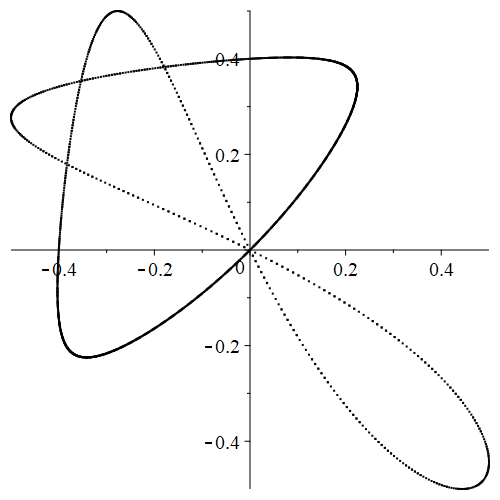}
  \caption{Plot of $(\sin\psi_{a,n},\sin\psi_{b,n})$.}
  \label{brahmacurve} 
\end{subfigure}%
\begin{subfigure}{.5\textwidth}
  \centering
  \includegraphics[width=.7\linewidth]{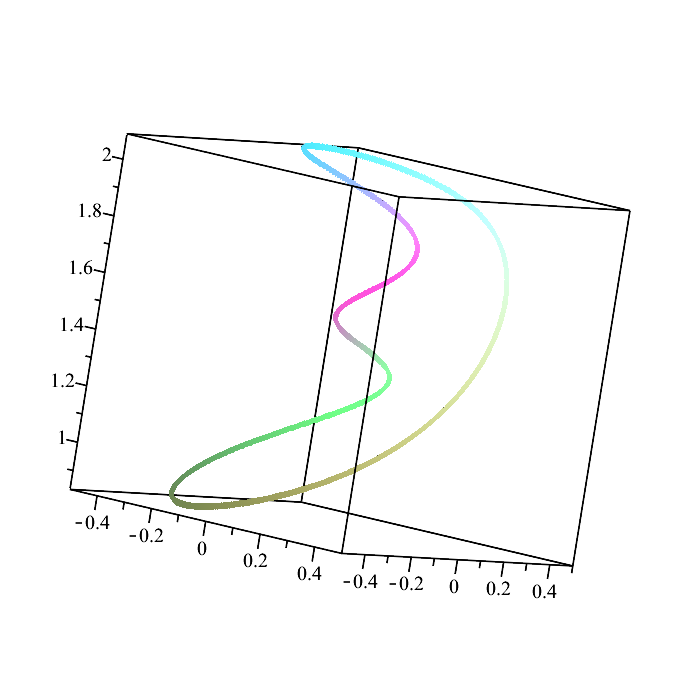}
  \caption{3D plot of $(\sin\psi_{a,n},\sin\psi_{b,n},u_n)$. }
  \label{brahmacurve3D} 
\end{subfigure}
\caption{500 points on the orbit of the Brahmagupta angles for the main sequence.}
\label{fig:test}
\end{figure}

The advantage of using these angles (which are allowed to be negative) instead of the actual positive acute angles is that their sines 
are given by fixed rational functions of the Schubert parameters, hence they define the same functions on the curve  for all $n$. 
By the same argument as in the proof of Theorem \ref{7cycle}, this means that the pairs $(\sin\psi_{a,n},\sin\psi_{b,n})$ lie on 
an algebraic  curve that is birationally equivalent to ${\cal C}$.   
The sequence of pairs of sines of (signed) Brahmagupta angles for $n=1,\ldots, 500$ is plotted in Fig.\ref{brahmacurve}, which shows 
that this plane curve has self-intersections, while in  Fig.\ref{brahmacurve3D} this is lifted to 3D by plotting the 
corresponding value of $u_n$ as a third component, showing a space curve  
without self-intersections.

As one more example of using elementary geometry to derive identities between elliptic functions, we start from 
a formula for the tangent of the bisected angle opposite side $a$, which is given on the first page of Schubert's monograph \cite{schubert}, 
and in our notation reads 
$$ 
\tan\left(\frac{\al+\be}{2}\right) =\frac{(s-b)(s-c)}{\Delta}.
$$

\begin{thm}\label{bis}
The squares of the coordinates of the points on the orbit of the QRT map $\varphi$ associated wth the sequence (\ref{s5eds}) 
are given as rational functions of five adjacent quantities $f_n$ by 
\beq\label{vsq} 
v_{n+2}^2=\frac{f_{n+1}^2f_{n+2}^3+f_{n+4}}{f_{n+3}(f_n^2-f_{n+2}f_{n+4})}, 
\qquad 
v_{n+3}^2=\frac{f_{n+2}^3f_{n+3}^2-f_{n}}{f_{n+1}(f_{n+4}^2+f_{n}f_{n+2})}.
\eeq 
\end{thm}
\begin{prf} 
In terms of the Schubert parameters, the tangent of the bisected angle for the $n$th triangle is 
$(P_{a,n}+M_{a,n})/(P_{a,n}M_{a,n}-1)$, while the right-hand side can be written in terms of ratios 
of $S_n$ and $T_n$ terms, and ultimately as $f_nf_{n+1}f_{n+2}/f_{n+4}$, so that this becomes an identity for 
signed quantities. Then substituting for the Schubert parameters from  (\ref{mpxfuva}) and writing everything in terms 
of $f_n$ and $v_n$, this produces the formula 
$$ 
\frac{1}{f_{n+4}}=\frac{f_{n+2}f_{n+3}\,v_{n+2}^2+1}{f_n^2f_{n+3}\,v_{n+2}^2-f_{n+1}^2f_{n+2}^3};
$$ 
the right-hand side above is a M\"obius transformation of $v_{n+2}^2$, which is inverted to obtain the 
first expression in (\ref{vsq}). An analogous calculation starting from the corresponding formula for 
the tangent of the bisected angle opposite side $b$ produces the second expression for $v_{n+3}^2$.  
\end{prf} 

\begin{remark}
In principle one can use repeated application of the recurrence (\ref{fqrt}) in (\ref{vsq})  to eliminate 
$f_{n+4},f_{n+3},f_{n+2}$ and 
 write $v_{n+2},v_{n+3}$ 
as rational functions of $f_n,f_{n+1}$ only, then use  (\ref{fcurve}) to simplify the resulting formulae;  
 this would be a rather tedious way to verify the 2-isogeny given by (\ref{2isog}) in terms of 
these pairs of coordinates. Since it involves the product $(v_{n+2}v_{n+3})^2$, substituting for the 
latter with (\ref{vsq}) and applying the same method provides an alternative proof of the identity  (\ref{vfid}), 
which is equivalent to the linear relation (\ref{lin}) for the signed lengths, although from the point of view of the logical progression of the paper this 
is a  somewhat circular argument, since Lemma \ref{linlem} was one of the key steps in showing that 
the quantities $f_n$ describe the geometry of the main sequence of Heron triangles as claimed. 
\end{remark} 

\section{Conclusions} 
\setcounter{equation}{0}

\begin{table}[h!]
  \begin{center}
    \caption{Prime factors of the Schubert parameters in the sporadic cases.}
    \label{tab:sporadicmpx}
\scalebox{1.1}{
    \begin{tabular}{ | r|| r| r| r| r | r| r|} %
\hline
      $n$ & $M_a$  &  $P_a$ & $X_a$ & $M_b$ & $P_b$ & $X_b$ \\
\hline 
&&&&&& \\ 
    *
& $\tfrac{2^3\cdot 7\cdot 13}{3\cdot 17}$ &  $17$ & $\tfrac{2^4\cdot 3}{7\cdot 13} $
 &  $\tfrac{3\cdot 7\cdot 11}{2^2\cdot 5\cdot 13}$ & $\tfrac{11\cdot 13\cdot 17}{2^2\cdot 3\cdot 5\cdot 7}$ &  $\tfrac{17}{5\cdot 11}$  \\
&&&&&& \\ 
      **  
& $\tfrac{3^2\cdot 5\cdot 31}{2^2\cdot 7\cdot 17}$ & $\tfrac{2^2\cdot 5\cdot 31}{3^2\cdot17}$ & $\tfrac{3^2\cdot 7}{5\cdot17}  $
& $\tfrac{3\cdot 7\cdot 17}{5\cdot 19}$ & $\tfrac{3\cdot 5\cdot 17\cdot 19}{2^3\cdot 7\cdot 31}$ & $\tfrac{3\cdot 19\cdot 31}{2^4\cdot 5\cdot 17}$ 
 \\
&&&&&& \\ 
	*** & $\tfrac{2^3\cdot 19\cdot 47}{3^2\cdot 11\cdot 23}$ & $\tfrac{3^2\cdot 11\cdot 17 \cdot 47}{2^3\cdot 7 \cdot 19 \cdot 23}$ & $\tfrac{2\cdot 7 \cdot 11 \cdot 47}{17\cdot 19\cdot 23}$ & 
$\tfrac{2^4\cdot 3^2\cdot7\cdot 17\cdot 23}{5\cdot 11\cdot 19 \cdot 97}$ &  $\tfrac{2\cdot7\cdot 11\cdot 17\cdot 19}{5\cdot 23\cdot 97}$ & $\tfrac{5\cdot 11 \cdot 19 \cdot 23}{2^5\cdot3^2\cdot 97}$ \\
&&&&&& \\ 
**** & $\tfrac{2^2\cdot 3\cdot 17 \cdot 43\cdot 59}{5\cdot 11\cdot 13\cdot 19\cdot 23} $  
& $\tfrac{2^3\cdot 3\cdot 5\cdot 13\cdot 23\cdot43}{11\cdot 19\cdot 41 \cdot 59} $  
& $\tfrac{19\cdot 23\cdot 41 \cdot 43}{5 \cdot 11 \cdot 13 \cdot 17 \cdot 59} $ 
 & $\tfrac{19\cdot 23 \cdot 41 \cdot 43}{2^4\cdot 3 \cdot 7 \cdot 11 \cdot 59} $  
& $\tfrac{11 \cdot 19 \cdot 41 \cdot 59}{7\cdot 17 \cdot 23 \cdot 43} $  
& $\tfrac{7 \cdot 19 \cdot 43 \cdot 59}{2^5\cdot 3 \cdot 11 \cdot 17 \cdot 23} $ \\
&&&&&& \\ 
\hline 
    \end{tabular}
}
  \end{center}
\end{table}

\begin{table}[h!]
  \begin{center}
    \caption{Prime factors of the semiperimeter, reduced side lengths and area in the sporadic cases.}
    \label{tab:sporadicsabcbar}
\scalebox{0.8}{
    \begin{tabular}{ | r|| r| r| r| r | r |} %
\hline
      $n$ & $s$  &  $s-a$ & $s-b$ & $s-c$ & $\Delta$ \\
\hline 
    * &  $3\cdot 7\cdot 13\cdot 17 $& $2^3\cdot 5^2\cdot 17$ & $3\cdot 7\cdot 13$ &  $2^3\cdot 11^2$ & $2^3\cdot3\cdot 5\cdot 7\cdot 11\cdot 13\cdot 17$  \\ 
      ** & $2^3\cdot 7\cdot 19^2$ &  $2^3\cdot 3^6$  & $5^2\cdot 7\cdot 31$ 
&  $17^2\cdot 31$
& $2^3\cdot 3^3\cdot 5\cdot 7\cdot 17\cdot 19\cdot 31$ 
 \\
	*** & $2^3\cdot 3^2\cdot 11^2\cdot 19^2$ & $2^5\cdot 3^2\cdot 5^2\cdot 7\cdot 17$& $23^2\cdot 47^2$ & $ 7\cdot 17\cdot 97^2$ & $2^4\cdot3^2\cdot 5\cdot7\cdot 11\cdot 17 \cdot 19\cdot 23\cdot 47\cdot 97$
  \\
**** &   $17\cdot 23^2\cdot59^2$    
& $5^2\cdot 7^2\cdot 13^2\cdot 41$ 
& $2^4\cdot 3\cdot11^2\cdot 43^2$ 
& $2^4\cdot 3\cdot 17\cdot 19^2\cdot 41$
& $2^4\cdot 3\cdot 5\cdot 7\cdot 11\cdot 13\cdot 17\cdot 19 \cdot 23 \cdot 41\cdot 43 \cdot 59$
\\
\hline 
    \end{tabular}
}
  \end{center}
\end{table}

We have proved all the empirical observations on the infinite sequence of Heron triangles with two rational medians 
found by Buchholz and Rathbun in \cite{br1}, and this led us to explicit formulae for the side lengths, rational medians 
and the area 
in terms of the Somos-5 sequences (\ref{s5seq}) and (\ref{s5eds}).  
The crux of our proof was to view the underlying elliptic curve as a complex torus, and construct an analytic embedding of 
this torus in the Schubert surface, such that the required algebraic identities for the two sets of Schubert parameters 
are valid for all complex values of the argument, and in particular at an infinite discrete set of points corresponding to 
the indices $n\in\Z$. The indices $n\geq 1$ provide infinitely many distinct triangles, while under the 
involution $n\to -n-4$ the indices $n\leq 0$ correspond to the same set of triangles repeated with reversed orientation via the symmetry 
$a\leftrightarrow b$, together with a finite number of singular values. 

It would also be possible to prove Theorem \ref{schubrels} by purely algebraic means, using (\ref{fqrt}) and the 2-isogeny 
(\ref{2isog}), but this would require extensive amounts of computer algebra, without providing much insight into the problem.
Our hope was that the analytical approach would give a better understanding of the structure of the triangles in the main 
sequence, and suggest whether this might allow the sporadic  solutions to be extended to other infinite families, 
perhaps by writing them in terms of different Somos sequences, as well as possibly 
shedding some light on the harder problem of showing why no perfect triangle exists.  
However, the detailed features of the solutions 
in the main sequence appear to rely on specific arithmetical  properties of the group of rational points on the curve 
${\cal C}$, and even the particular numerical values of the coordinates of the two generators, and we have been unable to 
generalize this to produce one or more other infinite families. Nevertheless, the prime factors of the Schubert parameters 
and reduced lengths for the sporadic triangles (see Tables \ref{tab:sporadicmpx} and \ref{tab:sporadicsabcbar}) give 
tantalizing hints that there might be other families that these solutions could belong to, only with a different structure 
compared with the main sequence. A more extensive computer search, and the discovery of more new solutions,  
would provide further evidence in this direction, or otherwise suggest that these four solutions are truly sporadic.

\noindent \textbf{Acknowledgments:} This research was supported by 
Fellowship EP/M004333/1  from the Engineering \& Physical Sciences Research Council, UK, 
and grant 
 IEC\textbackslash R3\textbackslash 193024 from the Royal Society.
It is a pleasure to thank Chris Athorne and Claire Gilson for the opportunity to present a preliminary version of 
this work at the joint 
British Mathematical/Applied Mathematics Colloquium in Glasgow, and to thank 
Harry Braden for some helpful comments during the ICMS 
Integrable Systems Virtual Seminar in April 2021. 
All data generated or analysed during this study are included in this published article.

\end{document}